\documentclass[11pt,a4paper,usenames,dvipsnames]{article}

\usepackage{url,amsmath,amssymb,latexsym,pstricks,mathrsfs,comment,amsthm,graphicx,tikz,tikz-cd,enumerate,accents,pgffor,cite,wrapfig,multicol,float,calc,geometry,tikz-3dplot}
\usepackage[colorlinks]{hyperref}
\usepackage[subnum]{cases}

\usepackage[margin=10pt,font=small,labelfont=bf, labelsep=period]{caption}

\newcommand{\nc}{\newcommand}
\nc{\rnc}{\renewcommand}

\nc\U{{\raise1.6 ex\hbox{\rotatebox{180}{$U$}}}}

\usepackage{xcolor}

\usepackage{fnpct}
\setfnpct{dont-mess-around}

\usepackage{enumitem}

\usepackage[T1]{fontenc}
\usepackage{ wasysym }
\usepackage{stmaryrd}

\usepackage{algcompatible}
\usepackage{algorithmicx}
\usepackage{algorithm}
\usepackage[noend]{algpseudocode}
\algrenewcommand{\algorithmiccomment}[1]{\hfill[{\it #1}]}

\usetikzlibrary{positioning, decorations}
\makeatletter
\tikzset{
  distance from start/.code={%
    \pgfgetpath\currentpath\pgfprocessround{\currentpath}{\currentpath}%
    \pgf@decorate@parsesoftpath{\currentpath}{\currentpath}%
    \pgfmathparse{#1/\pgf@decorate@totalpathlength}\tikzset{pos=\pgfmathresult}},
  distance from end/.code={%
    \pgfgetpath\currentpath\pgfprocessround{\currentpath}{\currentpath}%
    \pgf@decorate@parsesoftpath{\currentpath}{\currentpath}%
    \pgfmathparse{1-(#1/\pgf@decorate@totalpathlength)}\tikzset{pos=\pgfmathresult}}
}

\usetikzlibrary{calc}

\allowdisplaybreaks

\makeatletter

\DeclareMathSymbol{\widetildesym}{\mathord}{largesymbols}{"65}
\newcommand\lowerwidetildesym{%
  \text{\smash{\raisebox{-1.3ex}{%
    $\widetildesym$}}}}
\newcommand\fixwidetilde[1]{%
  \mathchoice
    {\accentset{\displaystyle\lowerwidetildesym}{#1}}
    {\accentset{\textstyle\lowerwidetildesym}{#1}}
    {\accentset{\scriptstyle\lowerwidetildesym}{#1}}
    {\accentset{\scriptscriptstyle\lowerwidetildesym}{#1}}
}

\nc\tup[6]{\big[(#1,#2),(#3,#4),(#5,#6)\big]}
\nc\tupp[6]{(#1,#2,#3,#4,#5,#6)}
\nc\fix{\operatorname{fix}}
\nc\Fix{\operatorname{Fix}}
\nc\Aut{\operatorname{Aut}}

\nc\tn{t_n}
\nc\td{{}^dt}
\nc\tnd{{}^dt_n}
\nc\Tn{\T_n}
\nc\Td{{}^d\T}
\nc\Tnd{{}^d\T_n}
\nc\mun{\mu_n}
\nc\mud{{}^d\!\mu}
\nc\rhon{\rho_n}
\nc\rhod{{}^d\!\rho}

\nc\graph[7]{
\begin{tikzpicture}[scale=#7]
\node[circle, inner sep=2pt, fill=blue!20] (1) at (90:4) {\footnotesize $1$};
\node[circle, inner sep=2pt, fill=blue!20] (2) at (210:4) {\footnotesize $2$};
\node[circle, inner sep=2pt, fill=blue!20] (3) at (330:4) {\footnotesize $3$};
\node[circle, inner sep=2pt, fill=blue!20] (4) at (0,0) {\footnotesize $4$};
\draw 
(2)  --node[circle,inner sep=2pt, fill=white]{\footnotesize $#1$}  (3)
(3)  --node[circle,inner sep=2pt, fill=white]{\footnotesize $#2$}  (1)
(1)  --node[circle,inner sep=2pt, fill=white]{\footnotesize $#3$}  (2)
(1)  --node[circle,inner sep=2pt, fill=white]{\footnotesize $#4$}  (4)
(2)  --node[circle,inner sep=2pt, fill=white]{\footnotesize $#5$}  (4)
(3)  --node[circle,inner sep=2pt, fill=white]{\footnotesize $#6$}  (4)
;
\end{tikzpicture}
}

\nc\graphhack[8]{
\begin{tikzpicture}[scale=#7]
\node[circle, inner sep=2pt, fill=blue!20] (1) at (90:4) {\footnotesize $1$};
\node[circle, inner sep=2pt, fill=blue!20] (2) at (330:4) {\footnotesize $2$};
\node[circle, inner sep=2pt, fill=blue!20] (3) at (210:4) {\footnotesize $3$};
\node[circle, inner sep=2pt, fill=blue!20] (4) at (0,0) {\footnotesize $4$};
\draw 
(2)  --node[circle,inner sep=2pt, fill=white]{\footnotesize $#1$}  (3)
(3)  --node[circle,inner sep=2pt, fill=white]{\footnotesize $#2$}  (1)
(1)  --node[circle,inner sep=2pt, fill=white]{\footnotesize $#3$}  (2)
(1)  --node[circle,inner sep=2pt, fill=white]{\footnotesize $#4$}  (4)
(2)  --node[circle,inner sep=2pt, fill=white]{\footnotesize $#5$}  (4)
(3)  --node[circle,inner sep=2pt, fill=white]{\footnotesize $#6$}  (4)
;
\node[circle, outer sep=+0pt, fill=white] () at (270:4) {};
\node () at (270:3.2) {$#8$};
\end{tikzpicture}
}

\rnc\L{\mathcal L}
\nc\R{\mathcal R}
\nc\bu{{\bf u}}
\nc\bv{{\bf v}}
\nc\bx{{\bf x}}
\nc\by{{\bf y}}
\nc\bz{{\bf z}}
\nc\bl{{\bf l}}
\nc\br{{\bf r}}
\nc\nt{\frac n2}
\nc\tnt{\tfrac n2}
\nc\fnt{\lfloor\frac n2\rfloor}
\nc\tfnt{\lfloor\tfrac n2\rfloor}
\nc\cnt{\lceil\frac n2\rceil}
\nc\s{\mathfrak s}
\rnc\t{\mathfrak t}

\nc\Rect{\mathscr R}

\nc\udotted[2]{\draw[dotted](#1+.4,2)--(#2-.4,2);}
\nc\ddotted[2]{\draw[dotted](#1+.4,0)--(#2-.4,0);}

\nc{\ubluebox}[2]{\bluebox{#1}{1.7}{#2}2\udotted{#1}{#2}}
\nc{\lbluebox}[2]{\bluebox{#1}0{#2}{.3}\ddotted{#1}{#2}}
\nc{\ublueboxes}[1]{{
\foreach \x/\y in {#1}
{ \ubluebox{\x}{\y}}}
}
\nc{\lblueboxes}[1]{{
\foreach \x/\y in {#1}
{ \lbluebox{\x}{\y}}}
}
\nc{\udottedsm}[2]{\draw [dotted] (#1+.4,2)--(#2-.4,2);}
\nc{\udottedsms}[1]{{
\foreach \x/\y in {#1}
{ \udottedsm{\x}{\y}}
}}
\nc{\ldottedsm}[2]{\draw [dotted] (#1+.4,0)--(#2-.4,0);}
\nc{\ldottedsms}[1]{{
\foreach \x/\y in {#1}
{ \ldottedsm{\x}{\y}}
}}
\nc{\bluebox}[4]{
\draw[color=blue!20, fill=blue!20] (#1,#2)--(#3,#2)--(#3,#4)--(#1,#4)--(#1,#2);
}
\nc{\uudotted}[2]{\draw [dotted] (#1+.4,4)--(#2-.4,4);}
\nc{\uudotteds}[1]{{
\foreach \x/\y in {#1}
{ \uudotted{\x}{\y}}
}}
\nc{\uubluebox}[2]{\bluebox{#1}{3.7}{#2}4\uudotted{#1}{#2}}
\nc{\uublueboxes}[1]{{
\foreach \x/\y in {#1}
{ \uubluebox{\x}{\y}}}
}
\nc{\uudottedsm}[2]{\draw [dotted] (#1+.4,4)--(#2-.4,4);}
\nc{\uudottedsms}[1]{{
\foreach \x/\y in {#1}
{ \uudottedsm{\x}{\y}}
}}

\nc\sib{\overline\si}
\nc\veb{\overline\ve}
\nc\taub{\overline\tau}
\nc\lamb{\overline\lam}
\nc\rhob{\overline\rho}
\nc\iob{\overline\io}
\nc\Sib{\overline\Si}
\nc\up+

\nc\es\emptyset
\nc\op\oplus
\nc\io\iota
\nc\kk\Bbbk

\nc\OrdP[1]{\mathscr P^+(#1)}
\nc\OrdPP[1]{\mathscr P^{++}(#1)} 
\nc\OrdR[1]{\mathscr P^-(#1)}
\nc\OrdPR[1]{\mathscr P^\pm(#1)}
\nc\OriP[1]{\mathscr A^+(#1)}
\nc\OriR[1]{\mathscr A^-(#1)}
\nc\OriPR[1]{\mathscr A^\pm(#1)}

\nc{\Sone}{S^1}
\nc\comp{\sim}
\nc\wt\fixwidetilde
\nc\nab{\wt\nabla}
\nc\itemit[1]{\item[\emph{(#1)}]}
\nc\ve\varepsilon
\nc\Lam\Lambda
\nc\IB{\mathcal{IB}}
\nc\ab{\overline{\al}}
\nc\bb{\overline{\be}}
\nc\Th\Theta
\nc\rev{\operatorname{rev}}
\nc\RR{\mathbb R}
\nc\E{\mathcal E}
\nc\G{\mathcal G}
\nc\Q{\mathcal Q}
\nc\PT{\mathcal{PT}}
\nc\A{\mathcal A}

\nc\LLL{\mathbb L}
\nc\RRR{\mathbb R}
\nc\bn{{\bf n}}
\nc\btwo{{\bf 2}}
\nc\bthree{{\bf 3}}
\nc\Mod[1]{\ (\operatorname{mod} #1)}
\nc\even{{\operatorname{even}}}
\nc\odd{{\operatorname{odd}}}
\nc\TL{\mathcal T\!\mathcal L}
\nc\TLpm{\TL^\pm}
\nc\TLm{\TL^-}
\nc\Jpm{\JJ^\pm}

\nc\LL{\mathcal L}
\nc\Fx{\F^\times}
\nc\PL{P\LL}
\nc\Ker{\operatorname{Ker}}
\nc\vu{{\bf u}}
\nc\vv{{\bf v}}
\nc\vw{{\bf w}}
\nc\vzero{{\bf 0}}
\nc\GL{\operatorname{GL}}
\nc\SL{\operatorname{SL}}
\nc\PGL{\operatorname{PGL}}
\nc\tr{{\operatorname{T}}}
\nc\BB{\mathscr B}
\nc{\mat}[4]{\left[\begin{matrix}#1&#2\\#3&#4\end{matrix}\right]}
\nc{\tmat}[4]{\left[\begin{smallmatrix}#1&#2\\#3&#4\end{smallmatrix}\right]}
\nc\sicongnode[6] {\node[rounded corners,rectangle,draw,fill=blue!20] (#1#2) at (#3,#4) {$\si_{#5,#6}$};}
\nc\redsicongnode[6] {\node[rounded corners,rectangle,draw,fill=red!20] (#1#2) at (#3,#4) {$\si_{#5,#6}$};}
\nc\taucongnode[4] {\node[rounded corners,rectangle,draw,fill=red!20] (#1) at (#2,#3) {$\tau_{#4}$};}
\nc\congnode[4] {\node[rounded corners,rectangle,draw,fill=blue!20] (#3) at (#1,#2) {$#4$};}
\nc\redcongnode[4] {\node[rounded corners,rectangle,draw,fill=red!20] (#3) at (#1,#2) {$#4$};}
\nc\Znode[4] {\node[rounded corners,rectangle,draw,fill=blue!20] (Z#1#2) at (#3,#4) {$Z_{#1,#2}$};}
\nc\Snode[6] {\node[rounded corners,rectangle,draw,fill=green!20] (S#1#2) at (#3+#5,#4+#6) {$S_{#1,#2}$};}
\nc\squarelattice[3]{
\draw[#3] (#1,#2)--(#1+1,#2+1)--(#1,#2+2)--(#1-1,#2+1)--(#1,#2);
\foreach \x in {(#1,#2),(#1+1,#2+1),(#1,#2+2),(#1-1,#2+1),(#1,#2)} {\fill[#3]  \x circle(.1);}
}
\nc\understring[2]{\draw[thick] (#1,4) .. controls (#1,2) and (#2,2) .. (#2,0);}
\nc\overstring[2]{\draw[white,line width=2mm] (#1,4) .. controls (#1,2) and (#2,2) .. (#2,0); \draw[thick] (#1,4) .. controls (#1,2) and (#2,2) .. (#2,0);}
\nc\redoverstring[2]{\draw[white,line width=2mm] (#1,4) .. controls (#1,2) and (#2,2) .. (#2,0); \draw[red,thick] (#1,4) .. controls (#1,2) and (#2,2) .. (#2,0);}
\nc\halfunderstring[2]{\draw[thick] (#1,2) .. controls (#1,1) and (#2,1) .. (#2,0);}
\nc\halfoverstring[2]{\draw[white,line width=2mm] (#1,2) .. controls (#1,1) and (#2,1) .. (#2,0); \draw[thick] (#1,2) .. controls (#1,1) and (#2,1) .. (#2,0);}
\nc\redhalfoverstring[2]{\draw[white,line width=2mm] (#1,2) .. controls (#1,1) and (#2,1) .. (#2,0); \draw[red,thick] (#1,2) .. controls (#1,1) and (#2,1) .. (#2,0);}
\nc\understringx[4]{\draw[thick] (#1,#2) .. controls (#1,#2/2+#4/2) and (#3,#2/2+#4/2) .. (#3,#4);}
\nc\overstringx[4]{\draw[white,line width=2mm] (#1,#2) .. controls (#1,#2/2+#4/2) and (#3,#2/2+#4/2) .. (#3,#4); \draw[thick] (#1,#2) .. controls (#1,#2/2+#4/2) and (#3,#2/2+#4/2) .. (#3,#4);}
\nc\redoverstringx[4]{\draw[white,line width=2mm] (#1,#2) .. controls (#1,#2/2+#4/2) and (#3,#2/2+#4/2) .. (#3,#4); \draw[red,thick] (#1,#2) .. controls (#1,#2/2+#4/2) and (#3,#2/2+#4/2) .. (#3,#4);}
\nc\dotts[3]{\draw[dotted](#1+.35,#3)--(#2-.35,#3);}
\nc\braidrest[1]{^{\downarrow(#1)}}
\nc\V{\mathcal V}
\nc\DD{\mathcal D}
\nc\C{\mathcal C}
\nc\bN{{\bf N}}
\nc\tran\ol
\nc\CongH{\Cong_H}

\nc\overuarcx[3]{
\draw[white,line width=2mm](#1,4)arc(180:270:#3) (#1+#3,4-#3)--(#2-#3,4-#3) (#2-#3,4-#3) arc(270:360:#3);
\draw[thick](#1,4)arc(180:270:#3) (#1+#3,4-#3)--(#2-#3,4-#3) (#2-#3,4-#3) arc(270:360:#3);
}

\nc\overdarcx[3]{
\draw[white,line width=2mm](#1,0)arc(180:90:#3) (#1+#3,#3)--(#2-#3,#3) (#2-#3,#3) arc(90:0:#3);
\draw[thick](#1,0)arc(180:90:#3) (#1+#3,#3)--(#2-#3,#3) (#2-#3,#3) arc(90:0:#3);
}

\nc{\set}[2]{\{#1:#2\}}
\nc{\bigset}[2]{\big\{#1:#2\big\}}
\nc{\pres}[2]{\la#1:#2\ra}
\nc{\bigpres}[2]{\big\la#1:#2\big\ra}
\nc\sub\subseteq
\nc\mt\mapsto
\nc\al\alpha
\nc\be\beta
\nc\ga\gamma
\nc\de\delta
\nc\si\sigma
\nc\ep\epsilon
\nc\lam\lambda
\nc\ze\zeta
\nc\vk\varkappa
\nc\vt{\widetilde{\nu}}
\nc\vs\varsigma
\nc\ka\kappa
\rnc\th\theta
\nc\om\omega
\nc\ups\upsilon
\nc\Ga\Gamma
\nc\De\Delta
\nc\Si\Sigma
\nc\Om\Omega
\rnc\O{\mathcal O}
\nc\bit{\begin{itemize}}
\nc\eit{\end{itemize}}
\nc\ben{\begin{enumerate}[label=\textup{(\roman*)},leftmargin=7mm]}
\nc\BEN{\begin{enumerate}[label=\textup{(\Roman*)},leftmargin=7mm]}
\nc\bena{\begin{enumerate}[label=\textup{(\alph*)},leftmargin=7mm]}
\nc\een{\end{enumerate}}
\nc\bmc{\begin{multicols}}
\nc\emc{\end{multicols}}
\nc{\leqR}{\leq_{\R}}
\nc{\leqL}{\leq_{\L}}
\nc{\leqJ}{\leq_{\J}}
\nc{\leqK}{\leq_{\K}}
\nc{\geqR}{\geq_{\R}}
\nc{\geqL}{\geq_{\L}}
\nc{\geqJ}{\geq_{\J}}

\nc\bp{{\bf p}}
\nc\bq{{\bf q}}
\rnc\iff{\ \Leftrightarrow\ }
\rnc\implies{\ \Rightarrow\ }

\nc\pf{\begin{proof}}
\nc\epf{\end{proof}}
\nc\epfres{\hfill\qed}
\nc\epfreseq{\tag*{\qed}}

\let\oldproofname=\proofname
\renewcommand{\proofname}{\rm\bf{\oldproofname}}

\nc\AND{\qquad\text{and}\qquad}
\nc\WHERE{\qquad\text{where}\qquad}
\nc\ANd{\quad\text{and}\quad}
\nc\anD{\ \ \ \text{and}\ \ \ }
\nc\ANDSIM{\qquad\text{and similarly}\qquad}
\nc{\COMMA}{,\qquad}
\nc{\COMMa}{,\quad}
\nc\ul\underline
\nc\ol\overline
\nc\permdec[1]{#1^{\natural}}
\nc\ext[1]{#1^\textup{E}}
\nc\wh\widehat
\nc\normal\unlhd
\rnc\emptyset\varnothing
\nc{\pfitem}[1]{\medskip\noindent #1.}
\nc{\firstpfitem}[1]{#1.}
\nc{\pfcase}[1]{\medskip\noindent {\bf Case #1.}}
\nc{\pfstep}[1]{\medskip\noindent {\bf Step #1.}}
\nc\aftercases{\medskip\noindent}
\nc{\pfclaim}[1]{\medskip\noindent{\bf Claim #1.}}
\nc{\pfclaimnn}{\medskip\noindent{\bf Claim.} } 
\nc\afterclaim{\medskip}
\nc{\pfsubcase}[1]{\medskip\noindent {\bf Subcase #1.}}
\nc\im{\operatorname{im}}
\nc\LSUB{\operatorname{LSUB}}

\nc\pre\preceq
\nc\Y{\mathcal Y}
\nc\B{\mathcal B}
\nc\Z{\mathcal Z}
\nc\ZZ{\mathbb Z}
\nc\F{\mathcal F}
\nc\T{\mathcal T}
\nc\TT{\mathscr T}
\nc\PlP{\mathscr P\mathcal P}
\nc\I{\mathcal I}
\nc\Eq{\mathfrak{Eq}}
\nc\Part{\mathbb{P}}
\nc\cg[2]{(#1,#2)^\sharp}
\nc\Rev{\operatorname{Rev}}
\nc\cR{\mathcal R}
\nc\Ptop{P^\top}
\nc\Qtop{Q^\top}
\nc\la\langle
\nc\ra\rangle
\nc\tb[1]{\operatorname{Seq}(#1)}
\nc\RevX{{\Rev}\big([0,|X|],[0,|X|^+]\big)}

\makeatletter
\DeclareRobustCommand\widecheck[1]{{\mathpalette\@widecheck{#1}}}
\def\@widecheck#1#2{%
    \setbox\z@\hbox{\m@th$#1#2$}%
    \setbox\tw@\hbox{\m@th$#1%
       \widehat{%
          \vrule\@width\z@\@height\ht\z@
          \vrule\@height\z@\@width\wd\z@}$}%
    \dp\tw@-\ht\z@
    \@tempdima\ht\z@ \advance\@tempdima2\ht\tw@ \divide\@tempdima\thr@@
    \setbox\tw@\hbox{%
       \raise\@tempdima\hbox{\scalebox{1}[-1]{\lower\@tempdima\box
\tw@}}}%
    {\ooalign{\box\tw@ \cr \box\z@}}}
\makeatother

\nc\ch\widecheck

\newcommand{\darcx}[3]{\draw(#1,0)arc(180:90:#3) (#1+#3,#3)--(#2-#3,#3) (#2-#3,#3) arc(90:0:#3);}
\newcommand{\darc}[2]{\darcx{#1}{#2}{.4}}
\newcommand{\uarcx}[3]{\draw(#1,2)arc(180:270:#3) (#1+#3,2-#3)--(#2-#3,2-#3) (#2-#3,2-#3) arc(270:360:#3);}

\newcommand{\uarc}[2]{\uarcx{#1}{#2}{.4}}

\nc{\buv}[1]{\fill (#1,2)circle(.18);}
\nc{\buvs}[1]{{
\foreach \x in {#1}
{ \buv{\x}}
}}
\nc{\blv}[1]{\fill (#1,0)circle(.18);}
\nc{\blvs}[1]{{
\foreach \x in {#1}
{ \blv{\x}}
}}

\nc{\uarcs}[1]{
{\foreach \x/\y in {#1}
{ \uarc{\x}{\y} }
}
}

\nc{\darcs}[1]{
{\foreach \x/\y in {#1}
{ \darc{\x}{\y} }
}
}
\nc{\darcxhalf}[3]{\draw(#1,0)arc(180:90:#3) (#1+#3,#3)--(#2,#3) ;}
\nc{\darchalf}[2]{\darcxhalf{#1}{#2}{.4}}
\nc{\uarcxhalf}[3]{\draw(#1,2)arc(180:270:#3) (#1+#3,1.5-#3)--(#2,1.5-#3) ;}
\nc{\uarchalf}[2]{\uarcxhalf{#1}{#2}{.4}}
\nc{\colv}[3]{\fill[#3] (#1,#2)circle(.17);}

\nc{\uvert}[1]{\fill (#1,2)circle(.2);}
\rnc{\lvert}[1]{\fill (#1,0)circle(.1);}

\nc{\custpartn}[3]{{\lower1.4 ex\hbox{
\begin{tikzpicture}[scale=.3]
\foreach \x in {#1}
{ \uvert{\x}  }
\foreach \x in {#2}
{ \lvert{\x}  }
#3 \end{tikzpicture}
}}}

\newcommand{\JJ}{\mathcal{J}} 
\renewcommand{\S}{\mathcal{S}}

\nc\MYZ{\mathcal M_{Y\cup Z}}

\newcommand{\N}{\mathbb{N}}

\nc\HH{\mathcal H}

\newcommand{\Cong}{\operatorname{Cong}}

\newcommand{\id}{\operatorname{id}}

\nc\congsquare[4]
{\fill[black!15] (#1*6,#2*6,#3*8)--({#1*6+(#4)*6},#2*6,#3*8)--({#1*6+(#4)*6},{#2*6+(#4)*6},#3*8)--(#1*6,{#2*6+(#4)*6},#3*8)--(#1*6,#2*6,#3*8);
\foreach \x in {0,1,...,#4} {\draw (#1*6,#2*6+\x*6,#3*8)--({#1*6+(#4)*6},#2*6+\x*6,#3*8); \draw (#1*6+\x*6,#2*6,#3*8)--(#1*6+\x*6,{#2*6+(#4)*6},#3*8);}}

\nc\congsquareconnections[7]
{\foreach \x in {0,...,#7} \foreach \y in {0,...,#7} {\draw(#1*6+6*\x,#2*6+6*\y,#3*8)--(#4*6+6*\x,#5*6+6*\y,#6*8);}}

\numberwithin{equation}{section}

\newtheorem{lemma}[equation]{Lemma}

\newtheorem{prop}[equation]{Proposition}
\newtheorem{con}[equation]{Conjecture}
\newtheorem*{thm*}{Theorem}

\theoremstyle{definition}

\newtheorem{rem}[equation]{Remark}

\newtheorem{prob}{Problem}

\newcommand{\sm}{\setminus}

\begin{document}

\title{\vspace{-.3cm}On the enumeration of integer tetrahedra}
\author{}
\date{}

\maketitle
~\vspace{-2cm}
\begin{center}
{\large 
James East,%
\hspace{-.3em}\footnote{\label{footnote:JE}Centre for Research in Mathematics and Data Science, Western Sydney University, Australia. {\tt J.East@WesternSydney.edu.au}, {\tt L.Park@WesternSydney.edu.au}.}\footnote{Supported by ARC Future Fellowship FT190100632.}
Michael Hendriksen,%
\hspace{-.3em}\footnote{School of Mathematics and Statistics, University of Melbourne, Australia. {\tt michael.hendriksen@unimelb.edu.au}.}\footnote{Parts of this research were carried out when this author was a postgraduate student at Western Sydney University, and a postdoctoral fellow at Heinrich Heine Universit\"at D\"usseldorf.}
Laurence Park\textsuperscript{\ref{footnote:JE}}}
\end{center}

\begin{abstract}
We consider the problem of enumerating integer tetrahedra of fixed perimeter (sum of side-lengths) and/or diameter (maximum side-length), up to congruence.  As we will see, this problem is considerably more difficult than the corresponding problem for triangles, which has long been solved.  We expect there are no closed-form solutions to the tetrahedron enumeration problems, but we explore the extent to which they can be approached via classical methods, such as orbit enumeration.  We also discuss algorithms for computing the numbers, and present several tables and figures that can be used to visualise the data.  Several intriguing patterns seem to emerge, leading to a number of natural conjectures.  The central conjecture is that the number of integer tetrahedra of perimeter $n$, up to congruence, is asymptotic to $n^5/C$ for some constant $C\approx 229000$.



\end{abstract}

\setcounter{tocdepth}{1}
\tableofcontents

\section{Introduction}\label{s:intro}

Geometry and combinatorics have many natural meeting points.  Arguably the most ancient known example is the application of right-angled integer triangles in Babylonian architecture and agriculture.  The classification of all such right triangles is recorded in Book X of Euclid's Elements.
Enumeration of arbitrary integer triangles goes back at least to the 1979 paper of Jordan, Walch and Wisner \cite{JWW1979}, and we have the following elegant result of Honsberger \cite{Honsberger1985}, which has been proved in a variety of ways \cite{Honsberger1985,JM2000,Hirschhorn2000,Hirschhorn2003,EN2019,EN2019_2}:

\begin{thm*}\label{t:triangles}
The number of integer triangles with perimeter $n$, up to congruence, is the nearest integer to $\frac{n^2}{48}$ if $n$ is even, or to $\frac{(n+3)^2}{48}$ if $n$ is odd.
\end{thm*}

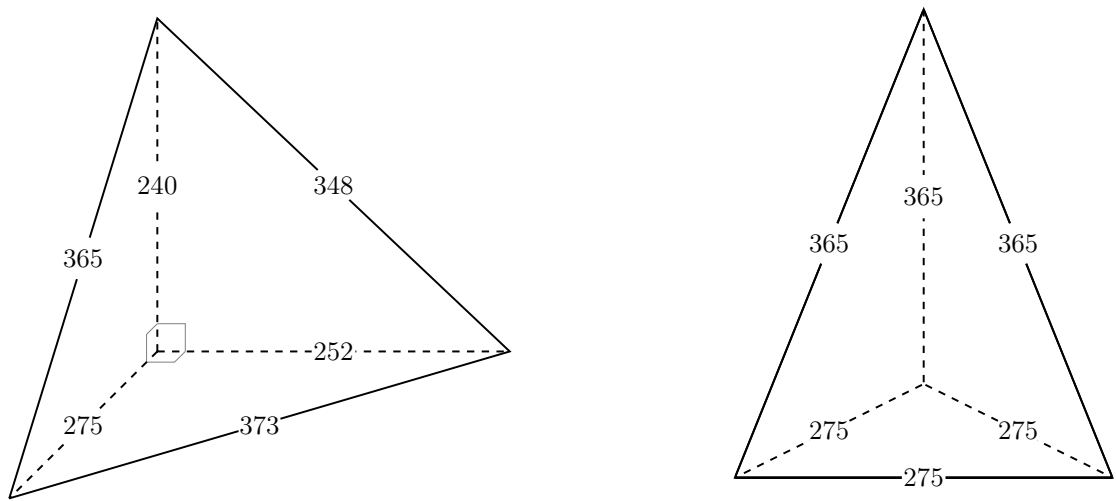
\begin{figure}[ht]
\begin{center}
\scalebox{0.92}{
\begin{tikzpicture}[scale=2,inner sep=0]
\draw[gray] (.20,0,0)--(.20,.20,0)--(0,.20,0);
\draw[gray] (0,.20,0)--(0,.20,.20)--(0,0,.20);
\draw[gray] (.20,0,0)--(.20,0,.20)--(0,0,.20);
\coordinate  (A) at (0,0,0);
\coordinate  (B) at (0,0,2.75);
\coordinate  (C) at (2.52,0,0);
\coordinate (D) at (0,2.40,0);
\draw[thick] (B) --node[circle,fill=white]{373} (C) --node[circle,fill=white]{348} (D) --node[circle,fill=white]{365} (B);
\draw[dashed,thick] (A) --node[circle,fill=white]{275} (B) (A) --node[circle,fill=white]{252} (C) (A) --node[circle,fill=white]{240} (D);
\end{tikzpicture}
\qquad\qquad\qquad\qquad
\begin{tikzpicture}[scale=2.7,inner sep=0]
\draw[thick] (1,-.5)--(-1,-.5)--(0,2)--(1,-.5)--(-1,-.5);
\node  (A) at (0,0) {};
\node  (B) at (1,-.5) {};
\node  (C) at (-1,-.5) {};
\node (D) at (0,2) {};
\draw[thick] (B) --node[circle,fill=white]{275} (C) --node[circle,fill=white]{365} (D) --node[circle,fill=white]{365} (B);
\draw[dashed,thick] (A) --node[circle,fill=white]{275} (B) (A) --node[circle,fill=white]{275} (C) (A) --node[circle,fill=white]{365} (D);
\end{tikzpicture}
}
\caption{Two integer tetrahedra.  The one on the left is tri-rectangular.  The one on the right is invariant under a $120^\circ$ rotation about a vertical axis.}
\label{f:T1}
\end{center}
\end{figure}

Two of the most natural extensions of this triangle enumeration problem are to consider integer polygons (increasing the number of sides), or integer tetrahedra (moving up a dimension).  The former was treated in \cite{EN2019}, and the current article considers the latter.  Our main guiding problem is the following.  By an \emph{integer tetrahedron} we mean a (non-degenerate) tetrahedron whose sides are all of integer length, as in Figure~\ref{f:T1}.  The \emph{perimeter} of a tetrahedron is the sum of its six side-lengths.

\begin{prob}\label{prob:tn}
Calculate the number $\tn$ of integer tetrahedra with perimeter $n$, up to congruence: i.e., combinations of rotations, translations and reflections.  
\end{prob}

One might hope that $\tn$ is given by a similar formula to the triangle sequence in Honsberger's Theorem above.  As we will see, however, this is \emph{very} far from the truth.  

Philip Benjamin has computed~$\tn$ for $n\leq30$, but as far as we know has not published his methods; see \cite[Sequence A208454]{OEIS}.
Sascha Kurz \cite{Kurz2007} has considered the related problem of enumeration by \emph{diameter}, defined to be the maximum of the six side-lengths:

\begin{prob}\label{prob:td}
Calculate the number $\td$ of integer tetrahedra with diameter $d$, up to congruence.
\end{prob}

Kurz has computed $\td$ for $d\leq1000$ in \cite{Kurz2007}; see also \cite[Sequence A097125]{OEIS}.  The article \cite{Kurz2007} gives a lot of detail about Kurz's methods and algorithms.  No formula is given or conjectured, though an exact expression is given for a related set of (orbits of) matrices; this leads to a conjectural asymptotic formula, as we discuss in Sections \ref{ss:as} and \ref{ss:max}.

It is not hard to show that the number of integer triangles with diameter $d$, up to congruence, is equal to $\big\lfloor\frac{(d+1)^2}4\big\rfloor$; see \cite[Sequence A002620]{OEIS}.

There is also of course the following natural problem, combining perimeter and diameter:

\begin{prob}\label{prob:tnd}
Calculate the number $\tnd$ of integer tetrahedra with perimeter $n$ and diameter~$d$, up to congruence.
\end{prob}

Clearly a solution to Problem \ref{prob:tnd} would yield solutions to Problems \ref{prob:tn} and \ref{prob:td}, since
\[
\tn = \sum_d\tnd \AND \td = \sum_n\tnd.
\]
The non-zero terms in these sums occur for $\lceil\frac n6\rceil \leq d \leq \lfloor\frac{n-3}3\rfloor$ and $3d+3\leq n\leq 6d$, as we will show in Lemma \ref{l:d}.

Apart from this introduction, and a brief conclusion, the paper contains three further sections.
In Section \ref{s:BL} we set up some ideas that (in principle) allow the calculation of $\tn$, $\td$ and $\tnd$ via Burnside's Lemma, and we make some partial progress by giving explicit formulas for some of the relevant parameters.
In Section \ref{s:computed} we discuss algorithms for computing the numbers $\tn$, $\td$ and $\tnd$, and give several tables and graphs of computed values; more data can be found at \cite{web}.
In Section \ref{s:obs} we explore some intriguing patterns that seem to emerge from an examination of the data.  A number of conjectures/open problems are stated.  We hope that these provide inspiration for future studies.

\section{Enumeration via Burnside's Lemma}\label{s:BL}

In this section we focus on the numbers $\tn$, though all we say can easily be adapted to work for the numbers $\td$ or $\tnd$ instead.  In Remark \ref{r:td_tnd} we indicate the modifications needed for these.

\begin{figure}[ht]
\begin{center}
\begin{tikzpicture}
\node () at (0,0) {\graphhack ABCabc{0.5}{G}};
\node () at (30:6) {\graphhack CABcab{0.5}{(1,2,3)\cdot G}};
\node () at (210:6) {\graphhack AbcaBC{0.5}{(1,4)(2,3)\cdot G}};
\node () at (330:6) {\graphhack acBACb{0.5}{(1,2,4,3)\cdot G}};
\node () at (150:6) {\graphhack ACBacb{0.5}{(2,3)\cdot G}};
\draw[-{latex}] (330:2.666)--(330:4.333);
\draw[-{latex}] (210:2.666)--(210:4.333);
\draw[-{latex}] (30:1.666)--(30:3.333);
\draw[-{latex}] (150:1.666)--(150:3.333);
\end{tikzpicture}
\vspace{-5truemm}
\caption{The graph $G\in\G$ (centre), as well as $\si\cdot G$ for various permutations $\si$ from $\S_4$.}
\label{f:act}
\end{center}
\end{figure}
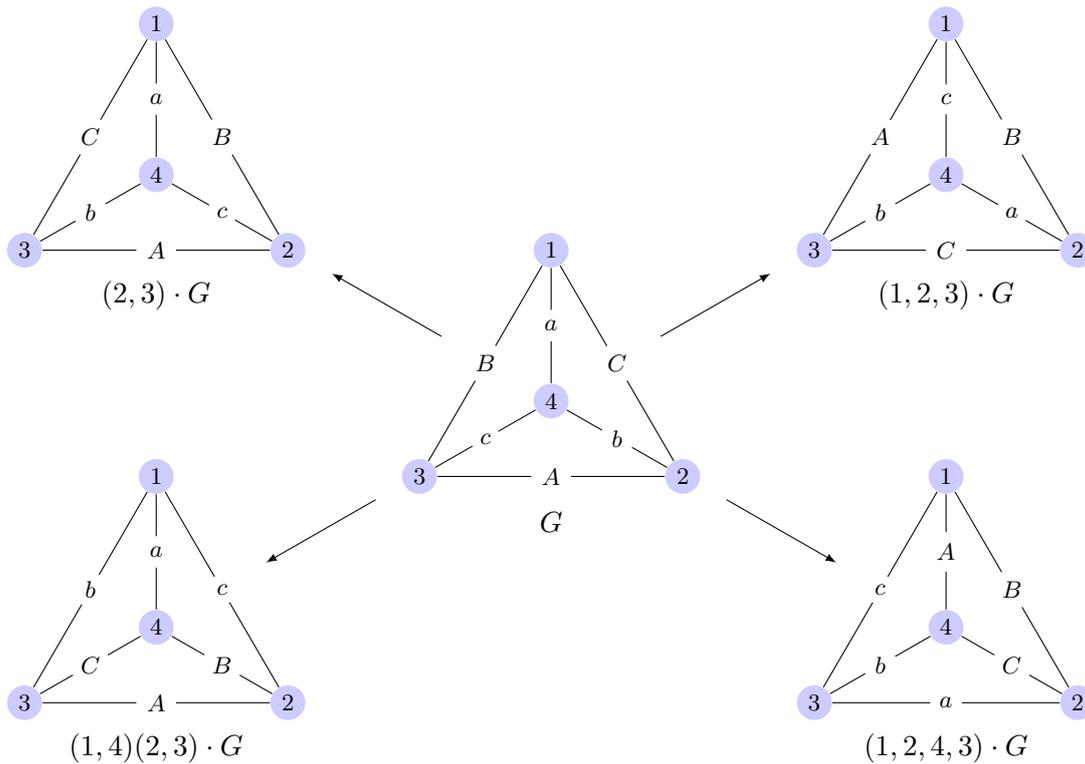

Let $\G$ be the set of all edge-labelled graphs whose underlying unlabelled graph is the complete graph on vertex set $\{1,2,3,4\}$, and whose labels all belong to $\N=\{1,2,3,\ldots\}$.  For $G\in\G$, we denote the label of the edge $\{i,j\}$ by $G(i,j)=G(j,i)$.  The symmetric group $\S_4$ has a natural action on $\G$, denoted $(\si,G)\mt\si\cdot G$, and induced by permuting the vertices.  For $G\in\G$ and $\si\in\S_4$, and for distinct $i,j\in\{1,2,3,4\}$, we have 
\begin{equation}\label{e:act}
(\si\cdot G)(i,j) = G(\si^{-1}(i),\si^{-1}(j)).
\end{equation}
For example, Figure \ref{f:act} illustrates the action of the permutations $(2,3)$, $(1,2,3)$, $(1,2,4,3)$ and $(1,4)(2,3)$, written in standard cycle notation.

\begin{figure}[ht]
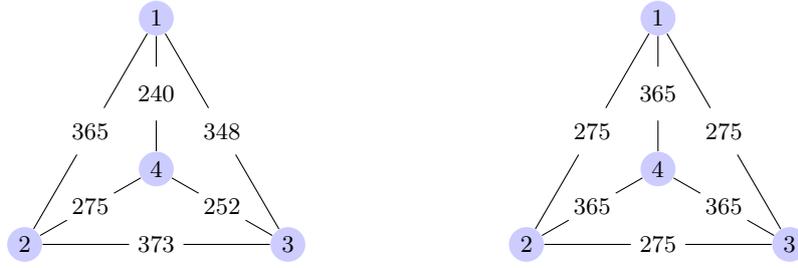

\begin{center}
\graph{373}{348}{365}{240}{275}{252}{.5}
\qquad\qquad\qquad
\graph{275}{275}{275}{365}{365}{365}{.5}
\caption{
Graphs from $\G$ representing the tetrahedra pictured in Figure \ref{f:T1}.}
\label{f:T2}
\end{center}
\end{figure}

An integer tetrahedron $T$ may be represented by a graph $G$ from $\G$ as follows.  Choose an ordering $1,2,3,4$ on the corners of~$T$, and let the label of the edge $\{i,j\}$ from $G$ be the corresponding side-length from $T$.  Different orderings on the vertices typically lead to different graphs, so that $T$ may be represented by up to $24$ such graphs.  In fact, these graphs representing $T$ are precisely those belonging to the orbit of $G$ under the action of $\S_4$ given in~\eqref{e:act}.  For example, the tetrahedra from Figure \ref{f:T1} may be represented by the graphs in Figure \ref{f:T2}; the orbits of these graphs have size $24$ and $4$, respectively.  
It is clear that two tetrahedra are congruent if and only if they are represented by the same set of graphs from $\G$.

\begin{figure}[ht]
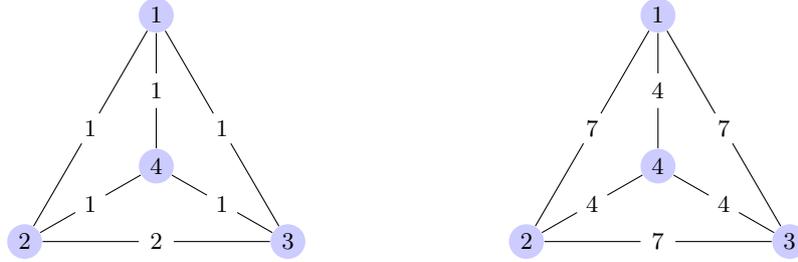

\begin{center}
\graph211111{.5}
\qquad\qquad\qquad
\graph777444{.5}
\caption{Graphs from $\G$ not corresponding to tetrahedra.}
\label{f:notT}
\end{center}
\end{figure}

\begin{figure}[ht]
\begin{center}
\begin{tikzpicture}[scale=0.5,inner sep=0]
\draw (90:4.04) --node[circle,inner sep=2pt, fill=white]{\footnotesize $7$} (330:4.04) --node[circle,inner sep=2pt, fill=white]{\footnotesize $7$} (210:4.04) --node[circle,inner sep=2pt, fill=white]{\footnotesize $7$} (90:4.04);
\draw (90:4.04) --node[circle,inner sep=2pt, fill=white]{\footnotesize $4$} (30:3.957) --node[circle,inner sep=2pt, fill=white]{\footnotesize $4$} (330:4.04);
\draw (90:4.04) --node[circle,inner sep=2pt, fill=white]{\footnotesize $4$} (150:3.957) --node[circle,inner sep=2pt, fill=white]{\footnotesize $4$} (210:4.04);
\draw (210:4.04) --node[circle,inner sep=2pt, fill=white]{\footnotesize $4$} (270:3.957) --node[circle,inner sep=2pt, fill=white]{\footnotesize $4$} (330:4.04);
\draw[-{latex}] (0:6)--(0:10);
\begin{scope}[shift={(15,0)}]
\draw (90:4.04) --node[circle,inner sep=2pt, fill=white]{\footnotesize $7$} (330:4.04) --node[circle,inner sep=2pt, fill=white]{\footnotesize $7$} (210:4.04) --node[circle,inner sep=2pt, fill=white]{\footnotesize $7$} (90:4.04);
\draw (90:4.04)--(30:0.085)--(330:4.04);
\draw (90:4.04)--(150:0.085)--(210:4.04);
\draw (210:4.04)--(270:0.085)--(330:4.04);
\node[circle,inner sep=2pt, fill=white] () at (90:1.7) {\footnotesize $4$};
\node[circle,inner sep=2pt, fill=white] () at (210:1.7) {\footnotesize $4$};
\node[circle,inner sep=2pt, fill=white] () at (330:1.7) {\footnotesize $4$};
\end{scope}
\end{tikzpicture}
\caption{Attempting to construct a tetrahedron from the graph in Figure \ref{f:notT} (right).}
\label{f:notT2}
\end{center}
\end{figure}
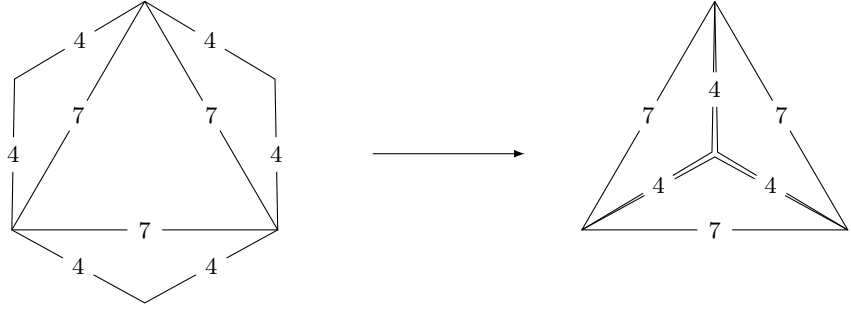

Not every graph from $\G$ corresponds to a tetrahedron in the above manner.  For example, we claim that this is the case for the two graphs pictured in Figure \ref{f:notT}.  This is clear for the left-hand graph, as there is no triangle with edges $(1,1,2)$.  However, the right-hand graph cannot be ruled out so easily, as the triples $(7,7,7)$ and $(7,4,4)$ do indeed correspond to triangles.  Rather, the problem here is that if one tries to fold the ``net'' shown in Figure \ref{f:notT2} (left) into a tetrahedron, then the tips of the three $(7,4,4)$ triangles do not meet.

For $n\in\N$, let~$\T_n$ be the set of graphs from $\G$ corresponding to some tetrahedron of perimeter~$n$.
The number of integer tetrahedra of perimeter $n$, up to congruence, is then given by the number~$\tn = |\T_n/\S_4|$ of orbits of $\T_n$ under the action of $\S_4$ given in \eqref{e:act}.  Burnside's Lemma (cf.~\cite[p.~246]{Cameron1994}) then gives
\begin{equation}\label{e:BS1}
\tn = |\T_n/\S_4| = \frac1{24}\sum_{\si\in\S_4}\fix(\si).
\end{equation}
Here, for $\si\in \S_4$, $\fix(\si)$ is the cardinality of the set 
\begin{equation}\label{e:Fix}
\Fix(\si) = \set{G\in\T_n}{\si\cdot G=G}.
\end{equation}
To enumerate tetrahedra up to rotations and translations only, we would be looking at orbits under the restricted action of the alternating group $\A_4\sub\S_4$, since even and odd permutations correspond to rotations and reflections, respectively.  The number $\tn'$ of integer tetrahedra of perimeter~$n$, up to this kind of restricted congruence, is given by
\begin{equation}\label{e:BS2}
\tn' = |\T_n/\A_4| = \frac1{12}\sum_{\si\in\A_4}\fix(\si).
\end{equation}
Thus, to calculate the numbers $\tn$ and $\tn'$, it is enough to identify the set $\T_n$, and calculate the parameters $\fix(\si)$ for each $\si\in\S_4$.  For the former, we have the following:

\begin{prop}\label{p:Tn}
The set $\T_n$ consists of all graphs $G$ from $\G$ whose edge-labels, as shown in Figure \ref{f:act} (centre), satisfy:
\begin{enumerate}[label=\textup{(T\arabic*)},leftmargin=12mm]
\bmc2
\item \label{T1} $A+B+C+a+b+c=n$,
\item \label{T2} $A+B+C>2\max(A,B,C)$,
\item \label{T3} $A+b+c>2\max(A,b,c)$,
\item \label{T4} $a+B+c>2\max(a,B,c)$,
\item \label{T5} $a+b+C>2\max(a,b,C)$,
\item \label{T6} $(a^2-B^2-c^2+2x_1x_2)^2 < 4Y_1Y_2$, where 
\bit
\item $x_1 = \tfrac{A^2+B^2-C^2}{2A}$,
\item $x_2 = \tfrac{A^2+c^2-b^2}{2A}$.
\item $Y_1 = B^2-x_1^2$,
\item $Y_2 = c^2-x_2^2$.
\eit
\emc\een
\end{prop}

\pf
First note that three positive real numbers $x,y,z$ can be the sides of a triangle if and only if the sum of the smaller two is greater than the largest; this is equivalent to $x+y+z>2\max(x,y,z)$.  Thus, consulting Figure \ref{f:act}, it is clear that $\T_n$ is a subset of
\[
\T_n' = \bigset{G\in\G}{\text{\ref{T1}--\ref{T5} hold}}.
\]
Now consider a graph $G$ from $\T_n'$, with edges labelled as in Figure \ref{f:act}.  Then $G$ corresponds to a tetrahedron if and only if the following procedure can be carried out:
\bit
\item Begin with the quadrilateral shown in Figure \ref{f:T3} (left).
\item Keeping the $(A,B,C)$-triangle fixed in place, fold along the $x$-axis by some angle $0<\th<\pi$ until the free tips of the two triangles are $a$ units apart, as shown in Figure \ref{f:T3} (middle).  (We have strict inequalities for $\th$ to ensure that the tetrahedron is not degenerate.)
\eit
Figure \ref{f:T3} (right) displays the two original triangles, as well as the result of folding by the angle $\th=\pi$, and defines three points, $P,Q,R\in\mathbb R^2$.  Let the coordinates of these points be $P=(x_1,y_1)$, $Q=(x_2,y_2)$ and $R=(x_2,-y_2)$.  Also, put $Y_1=y_1^2$ and $Y_2=y_2^2$.  It is easy to check that $x_1,x_2,Y_1,Y_2$ are as given in \ref{T6}: e.g., consider $P$ as the intersection of the circles $x^2+y^2=B^2$ and $(x-A)^2+y^2=C^2$.

We claim that the above folding procedure can be carried out if and only if ${|PQ|<a<|PR|}$.  Indeed, if we denote by $R_\th\in\mathbb R^3$ the tip of the moving triangle after folding by the angle $0\leq\th\leq\pi$, then $R_\th=(x_2,-y_2\cos\th,y_2\sin\th)$.  It follows that $|PR_\th|^2=(x_1-x_2)^2+y_1^2+y_2^2+2y_1y_2\cos\th$ is a smooth, decreasing function of $0\leq\th\leq\pi$.

With the claim established, it remains to observe that the inequality $|PQ|<a<|PR|$ is equivalent to that in~\ref{T6}, as
\begin{align*}
|PQ|<a<|PR| &\iff |PQ|^2<a^2<|PR|^2 \\
&\iff (x_1-x_2)^2+(y_1-y_2)^2 < a^2 < (x_1-x_2)^2+(y_1+y_2)^2\\
&\iff y_1^2+y_2^2-2y_1y_2 < a^2 - (x_1-x_2)^2< y_1^2+y_2^2+2y_1y_2\\
&\iff -2y_1y_2 < a^2 - (x_1-x_2)^2 - y_1^2 - y_2^2 < 2y_1y_2\\
&\iff (a^2 - (x_1-x_2)^2 - y_1^2 - y_2^2)^2 < 4y_1^2y_2^2\\
&\iff (a^2-B^2-c^2+2x_1x_2)^2 < 4y_1^2y_2^2.  \qedhere
\end{align*}
\epf

\begin{figure}[ht]
\begin{center}
\begin{tikzpicture}[scale=.6,inner sep=1.0]
\nc\AAA5
\nc\BBB6
\nc\CCC7
\nc\bbb5
\nc\ccc6
\begin{scope}[shift={(0,0)}]
\coordinate (P) at (0,0);
\coordinate (R) at (\AAA,0);
\coordinate (S) at ({(\AAA^2+\BBB^2-\CCC^2)/(2*\AAA)},{sqrt(\BBB^2-((\AAA^2+\BBB^2-\CCC^2)/(2*\AAA))^2)});
\coordinate (Q) at ({(\AAA^2+\ccc^2-\bbb^2)/(2*\AAA)},{-sqrt(\ccc^2-((\AAA^2+\ccc^2-\bbb^2)/(2*\AAA))^2)});
\coordinate (Q') at ({(\AAA^2+\ccc^2-\bbb^2)/(2*\AAA)},{sqrt(\ccc^2-((\AAA^2+\ccc^2-\bbb^2)/(2*\AAA))^2)});
\coordinate (Q'') at (4.5,3.5);
\fill[blue!20] (P) -- (Q) --(R); \fill[red!20] (R) -- (S) -- (P); 
\draw[blue,ultra thick] (P) --node[circle,fill=white]{$c$} (Q) --node[circle,fill=white]{$b$} (R); \draw[red,ultra thick] (R) --node[circle,fill=white]{$C$} (S) --node[circle,fill=white]{$B$} (P); 
\draw[->] (-1,0)--(6,0); \node () at (6,-.4) {\footnotesize $x$};
\draw[->] (0,-1)--(0,5); \node () at (-.4,5) {\footnotesize $y$};
\draw[ultra thick] (P) --node[circle,fill=white]{$A$} (R);
\end{scope}
\begin{scope}[shift={(10,0)}]
\coordinate (P) at (0,0);
\coordinate (R) at (\AAA,0);
\coordinate (S) at ({(\AAA^2+\BBB^2-\CCC^2)/(2*\AAA)},{sqrt(\BBB^2-((\AAA^2+\BBB^2-\CCC^2)/(2*\AAA))^2)});
\coordinate (Q) at ({(\AAA^2+\ccc^2-\bbb^2)/(2*\AAA)},{-sqrt(\ccc^2-((\AAA^2+\ccc^2-\bbb^2)/(2*\AAA))^2)});
\coordinate (Q') at ({(\AAA^2+\ccc^2-\bbb^2)/(2*\AAA)},{sqrt(\ccc^2-((\AAA^2+\ccc^2-\bbb^2)/(2*\AAA))^2)});
\coordinate (Q'') at (4.7,3);
\coordinate (Q''') at (2.35,1.5);
\coordinate (Q'''') at (4.23,2.7);
 \fill[red!20] (R) -- (S) -- (P);  
\draw[red,ultra thick] (R) --node[circle,fill=white]{$C$} (S) --node[circle,fill=white]{$B$} (P); 
 \fill[blue!20] (P) -- (Q'') --(R);
\draw[dashed,red!50,ultra thick] (R) -- ($ (R) !.385! (S) $); 
\draw[dotted, ultra thick] (Q'') --node[circle,fill=white]{$a$} (S);
\draw[blue,ultra thick] (P) --node[circle,fill=white]{$c$} (Q'') --node[circle,fill=white]{$b$} (R); 
\draw[->] (-1,0)--(6,0); \node () at (6,-.4) {\footnotesize $x$};
\draw[->] (0,-1)--(0,5); \node () at (-.4,5) {\footnotesize $y$};
\draw[ultra thick] (P) --node[circle,fill=white]{$A$} (R);
\end{scope}
\begin{scope}[shift={(20,0)}]
\coordinate (P) at (0,0);
\coordinate (R) at (\AAA,0);
\coordinate (S) at ({(\AAA^2+\BBB^2-\CCC^2)/(2*\AAA)},{sqrt(\BBB^2-((\AAA^2+\BBB^2-\CCC^2)/(2*\AAA))^2)});
\coordinate (Q) at ({(\AAA^2+\ccc^2-\bbb^2)/(2*\AAA)},{-sqrt(\ccc^2-((\AAA^2+\ccc^2-\bbb^2)/(2*\AAA))^2)});
\coordinate (Q') at ({(\AAA^2+\ccc^2-\bbb^2)/(2*\AAA)},{sqrt(\ccc^2-((\AAA^2+\ccc^2-\bbb^2)/(2*\AAA))^2)});
\coordinate (Q'') at (4.5,3.5);
\fill[blue!20] (P) -- (Q) --(R); \fill[red!20] (R) -- (S) -- (P); 
\draw[blue,ultra thick] (P) -- (Q) -- (R); 
\draw[red,ultra thick] (R) -- (S) -- (P); 
\fill[blue!20] (P) -- (Q') --(R); 
\draw[dashed,red!50,ultra thick] (R) -- ($ (R) !.62! (S) $); 
\draw[ultra thick] (P) -- (R);
\draw[blue,ultra thick] (P) -- (Q') -- (R); 
\draw[->] (-1,0)--(6,0); \node () at (6,-.4) {\footnotesize $x$};
\draw[->] (0,-1)--(0,5); \node () at (-.4,5) {\footnotesize $y$};
\node () at (1,6.2) {$P$};
\node () at (4,5) {$Q$};
\node () at (4,-5) {$R$};
\end{scope}
\end{tikzpicture}
\caption{Creating an integer tetrahedron; see the proof of Proposition \ref{p:Tn} for further details.}
\label{f:T3}
\end{center}
\end{figure}

\begin{rem}
Items \ref{T4} and \ref{T5} from Proposition \ref{p:Tn} are actually unnecessary.  Indeed, given \ref{T1}--\ref{T3}, and in the notation of the above proof (cf.~Figure \ref{f:T3}), \ref{T6} is equivalent to the inequality $|PQ|<a<|PR|$.  This ensures that the triples $(a,B,c)$ and $(a,b,C)$ correspond to triangles, which implies \ref{T4} and \ref{T5}.

Alternative formulations of \ref{T6} also exist.  On \cite[p.~3]{Kurz2007}, Kurz attributes to Menger \cite{Menger1928} the equivalent inequality involving a determinant:
\[
\left|
\begin{matrix}
0&A^2&C^2&b^2&1\\
A^2&0&B^2&c^2&1\\
C^2&B^2&0&a^2&1\\
b^2&c^2&a^2&0&1\\
1&1&1&1&0
\end{matrix}
\right| > 0.
\]
See also \cite{WD2009} for an elementary proof of this inequality, and a discussion of its history.  Although the inequalities in Proposition \ref{p:Tn} are less elegant (and symmetrical) than the above determinant inequality, they seem to be quicker to check computationally.
\end{rem}

Now that we have characterised the set $\T_n$, it remains to calculate the values of $\fix(\si)$ for each~${\si\in \S_4}$.  By symmetry we only have to do this for $\si$ being one of the following:
%
\ben\bmc5
\item \label{fix1} $\id_4$,
\item \label{fix2} $(1,2,3)$,
\item \label{fix3} $(1,4)(2,3)$,
\item \label{fix4} $(2,3)$,
\item \label{fix5} $(1,2,4,3)$.
\emc\een
In these cases, the number of elements of $\S_4$ with the same values of $\fix(\si)$ are, respectively,
\ref{fix1}~1,
\ref{fix2}~8,
\ref{fix3}~3,
\ref{fix4}~6,
\ref{fix5}~6.
Moreover, $\Fix(\si)$ may be characterised as the graphs $G$ from $\T_n$ whose labels (as in Figure \ref{f:act}) satisfy the following constraints, respectively:
\ben
\bmc2
\item NA,
\item $A=B=C$ and $a=b=c$,
\item $B=b$ and $C=c$,
\item $B=C$ and $b=c$,
\item $A=a$ and $B=C=b=c$.
\item[] ~
\emc
\een
In case \ref{fix1} we have $\fix(\id_4)=|\T_n|$.  
The next two results give exact formulas for $\fix(\si)$ when $\si$ has type~\ref{fix2} or \ref{fix5}.  We do not currently have formulas for types \ref{fix1}, \ref{fix3} or \ref{fix4}; nor do we know if ``neat'' such formulas exist.

\newpage

\begin{lemma}\label{l:fix_123}
If $\si=(1,2,3)$, then
\[
\fix(\si) = \begin{cases}
0 &\text{if $n\not\equiv0\Mod3$}\\
\big\lfloor\frac n{3+\sqrt3}\big\rfloor &\text{if $n\equiv0\Mod3$.}
\end{cases}
\]
\end{lemma}

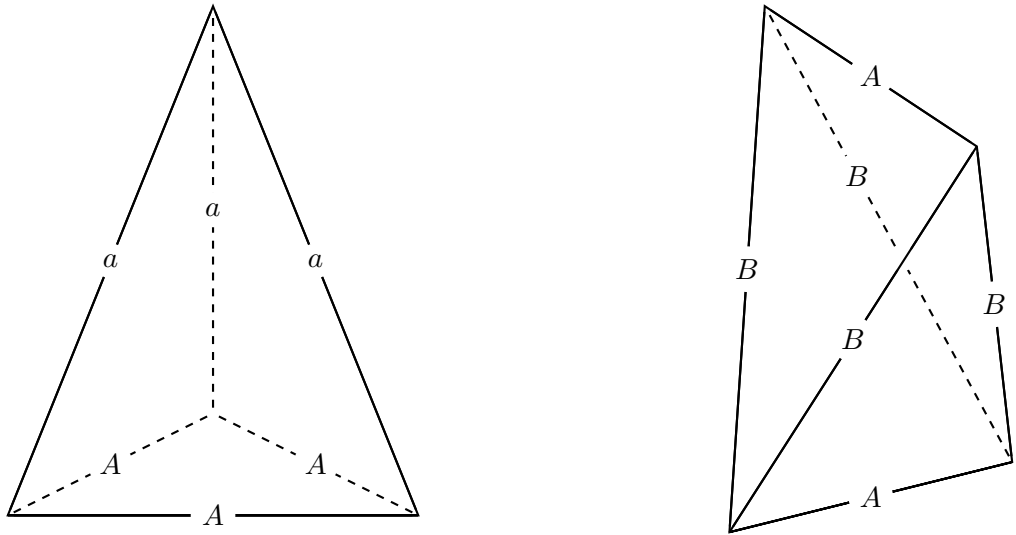
\begin{figure}[ht]
\begin{center}
\begin{tikzpicture}[scale=2.7,inner sep=2]
\draw[thick] (1,-.5)--(-1,-.5)--(0,2)--(1,-.5)--(-1,-.5);
\node  (A) at (0,0) {};
\node  (B) at (1,-.5) {};
\node  (C) at (-1,-.5) {};
\node (D) at (0,2) {};
\draw[thick] (B) --node[circle,fill=white]{$A$} (C) --node[circle,fill=white]{$a$} (D) --node[circle,fill=white]{$a$} (B);
\draw[dashed,thick] (A) --node[circle,fill=white]{$A$} (B) (A) --node[circle,fill=white]{$A$} (C) (A) --node[circle,fill=white]{$a$} (D);
\end{tikzpicture}
\qquad\qquad\qquad\qquad\qquad
\begin{tikzpicture}[scale=0.93,inner sep=2]
\draw[thick] (-2,-.5)--(2,.5)--(2-.5,6-1)--(-2+.5,6+1)--(-2,-.5)--(2,.5);
\draw[thick] (-2,-.5)--(2-.5,6-1);
\node  (A) at (-2,-.5) {};
\node  (B) at (2,.5) {};
\node (C) at (2-.5,6-1) {};
\node  (D) at (-2+.5,6+1) {};
\draw[dashed,thick] (B) --node[distance from start=120, circle,fill=white]{$B$} (D);
\fill[white] (0.45,3.4) circle (0.18);
\draw[thick] (A) --node[circle,fill=white]{$A$} (B) --node[circle,fill=white]{$B$} (C) --node[circle,fill=white]{$A$} (D) --node[circle,fill=white]{$B$} (A) --node[circle,fill=white]{$B$} (C);
\end{tikzpicture}
\caption{Tetrahedra from the proofs of Lemmas \ref{l:fix_123} and \ref{l:fix_1243}.}
\label{f:T4}
\end{center}
\end{figure}

\pf
As noted above, $\Fix(\si)$ consists of the graphs from $\T_n$ whose edge-labels (as in Figure \ref{f:act}) satisfy $A=B=C$ and $a=b=c$.  Such a graph corresponds to a tetrahedron of the form shown in Figure~\ref{f:T4} (left).  If such a tetrahedron exists, then from $n=3A+3a$ we must have $n\equiv0\Mod3$, and $a=\frac n3-A$.  Then we note that \ref{T2} is trivial, while \ref{T3} is equivalent to the assertion that there is an $(A,a,a)$-triangle; the latter is equivalent to $A<2a$, and hence by integrality to $A+1\leq 2a=\frac{2n}3-2A$.  Keeping in mind that $A$ is an integer, the latter gives
\begin{align}
\label{e:231_1} A&\leq\left\lfloor\frac{2n-3}9\right\rfloor.
\intertext{For \ref{T6}, we have $x_1=x_2=\frac A2$, $Y_1=\frac{3A^2}4$ and $Y_2=a^2-\frac{A^2}4$, and the inequality in \ref{T6} becomes ${\frac{A^4}4<3A^2(a^2-\frac{A^2}4)}$.  Rearranging, and remembering that $A,a>0$, this gives $A<\sqrt3a=\sqrt3(\frac n3-A)$, and so $A<\frac n{3+\sqrt3}$.  Since $A$ is an integer, and since $\frac n{3+\sqrt3}$ is irrational, it follows that}
\label{e:231_2} A&\leq\left\lfloor\frac n{3+\sqrt3}\right\rfloor.
\end{align}
One may show that $\frac{2n-3}9\geq\frac n{3+\sqrt3} \iff n\geq31$, and also that $\left\lfloor\frac{2n-3}9\right\rfloor=\big\lfloor\frac n{3+\sqrt3}\big\rfloor$ if $n\leq30$ is a multiple of $3$.  The result then follows from \eqref{e:231_1} and \eqref{e:231_2}.
\epf

\begin{lemma}\label{l:fix_1243}
If $\si=(1,2,4,3)$, then
\[
\fix(\si) = \begin{cases}
0 &\text{if $n\not\equiv0\Mod2$}\\
\big\lfloor\frac n{4+4\sqrt2}\big\rfloor &\text{if $n\equiv0\Mod4$}\\
\big\lfloor\frac {n+2+2\sqrt2}{4+4\sqrt2}\big\rfloor &\text{if $n\equiv2\Mod4$.}
\end{cases}
\]
\end{lemma}

\pf
This time $\Fix(\si)$ consists of the graphs from $\T_n$ whose edge-labels (as in Figure \ref{f:act}) satisfy $A=a$ and $B=C=b=c$.  Such a graph corresponds to a tetrahedron of the form shown in Figure \ref{f:T4} (right).  If such a tetrahedron exists, then from $n=2A+4B$ we must have $n\equiv0\Mod2$, $A\equiv\frac n2\Mod2$ and $B=\frac{n-2A}4$.  Then we note that \ref{T2} and \ref{T3} are both equivalent to $A<2B$, and hence to $A+1\leq 2B=\frac n2-A$.  The latter gives
\begin{align}
\label{e:1243_1} A&\leq\left\lfloor\frac{n-2}4\right\rfloor.
\intertext{Rather than \ref{T6}, it is more convenient to work with the equivalent $|PQ|^2<a^2<|PR|^2$ from the proof of Proposition \ref{p:Tn}; cf.~Figure \ref{f:T3}.  Since $P=Q$ (and remembering that $A=a$, etc.), this is equivalent to $A^2<4B^2-A^2$, which gives $A<\sqrt2B=\frac{n-2A}{2\sqrt2}$, and ultimately}
\label{e:1243_2} A&\leq\left\lfloor\frac n{2+2\sqrt2}\right\rfloor.
\end{align}  
This time, $\frac{n-2}4\geq\frac n{2+2\sqrt2} \iff n\geq12$, and also $\left\lfloor\frac{n-2}4\right\rfloor=\big\lfloor\frac n{2+2\sqrt2}\big\rfloor$ if $n\leq10$ is even.  Thus, if we write $K=\big\lfloor\frac n{2+2\sqrt2}\big\rfloor$, then it follows from \eqref{e:1243_1} and \eqref{e:1243_2} that $\fix(\si)$ is equal to the cardinality of the set
\[
X = \bigset{A\in\N}{A\leq K,\ A\equiv\tfrac n2\Mod2}.
\]
Now,
\[
|X| = \begin{cases}
\lfloor\frac K2\rfloor &\text{if $\tfrac n2\equiv0\Mod2$: i.e., if $n\equiv0\Mod4$}\\
\lfloor\frac{K+1}2\rfloor &\text{if $\tfrac n2\equiv1\Mod2$: i.e., if $n\equiv2\Mod4$.}
\end{cases}
\]
The result then follows from the definition of $K$, and the fact that $\big\lfloor\frac{\lfloor x\rfloor}2\big\rfloor=\left\lfloor\frac x2\right\rfloor$ for real $x$.
\epf

\begin{figure}[ht]
\begin{center}
\includegraphics[height=5cm]{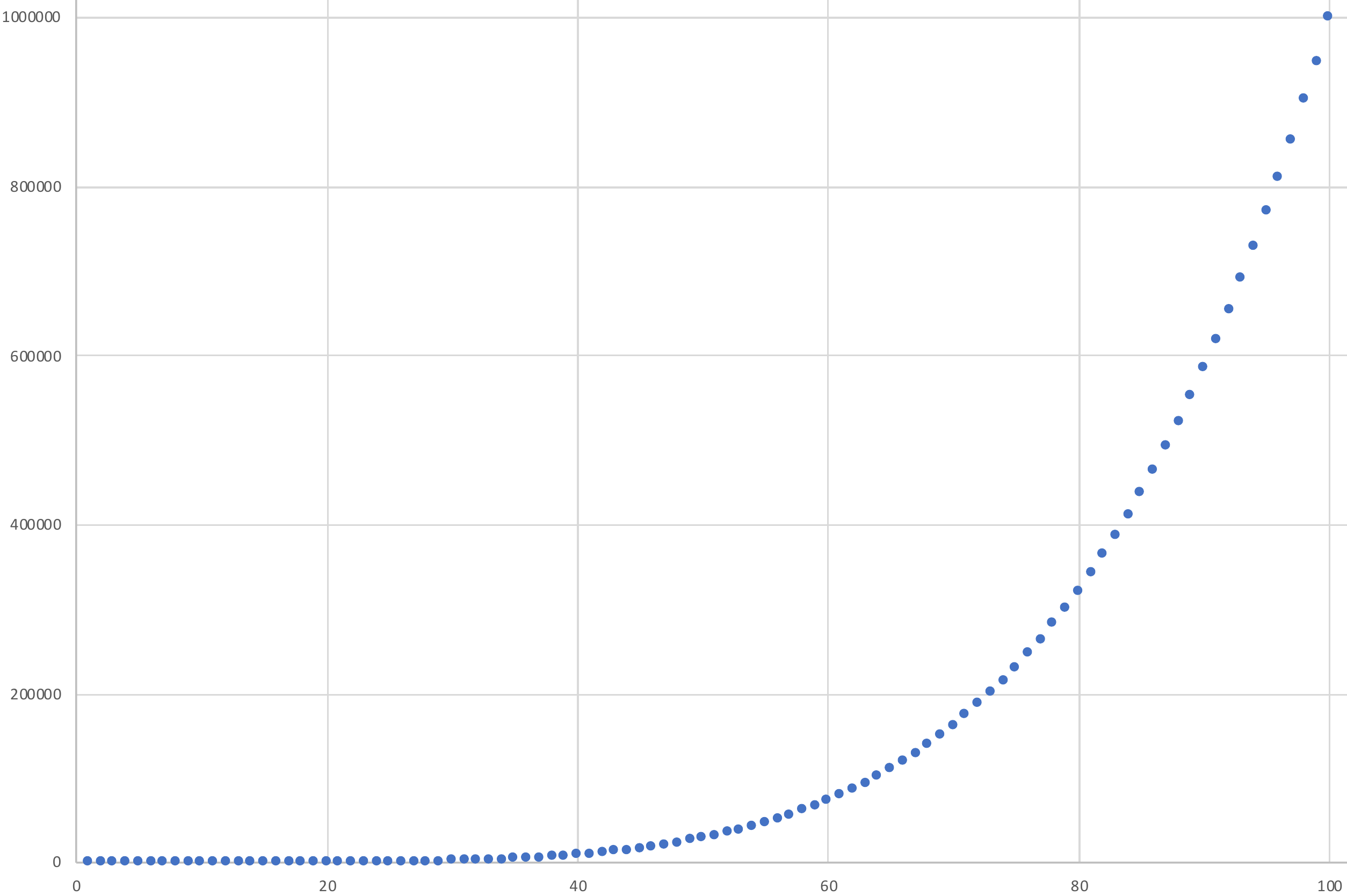}
\caption{Calculated values of $\fix(\si)$, for $\si$ of type \ref{fix1}, $1\leq n\leq100$.}
\label{f:i}
\end{center}
\end{figure}

\begin{figure}[ht]
\begin{center}
\includegraphics[height=4.25cm]{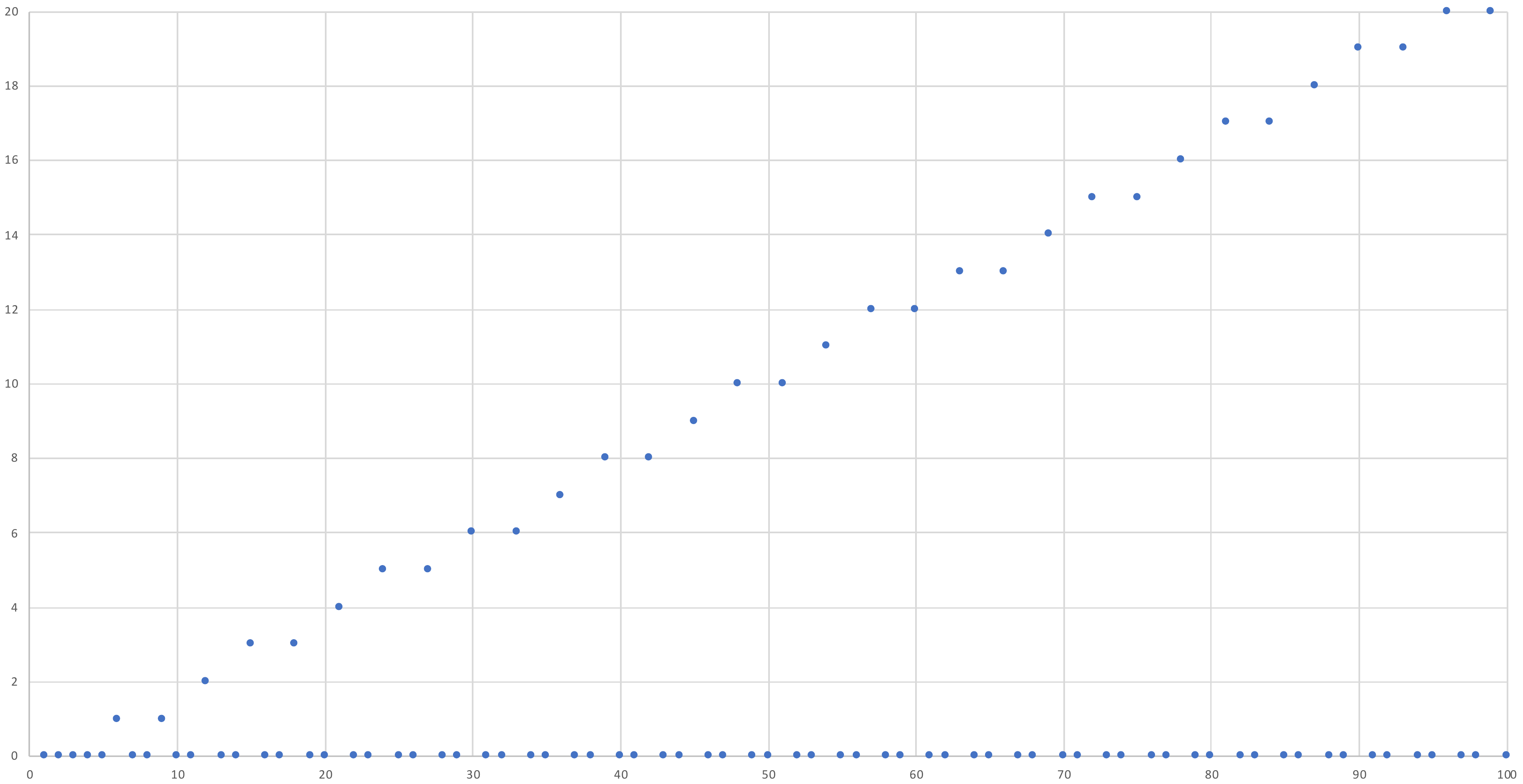}
\qquad
\includegraphics[height=4.25cm]{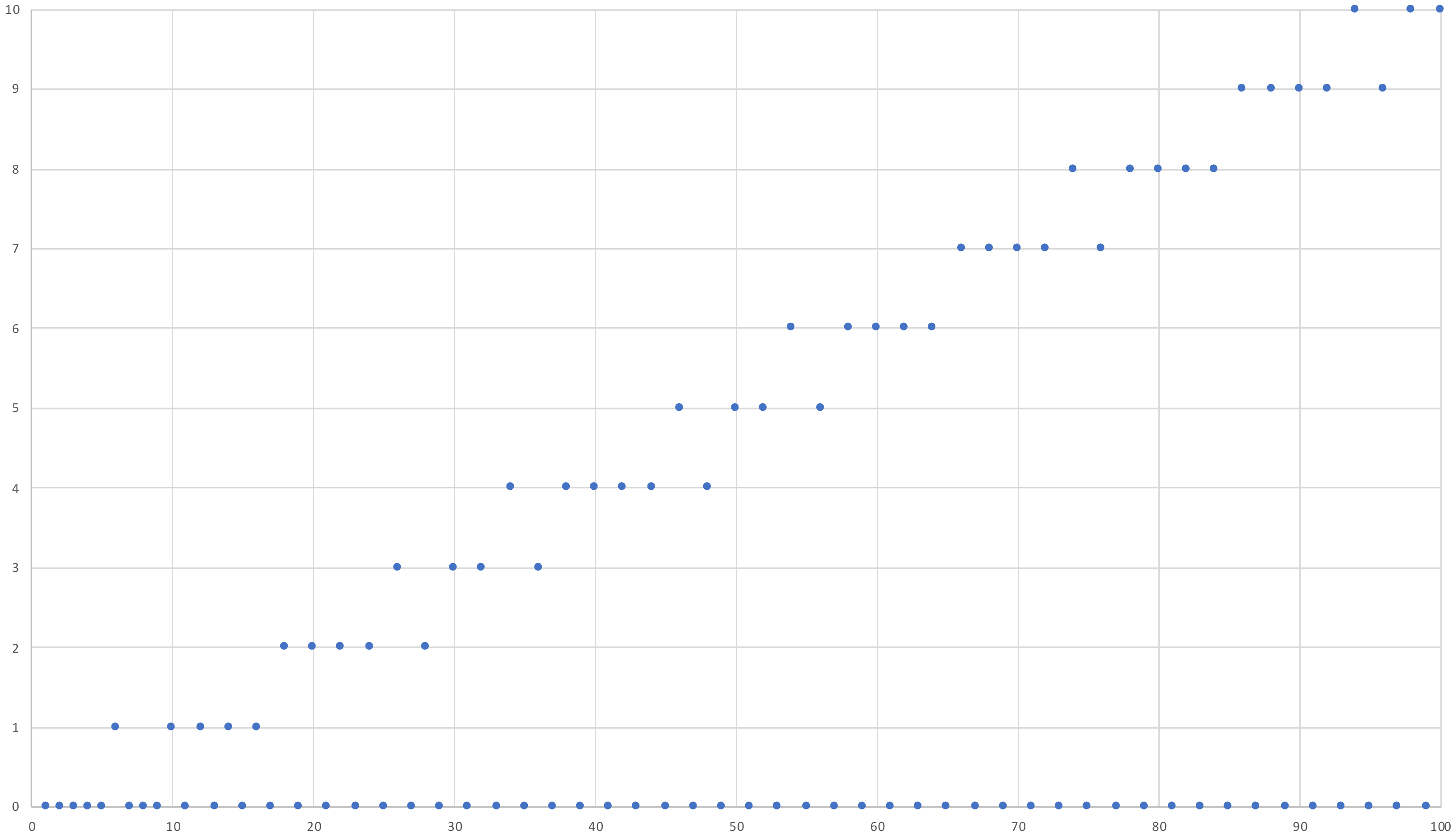}
\caption{Left and right: calculated values of $\fix(\si)$, for $\si$ of type \ref{fix2} and \ref{fix5}, $1\leq n\leq100$; cf.~Lemmas~\ref{l:fix_123} and~\ref{l:fix_1243}.}
\label{f:ii_v}
\end{center}
\end{figure}

\begin{figure}[ht]
\begin{center}
\includegraphics[height=4.8cm]{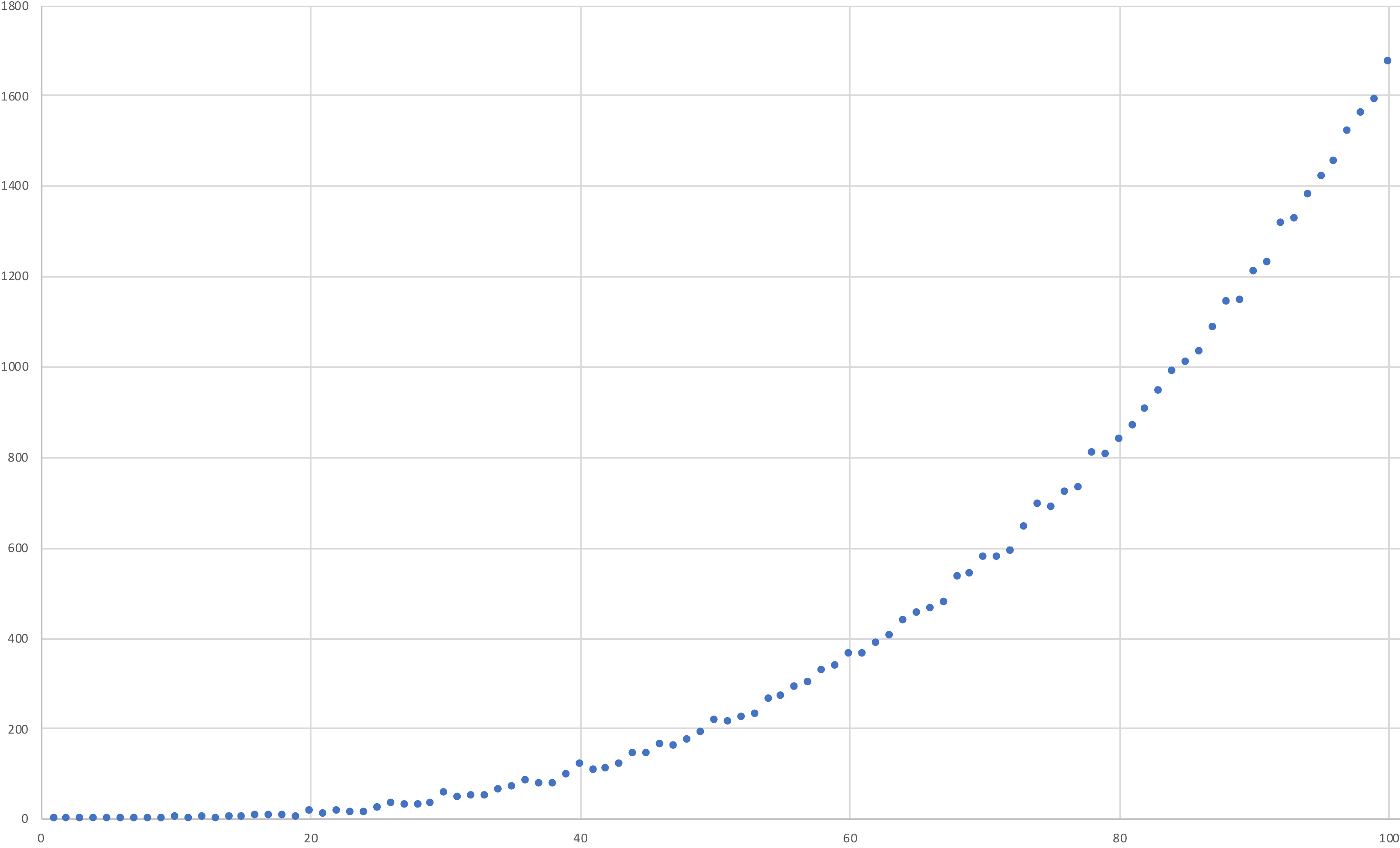}
\qquad
\includegraphics[height=4.8cm]{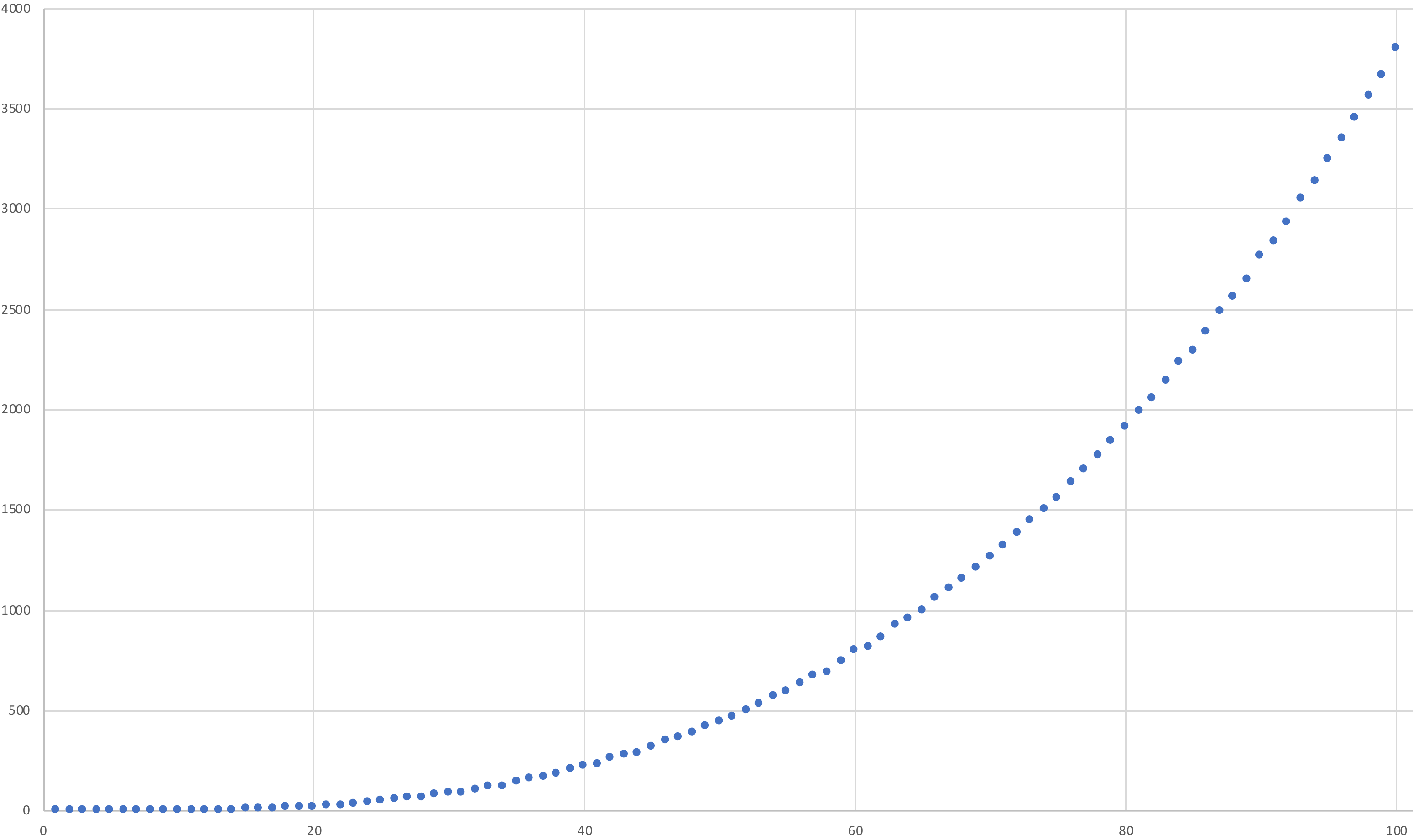}
\caption{Left and right: calculated values of $\fix(\si)$, for $\si$ of type \ref{fix3} and \ref{fix4}, $1\leq n\leq100$.}
\label{f:iii_iv}
\end{center}
\end{figure}

As noted above, we are currently unable to calculate $\fix(\si)$ for $\si$ of types \ref{fix1}, \ref{fix3} and~\ref{fix4}.  It is possible that types \ref{fix3} and \ref{fix4} could be treated by more elaborate versions of the arguments given in Lemmas \ref{l:fix_123} and \ref{l:fix_1243}.  Calculated values of all $\fix(\si)$ parameters for $1\leq n\leq100$ are given in Table~\ref{t:fix}; see also Figures \ref{f:i}--\ref{f:iii_iv}.
If one could calculate $\fix(\id_4)=|\T_n|$, this would have important consequences for the asymptotics of $\tn$ itself; see Section \ref{ss:as}.

\begin{rem}\label{r:td_tnd}
The arguments of this section can be modified to calculate the numbers $\td$ or $\tnd$.  Consider the following condition on a graph $G$ from $\G$, as shown in Figure \ref{f:act} (centre):
\begin{enumerate}[label=\textup{(T\arabic*)$'$},leftmargin=12mm]
\item \label{T1'} $\max\{A,B,C,a,b,c\}=d$.
\een
Then the set $\Td$ of graphs from $\G$ corresponding to integer tetrahedra of diameter $d$ are precisely those satisfying \ref{T1'} and \ref{T2}--\ref{T6}; cf.~Proposition \ref{p:Tn}.  The set $\Tnd=\Tn\cap\Td$ of graphs from $\G$ corresponding to integer tetrahedra of perimeter $n$ and diameter $d$ are precisely those satisfying \ref{T1'} and \ref{T1}--\ref{T6}.  Then
\[
\td = |\Td/\S_4| \AND \tnd = |\Tnd/\S_4|
\]
are given by counting orbits of the action of $\S_4$ given in \eqref{e:act}.  We also have the numbers $\td'$ and $\tnd'$, counting tetrahedra (of relevant parameters) up to rotations only, and these are given in terms of the restricted action of $\A_4$.
\end{rem}

\section{Computation and data}\label{s:computed}

In the previous section we gave a method for (in principle) computing the numbers $\tn$, $\td$ and $\tnd$.  Given that we are currently unable to give explicit formulas for the $\fix(\si)$ parameters in certain cases, and hence for the $\tn$, $\td$ and $\tnd$ sequences themselves, we now turn to some computations to make further progress.
In this section we discuss algorithms/code, and give several tables and graphs of calculated values.  In the next section we make some observations and conjectures based on the computational data.

\subsection{Enumeration algorithms}\label{ss:code}

To compute or enumerate all tetrahedra of a given perimeter $n$, we identify the graph $G$ from Figure~\ref{f:act} with the 6-tuple $\tup AaBbCc$.  This is clearly a tuple over $\{1,\ldots,n\}$, though we will see in Lemma \ref{l:d} that the entries of the tuple belong to a much smaller range; at this point, it is obvious at least that no entry could be bigger than $n-5$.  The orbits of $G$ (as in Figure \ref{f:act}) under the action of $\S_4$ given in \eqref{e:act} correspond to the following tuples:
\begin{align}
\nonumber &\tup AaBbCc, &&\tup AabBcC, &&\tup aABbcC, &&\tup aAbBCc, \\
\nonumber &\tup CcAaBb, &&\tup cCAabB, &&\tup cCaABb, &&\tup CcaAbB, \\
\nonumber &\tup BbCcAa, &&\tup bBcCAa, &&\tup BbcCaA, &&\tup bBCcaA, \\
\nonumber &\tup AaCcBb, &&\tup AacCbB, &&\tup aAcCBb, &&\tup aACcbB, \\
\nonumber &\tup BbAaCc, &&\tup bBAacC, &&\tup BbaAcC, &&\tup bBaACc, \\
\nonumber &\tup CcBbAa, &&\tup cCbBAa, &&\tup cCBbaA, &&\tup CcbBaA. \\
\label{e:tuples}
\end{align}
Each of these can be alternatively obtained from $\tup AaBbCc$ by the action of a signed permutation $\tau$ of the set $\{1,2,3\}$.  Roughly speaking, the non-signed part of $\tau$ tells us how to move the pairs $(A,a)$, $(B,b)$ and $(C,c)$, and the signed part tells us whether $(X,x)$ is listed in the final arrangement as $(X,x)$ or is flipped to $(x,X)$.  Moreover, examining \eqref{e:tuples}, only those signed permutations with an even number of flips occur.  This means that we are dealing with an action of the Coxeter group of type $D_3$, which is well-known to be isomorphic to the symmetric group $\S_4$ (type $A_3$).  For more on Coxeter groups, see \cite{Humphreys1990}.

The most obvious algorithm for creating the set $\Tn$ would be to take all 6-tuples $\tup AaBbCc$ over $\{1,\ldots,n-5\}$, and keep those satisfying conditions~\ref{T1}--\ref{T6} from Proposition \ref{p:Tn}.  Once we have the set $\Tn$, we can quickly compute the tuples fixed under the action of permutations $\si$ from $\S_4$ of types \ref{fix1}--\ref{fix5}, and thereby compute $\tn=\frac1{24}\sum_\si\fix(\si)$.

On the other hand, since we are ultimately concerned with counting \emph{orbits}, we do not need to store \emph{all} tuples from $\T_n$; we only need to store a single representative: e.g., the lex-greatest such tuple.

Some simple geometrical considerations allow us to reduce the search space, and also lead to other shortcuts.

\begin{lemma}\label{l:d}
If an integer tetrahedron has perimeter $n$ and diameter $d$, then
\[
\Big\lceil\frac n6\Big\rceil \leq d \leq \Big\lfloor\frac{n-3}3\Big\rfloor \qquad\text{or equivalently}\qquad 3d+3\leq n\leq 6d.
\]
\end{lemma}

\pf
It suffices to prove the second system of inequalities, and we note that $n\leq6d$ is clear.

Consider the graph $G$ pictured in Figure \ref{f:act}, and assume by symmetry that $d=A$.  From the $(A,B,C)$- and $(A,b,c)$-triangles we obtain $B+C\geq A+1$ and $b+c\geq A+1$, and so
\[
n = A+(B+C)+a+(b+c) \geq A+(A+1)+1+(A+1) = 3d+3.  \qedhere
\]
\epf

\begin{rem}
For fixed $d$, the maximum value of $n=6d$ occurs uniquely, of course, for the equilateral tetrahedron $\tup dddddd$.  The minimum value of $n=3d+3$ corresponds to $\tup 1d1d1d$, for example, a very ``thin'' tetrahedron of the form shown in Figure \ref{f:T4} (left).  In general, other tetrahedra have $n=3d+3$: e.g., $\tup1d2{d-1}2{d-1}$ for $d\geq3$.
\end{rem}

\begin{rem}
Although every edge of an integer teterahedron of perimeter $n$ is bounded above by $\lfloor\frac{n-3}3\rfloor$, it is not true that $\lceil\frac n6\rceil$ is a lower bound for all edges.  Indeed, we have the very ``thin'' tetrahedra $\tup 1d1d1d$ and $\tup dddd11$, as in Figure \ref{f:T4}.
\end{rem}

Similar considerations apply to the triangular faces:

\begin{lemma}\label{l:M}
If an integer tetrahedron has perimeter $n$ and maximum face-perimeter $M$, then
\[
\Big\lceil\frac n2\Big\rceil\leq M \leq \Big\lfloor\frac{2n-3}3 \Big\rfloor \qquad\text{or equivalently}\qquad \Big\lceil\frac {3M+3}2\Big\rceil\leq n\leq 2M.
\]
\end{lemma}

\pf
It suffices to show that $\frac {3M+3}2\leq n\leq 2M$.  To do so, consider the graph $G$ pictured in Figure~\ref{f:act}, assuming by symmetry that $M=A+B+C$.  Adding the perimeters of all four faces gives twice the perimeter of the tetrahedron.  The first consequence of this is that $2n\leq4M$: i.e., $n\leq2M$.  The second consequence (also using $b+c\geq A+1$, etc., from the other triangular faces) is that
\begin{align*}
2n &= (A+B+C) + (A+b+c) + (B+a+c) + (C+a+b) \\
&\geq (A+B+C) + (A+A+1) + (B+B+1) + (C+C+1) = 3(A+B+C) + 3 = 3M+3,
\end{align*}
which gives $n\geq\frac{3M+3}3$.
\epf

With the above considerations in mind, a simple algorithm for calculating a set of representatives of $\Tnd$ is as follows.  Here we write $d=A$ to keep the notation as in Figure \ref{f:act}.  Roughly speaking, we create tetrahedra as in Figure \ref{f:T3}, with the diameter on the $x$-axis.
\BEN
\item \label{a:I} Define the set $S=\emptyset$.
\item \label{a:II} Create the set of all triangles with maximum side $A({}=d)$ and with perimeter at most $\left\lfloor\tfrac{2n-3}3 \right\rfloor$.  (This is fairly routine, so the details are omitted.)  
\item \label{a:III} Then for each pair $(A, B, C)$ and $(A, b, c)$ of such triangles, we set $a={n-(A+B+C+b+c)}$.  
\item \label{a:IV} If $a\leq A$ and if the tuple $T=\tup AaBbCc$ satisfies condition \ref{T6} from Proposition \ref{p:Tn}, we add to the set $S$ the lex-greatest representative of $T$ from the list \eqref{e:tuples}.  
\een
The number $\tnd={}^At_n$ is then the size of the set $S$ created.  Summing over all $d$ (or $n$) gives~$\tn$ (or~$\td$), respectively.  There are many easy ways to simplify the above.  For example, the sixth entry of the tuple $\tup AaBbCc$ does not need to be stored, as it can be inferred from the others; neither do we need to store the first entry as this is always $d$ (this is why we chose to store the lex-\emph{greatest} representative).  

We carried out the above algorithm in a variety of languages.  We were able to calculate $\tn$ up to $n=1100$ before running out of memory on standard laptops (even for single values of $d$).  It was around this stage in our investigations that we discovered Kurz's article \cite{Kurz2007}.  We regret to say that the ``inte\emph{gral}'' rather than ``inte\emph{ger}'' in his title defeated our searching abilities.  Additionally, since we were initially interested only in the $\tn$ sequence, we initially found only Sequence A208454 on the OEIS \cite{OEIS}, but not A097125, which contains many terms in the $\td$ sequence.  

The most significant advantage of Kurz's algorithm for enumeration by diameter is that it does not require the creation and storage of vast numbers of tuples.  Rather, it moves through ``tuple space'' in such an orderly fashion that each canonical tuple is visited once, meaning that these only have to be counted rather than stored.  Our current algorithms are modifications of Kurz's, implemented in
C, and with a number of optimisations, including parallelisation and use of the GNU Multiple Precision Arithmetic Library to
store and manipulate large integers.  All of our calculations of the numbers~$\td$ match those given by Kurz, though we have gone somewhat further; see Table \ref{t:rhod}.

\subsection{Values}\label{ss:values}

So far we have calculated $\tn$ for every $1\leq n\leq3000$, and for several other values up to $n=20000$.  Various tables and figures in this section summarise some of these calculations, but more complete lists are available at \cite{web}.  Some of the tables below also include values of $\td$, $\tnd$ and $\tn'$ (tetrahedra up to rotations and translations only), and also $\fix(\si)$ for certain permutations $\si\in\S_4$.  (More extensive tables containing values of $\td$ are given by Kurz \cite{Kurz2007}.)  For convenience of reference, here is a summary of the content of the tables and figures given in this section (more are given in Section \ref{s:obs}, and we note that some of the tables in the current section also include additional parameters defined in Section \ref{s:obs}):
\begin{center}
\begin{tabular}{|l|c|l|c|}
\hline
\multicolumn{1}{|c|}{Location} & Numbers & \multicolumn{1}{|c|}{Range} & Comments \\
\hline
Table \ref{t:tn_200} & $\tn$ & $1\leq n\leq 200$ & \\
Table \ref{t:tn_10000} & $\tn$ & $1\leq n\leq 10000$ & by hundreds \\
Table \ref{t:tnd_50} & $\tnd$ & $1\leq n\leq 50$ & \\
Table \ref{t:tnd_200} & $\tnd$ & $1\leq n\leq 200$ & \\
Table \ref{t:fix} & $\fix(\si)$, $\tn$, $\tn'$ & $1\leq n\leq 100$ & \\
Table \ref{t:rhon} & $\tn$ & $1\leq n\leq 20000$ & by thousands \\
Table \ref{t:rhod} & $\td$ & $1\leq d\leq 2000$ & by hundreds\\
\hline
Figure \ref{f:tn_30} & $\tn$ & $1\leq n\leq 30$ & \\
Figure \ref{f:tn} & $\tn$ & $1\leq n\leq 200$ & \\
Figure \ref{f:tn_log} & $\tn$ & $1\leq n\leq 200$ & log-log plot \\
Figure \ref{f:tn_100s} & $\tn$ & $1\leq n\leq 10000$ & by hundreds \\
Figure \ref{f:HM} & $\tnd$ & $1\leq n\leq 2000$ & heat map \\
Figure \ref{f:tnd_100_200} & $\tnd$ & $n=100$ and $200$ & \\
Figure \ref{f:tnd_10000} & $\tnd$ & $n=10000$ & and discrete derivative \\
Figure \ref{f:tnd_50_100} & $\tnd$ & $d=50$ and $100$ & \\
Figure \ref{f:tnd_1000} & $\tnd$ & $d=1000$ & and discrete derivative \\
Figures \ref{f:i}--\ref{f:iii_iv} & $\fix(\si)$ & $1\leq n\leq 100$ & $\si\in\S_4$ \\
\hline
\end{tabular}
\end{center}

\begin{table}[ht]
\begin{center}
\begin{tabular}{|c|c|}
\hline
$n$ & $\tn$ \\
\hline
$	1	$ & $	0	$ \\
$	2	$ & $	0	$ \\
$	3	$ & $	0	$ \\
$	4	$ & $	0	$ \\
$	5	$ & $	0	$ \\
$	6	$ & $	1	$ \\
$	7	$ & $	0	$ \\
$	8	$ & $	0	$ \\
$	9	$ & $	1	$ \\
$	10	$ & $	1	$ \\
\hline
$	11	$ & $	1	$ \\
$	12	$ & $	3	$ \\
$	13	$ & $	2	$ \\
$	14	$ & $	3	$ \\
$	15	$ & $	6	$ \\
$	16	$ & $	6	$ \\
$	17	$ & $	7	$ \\
$	18	$ & $	12	$ \\
$	19	$ & $	11	$ \\
$	20	$ & $	18	$ \\
\hline
$	21	$ & $	21	$ \\
$	22	$ & $	25	$ \\
$	23	$ & $	31	$ \\
$	24	$ & $	38	$ \\
$	25	$ & $	46	$ \\
$	26	$ & $	56	$ \\
$	27	$ & $	66	$ \\
$	28	$ & $	76	$ \\
$	29	$ & $	90	$ \\
$	30	$ & $	117	$ \\
\hline
$	31	$ & $	123	$ \\
$	32	$ & $	151	$ \\
$	33	$ & $	175	$ \\
$	34	$ & $	196	$ \\
$	35	$ & $	234	$ \\
$	36	$ & $	264	$ \\
$	37	$ & $	297	$ \\
$	38	$ & $	346	$ \\
$	39	$ & $	391	$ \\
$	40	$ & $	448	$ \\
\hline
\end{tabular}
\qquad
\begin{tabular}{|c|c|}
\hline
$n$ & $\tn$ \\
\hline
$	41	$ & $	492	$ \\
$	42	$ & $	568	$ \\
$	43	$ & $	630	$ \\
$	44	$ & $	702	$ \\
$	45	$ & $	797	$ \\
$	46	$ & $	884	$ \\
$	47	$ & $	977	$ \\
$	48	$ & $	1089	$ \\
$	49	$ & $	1217	$ \\
$	50	$ & $	1338	$ \\
\hline
$	51	$ & $	1469	$ \\
$	52	$ & $	1624	$ \\
$	53	$ & $	1771	$ \\
$	54	$ & $	1970	$ \\
$	55	$ & $	2146	$ \\
$	56	$ & $	2343	$ \\
$	57	$ & $	2579	$ \\
$	58	$ & $	2782	$ \\
$	59	$ & $	3042	$ \\
$	60	$ & $	3322	$ \\
\hline
$	61	$ & $	3586	$ \\
$	62	$ & $	3912	$ \\
$	63	$ & $	4221	$ \\
$	64	$ & $	4568	$ \\
$	65	$ & $	4953	$ \\
$	66	$ & $	5339	$ \\
$	67	$ & $	5731	$ \\
$	68	$ & $	6204	$ \\
$	69	$ & $	6657	$ \\
$	70	$ & $	7169	$ \\
\hline
$	71	$ & $	7683	$ \\
$	72	$ & $	8230	$ \\
$	73	$ & $	8857	$ \\
$	74	$ & $	9446	$ \\
$	75	$ & $	10095	$ \\
$	76	$ & $	10846	$ \\
$	77	$ & $	11513	$ \\
$	78	$ & $	12345	$ \\
$	79	$ & $	13125	$ \\
$	80	$ & $	13969	$ \\
\hline
\end{tabular}
\qquad
\begin{tabular}{|c|c|}
\hline
$n$ & $\tn$ \\
\hline
$	81	$ & $	14903	$ \\
$	82	$ & $	15790	$ \\
$	83	$ & $	16811	$ \\
$	84	$ & $	17854	$ \\
$	85	$ & $	18940	$ \\
$	86	$ & $	20107	$ \\
$	87	$ & $	21261	$ \\
$	88	$ & $	22511	$ \\
$	89	$ & $	23831	$ \\
$	90	$ & $	25251	$ \\
\hline
$	91	$ & $	26631	$ \\
$	92	$ & $	28173	$ \\
$	93	$ & $	29744	$ \\
$	94	$ & $	31341	$ \\
$	95	$ & $	33116	$ \\
$	96	$ & $	34849	$ \\
$	97	$ & $	36696	$ \\
$	98	$ & $	38695	$ \\
$	99	$ & $	40619	$ \\
$	100	$ & $	42817	$ \\
\hline
$	101	$ & $	44951	$ \\
$	102	$ & $	47257	$ \\
$	103	$ & $	49641	$ \\
$	104	$ & $	52008	$ \\
$	105	$ & $	54685	$ \\
$	106	$ & $	57310	$ \\
$	107	$ & $	60046	$ \\
$	108	$ & $	62944	$ \\
$	109	$ & $	65896	$ \\
$	110	$ & $	69029	$ \\
\hline
$	111	$ & $	72152	$ \\
$	112	$ & $	75549	$ \\
$	113	$ & $	78976	$ \\
$	114	$ & $	82574	$ \\
$	115	$ & $	86199	$ \\
$	116	$ & $	90052	$ \\
$	117	$ & $	94053	$ \\
$	118	$ & $	98089	$ \\
$	119	$ & $	102371	$ \\
$	120	$ & $	106711	$ \\
\hline
\end{tabular}
\qquad
\begin{tabular}{|c|c|}
\hline
$n$ & $\tn$ \\
\hline
$	121	$ & $	111246	$ \\
$	122	$ & $	116006	$ \\
$	123	$ & $	120768	$ \\
$	124	$ & $	125837	$ \\
$	125	$ & $	130985	$ \\
$	126	$ & $	136402	$ \\
$	127	$ & $	141853	$ \\
$	128	$ & $	147575	$ \\
$	129	$ & $	153500	$ \\
$	130	$ & $	159470	$ \\
\hline
$	131	$ & $	165757	$ \\
$	132	$ & $	172253	$ \\
$	133	$ & $	178772	$ \\
$	134	$ & $	185682	$ \\
$	135	$ & $	192674	$ \\
$	136	$ & $	200070	$ \\
$	137	$ & $	207440	$ \\
$	138	$ & $	215111	$ \\
$	139	$ & $	223151	$ \\
$	140	$ & $	231179	$ \\
\hline
$	141	$ & $	239687	$ \\
$	142	$ & $	248292	$ \\
$	143	$ & $	257130	$ \\
$	144	$ & $	266522	$ \\
$	145	$ & $	275691	$ \\
$	146	$ & $	285536	$ \\
$	147	$ & $	295385	$ \\
$	148	$ & $	305605	$ \\
$	149	$ & $	316146	$ \\
$	150	$ & $	326952	$ \\
\hline
$	151	$ & $	337891	$ \\
$	152	$ & $	349316	$ \\
$	153	$ & $	361051	$ \\
$	154	$ & $	372976	$ \\
$	155	$ & $	385325	$ \\
$	156	$ & $	398032	$ \\
$	157	$ & $	410750	$ \\
$	158	$ & $	424249	$ \\
$	159	$ & $	437744	$ \\
$	160	$ & $	451759	$ \\
\hline
\end{tabular}
\qquad
\begin{tabular}{|c|c|}
\hline
$n$ & $\tn$ \\
\hline
$	161	$ & $	466128	$ \\
$	162	$ & $	480898	$ \\
$	163	$ & $	495924	$ \\
$	164	$ & $	511275	$ \\
$	165	$ & $	527179	$ \\
$	166	$ & $	543486	$ \\
$	167	$ & $	559950	$ \\
$	168	$ & $	576987	$ \\
$	169	$ & $	594453	$ \\
$	170	$ & $	612306	$ \\
\hline
$	171	$ & $	630619	$ \\
$	172	$ & $	649344	$ \\
$	173	$ & $	668273	$ \\
$	174	$ & $	688082	$ \\
$	175	$ & $	708017	$ \\
$	176	$ & $	728557	$ \\
$	177	$ & $	749759	$ \\
$	178	$ & $	770893	$ \\
$	179	$ & $	792983	$ \\
$	180	$ & $	815666	$ \\
\hline
$	181	$ & $	838151	$ \\
$	182	$ & $	862084	$ \\
$	183	$ & $	885874	$ \\
$	184	$ & $	910270	$ \\
$	185	$ & $	935691	$ \\
$	186	$ & $	960959	$ \\
$	187	$ & $	987344	$ \\
$	188	$ & $	1014149	$ \\
$	189	$ & $	1041441	$ \\
$	190	$ & $	1069393	$ \\
\hline
$	191	$ & $	1097727	$ \\
$	192	$ & $	1127002	$ \\
$	193	$ & $	1156747	$ \\
$	194	$ & $	1186965	$ \\
$	195	$ & $	1217963	$ \\
$	196	$ & $	1249657	$ \\
$	197	$ & $	1281804	$ \\
$	198	$ & $	1314840	$ \\
$	199	$ & $	1348514	$ \\
$	200	$ & $	1382630	$ \\
\hline
\end{tabular}
\caption{Calculated values of $\tn$.}
\label{t:tn_200}
\end{center}
\end{table}

\begin{table}[ht]
\begin{center}
\rotatebox{90}{
\begin{tabular}{|c|c|}
\hline
$n$ & $\tn$ \\
\hline
$	100	$ & $	42817	$ \\
$	200	$ & $	1382630	$ \\
$	300	$ & $	10542791	$ \\
$	400	$ & $	44512930	$ \\
$	500	$ & $	135995968	$ \\
$	600	$ & $	338647149	$ \\
$	700	$ & $	732301457	$ \\
$	800	$ & $	1428244253	$ \\
$	900	$ & $	2574424362	$ \\
$	1000	$ & $	4360687860	$ \\
\hline
$	1100	$ & $	7024066455	$ \\
$	1200	$ & $	10853927088	$ \\
$	1300	$ & $	16197354205	$ \\
$	1400	$ & $	23464200340	$ \\
$	1500	$ & $	33132516713	$ \\
$	1600	$ & $	45753651322	$ \\
$	1700	$ & $	61957547844	$ \\
$	1800	$ & $	82458026396	$ \\
$	1900	$ & $	108057877719	$ \\
$	2000	$ & $	139654346301	$ \\
\hline
\end{tabular}
\qquad
\begin{tabular}{|c|c|}
\hline
$n$ & $\tn$ \\
\hline
$	2100	$ & $	178244109766	$ \\
$	2200	$ & $	224928623704	$ \\
$	2300	$ & $	280919515766	$ \\
$	2400	$ & $	347543441363	$ \\
$	2500	$ & $	426247856357	$ \\
$	2600	$ & $	518605683615	$ \\
$	2700	$ & $	626321040933	$ \\
$	2800	$ & $	751234199645	$ \\
$	2900	$ & $	895326815879	$ \\
$	3000	$ & $	1060727392377	$ \\
\hline
$	3100	$ & $	1249716260424	$ \\
$	3200	$ & $	1464730980692	$ \\
$	3300	$ & $	1708371608956	$ \\
$	3400	$ & $	1983405792552	$ \\
$	3500	$ & $	2292773954230	$ \\
$	3600	$ & $	2639594896693	$ \\
$	3700	$ & $	3027170849616	$ \\
$	3800	$ & $	3458992238494	$ \\
$	3900	$ & $	3938743901688	$ \\
$	4000	$ & $	4470309160343	$ \\
\hline
\end{tabular}
\qquad
\begin{tabular}{|c|c|}
\hline
$n$ & $\tn$ \\
\hline
$	4100	$ & $	5057776483155	$ \\
$	4200	$ & $	5705443122022	$ \\
$	4300	$ & $	6417821623924	$ \\
$	4400	$ & $	7199644279518	$ \\
$	4500	$ & $	8055868989730	$ \\
$	4600	$ & $	8991684402358	$ \\
$	4700	$ & $	10012513881131	$ \\
$	4800	$ & $	11124023171269	$ \\
$	4900	$ & $	12332123259033	$ \\
$	5000	$ & $	13642977397892	$ \\
\hline
$	5100	$ & $	15063004902111	$ \\
$	5200	$ & $	16598887173873	$ \\
$	5300	$ & $	18257573582602	$ \\
$	5400	$ & $	20046284952842	$ \\
$	5500	$ & $	21972520100932	$ \\
$	5600	$ & $	24044061595649	$ \\
$	5700	$ & $	26268978527115	$ \\
$	5800	$ & $	28655635295636	$ \\
$	5900	$ & $	31212693377826	$ \\
$	6000	$ & $	33949118928429	$ \\
\hline
\end{tabular}
\qquad
\begin{tabular}{|c|c|}
\hline
$n$ & $\tn$ \\
\hline
$	6100	$ & $	36874187631384	$ \\
$	6200	$ & $	39997488181501	$ \\
$	6300	$ & $	43328930438429	$ \\
$	6400	$ & $	46878747743940	$ \\
$	6500	$ & $	50657504819100	$ \\
$	6600	$ & $	54676100750626	$ \\
$	6700	$ & $	58945775513454	$ \\
$	6800	$ & $	63478115036529	$ \\
$	6900	$ & $	68285056161851	$ \\
$	7000	$ & $	73378891579018	$ \\
\hline
$	7100	$ & $	78772276600831	$ \\
$	7200	$ & $	84478231957643	$ \\
$	7300	$ & $	90510152207309	$ \\
$	7400	$ & $	96881807240927	$ \\
$	7500	$ & $	103607350139966	$ \\
$	7600	$ & $	110701321888040	$ \\
$	7700	$ & $	118178657158659	$ \\
$	7800	$ & $	126054687371188	$ \\
$	7900	$ & $	134345147515787	$ \\
$	8000	$ & $	143066182224551	$ \\
\hline
\end{tabular}
\qquad
\begin{tabular}{|c|c|}
\hline
$n$ & $\tn$ \\
\hline
$	8100	$ & $	152234350689174	$ \\
$	8200	$ & $	161866628410449	$ \\
$	8300	$ & $	171980417994570	$ \\
$	8400	$ & $	182593551239230	$ \\
$	8500	$ & $	193724293506541	$ \\
$	8600	$ & $	205391352861569	$ \\
$	8700	$ & $	217613878464370	$ \\
$	8800	$ & $	230411473926234	$ \\
$	8900	$ & $	243804197732022	$ \\
$	9000	$ & $	257812568218126	$ \\
\hline
$	9100	$ & $	272457573348976	$ \\
$	9200	$ & $	287760666512263	$ \\
$	9300	$ & $	303743784383639	$ \\
$	9400	$ & $	320429342544289	$ \\
$	9500	$ & $	337840242823660	$ \\
$	9600	$ & $	355999882739690	$ \\
$	9700	$ & $	374932153990010	$ \\
$	9800	$ & $	394661455344190	$ \\
$	9900	$ & $	415212690326932	$ \\
$	10000	$ & $	436611276762080	$ \\
\hline
\end{tabular}
}
\caption{Calculated values of $\tn$.}
\label{t:tn_10000}
\end{center}
\end{table}

\begin{table}[ht]
\begin{center}
{\small
\begin{tabular}{c|ccccccccccccccc}
$	n/d	$ & $	1	$ & $	2	$ & $	3	$ & $	4	$ & $	5	$ & $	6	$ & $	7	$ & $	8	$ & $	9	$ & $	10	$ & $	11	$ & $	12	$ & $	13	$ & $	14	$ & $	15	$ \\
\hline
$	6	$ & $		{\red1}	$ & $	.	$ & $	.	$ & $	.	$ & $	.	$ & $	.	$ & $	.	$ & $	.	$ & $	.	$ & $	.	$ & $	.	$ & $	.	$ & $	.	$ & $	.	$ & $	.	$ \\
$	7	$ & $.	$ & $	.	$ & $	.	$ & $	.	$ & $	.	$ & $	.	$ & $	.	$ & $	.	$ & $	.	$ & $	.	$ & $	.	$ & $	.	$ & $	.	$ & $	.	$ & $	.	$ \\
$	8	$ & $.	$ & $	.	$ & $	.	$ & $	.	$ & $	.	$ & $	.	$ & $	.	$ & $	.	$ & $	.	$ & $	.	$ & $	.	$ & $	.	$ & $	.	$ & $	.	$ & $	.	$ \\
$	9	$ & $.	$ & $	{\blue1}	$ & $	.	$ & $	.	$ & $	.	$ & $	.	$ & $	.	$ & $	.	$ & $	.	$ & $	.	$ & $	.	$ & $	.	$ & $	.	$ & $	.	$ & $	.	$ \\
$	10	$ & $.	$ & $\textcolor{brown}{1}	$ & $	.	$ & $	.	$ & $	.	$ & $	.	$ & $	.	$ & $	.	$ & $	.	$ & $	.	$ & $	.	$ & $	.	$ & $	.	$ & $	.	$ & $	.	$ \\
$	11	$ & $.	$ & $		{\red1}	$ & $	.	$ & $	.	$ & $	.	$ & $	.	$ & $	.	$ & $	.	$ & $	.	$ & $	.	$ & $	.	$ & $	.	$ & $	.	$ & $	.	$ & $	.	$ \\
$	12	$ & $.	$ & $		{\red1}	$ & $	{\blue2}	$ & $	.	$ & $	.	$ & $	.	$ & $	.	$ & $	.	$ & $	.	$ & $	.	$ & $	.	$ & $	.	$ & $	.	$ & $	.	$ & $	.	$ \\
$	13	$ & $.	$ & $.	$ & $\textcolor{brown}{2}	$ & $	.	$ & $	.	$ & $	.	$ & $	.	$ & $	.	$ & $	.	$ & $	.	$ & $	.	$ & $	.	$ & $	.	$ & $	.	$ & $	.	$ \\
$	14	$ & $.	$ & $.	$ & ${\green 3}	$ & $	.	$ & $	.	$ & $	.	$ & $	.	$ & $	.	$ & $	.	$ & $	.	$ & $	.	$ & $	.	$ & $	.	$ & $	.	$ & $	.	$ \\
$	15	$ & $.	$ & $.	$ & $	4	$ & $	{\blue2}	$ & $	.	$ & $	.	$ & $	.	$ & $	.	$ & $	.	$ & $	.	$ & $	.	$ & $	.	$ & $	.	$ & $	.	$ & $	.	$ \\
$	16	$ & $.	$ & $.	$ & $	{\red3}	$ & $\textcolor{brown}{3}	$ & $	.	$ & $	.	$ & $	.	$ & $	.	$ & $	.	$ & $	.	$ & $	.	$ & $	.	$ & $	.	$ & $	.	$ & $	.	$ \\
$	17	$ & $.	$ & $.	$ & $	{\red1}	$ & ${\green 6}	$ & $	.	$ & $	.	$ & $	.	$ & $	.	$ & $	.	$ & $	.	$ & $	.	$ & $	.	$ & $	.	$ & $	.	$ & $	.	$ \\
$	18	$ & $.	$ & $.	$ & $	{\red1}$ & $	8	$ & $	{\blue3}	$ & $	.	$ & $	.	$ & $	.	$ & $	.	$ & $	.	$ & $	.	$ & $	.	$ & $	.	$ & $	.	$ & $	.	$ \\
$	19	$ & $.	$ & $.	$ & $.	$ & $	7	$ & $\textcolor{brown}{4}	$ & $	.	$ & $	.	$ & $	.	$ & $	.	$ & $	.	$ & $	.	$ & $	.	$ & $	.	$ & $	.	$ & $	.	$ \\
$	20	$ & $.	$ & $.	$ & $.	$ & $	8	$ & ${\green 10}	$ & $	.	$ & $	.	$ & $	.	$ & $	.	$ & $	.	$ & $	.	$ & $	.	$ & $	.	$ & $	.	$ & $	.	$ \\
$	21	$ & $.	$ & $.	$ & $.	$ & $	{\red6}	$ & $	12	$ & $	{\blue3}	$ & $	.	$ & $	.	$ & $	.	$ & $	.	$ & $	.	$ & $	.	$ & $	.	$ & $	.	$ & $	.	$ \\
$	22	$ & $.	$ & $.	$ & $.	$ & $	{\red3}	$ & $	17	$ & $\textcolor{brown}{5}	$ & $	.	$ & $	.	$ & $	.	$ & $	.	$ & $	.	$ & $	.	$ & $	.	$ & $	.	$ & $	.	$ \\
$	23	$ & $.	$ & $.	$ & $.	$ & $	{\red1}	$ & $	17	$ & ${\green 13}	$ & $	.	$ & $	.	$ & $	.	$ & $	.	$ & $	.	$ & $	.	$ & $	.	$ & $	.	$ & $	.	$ \\
$	24	$ & $.	$ & $.	$ & $.	$ & $	{\red1}	$ & $	16	$ & $	17	$ & $	{\blue4}	$ & $	.	$ & $	.	$ & $	.	$ & $	.	$ & $	.	$ & $	.	$ & $	.	$ & $	.	$ \\
$	25	$ & $.	$ & $.	$ & $.	$ & $.	$ & $	15	$ & $	25	$ & $\textcolor{brown}{6}	$ & $	.	$ & $	.	$ & $	.	$ & $	.	$ & $	.	$ & $	.	$ & $	.	$ & $	.	$ \\
$	26	$ & $.	$ & $.	$ & $.	$ & $.	$ & $	{\red11}	$ & $	29	$ & ${\green 16}	$ & $	.	$ & $	.	$ & $	.	$ & $	.	$ & $	.	$ & $	.	$ & $	.	$ & $	.	$ \\
$	27	$ & $.	$ & $.	$ & $.	$ & $.	$ & ${\red	6}	$ & $	34	$ & $	22	$ & $	{\blue4}	$ & $	.	$ & $	.	$ & $	.	$ & $	.	$ & $	.	$ & $	.	$ & $	.	$ \\
$	28	$ & $.	$ & $.	$ & $.	$ & $.	$ & $	{\red3}	$ & $	32	$ & $	34	$ & $\textcolor{brown}{7}	$ & $	.	$ & $	.	$ & $	.	$ & $	.	$ & $	.	$ & $	.	$ & $	.	$ \\
$	29	$ & $.	$ & $.	$ & $.	$ & $.	$ & $	{\red1}	$ & $	29	$ & $	41	$ & ${\green 19}	$ & $	.	$ & $	.	$ & $	.	$ & $	.	$ & $	.	$ & $	.	$ & $	.	$ \\
$	30	$ & $.	$ & $.	$ & $.	$ & $.	$ & $	{\red1}	$ & $	27	$ & $	55	$ & $	29	$ & $	{\blue5}	$ & $	.	$ & $	.	$ & $	.	$ & $	.	$ & $	.	$ & $	.	$ \\
$	31	$ & $.	$ & $.	$ & $.	$ & $.	$ & $.	$ & $	{\red18}	$ & $	55	$ & $	42	$ & $\textcolor{brown}{8}	$ & $	.	$ & $	.	$ & $	.	$ & $	.	$ & $	.	$ & $	.	$ \\
$	32	$ & $.	$ & $.	$ & $.	$ & $.	$ & $.	$ & $	{\red11}	$ & $	61	$ & $	57	$ & ${\green 22}	$ & $	.	$ & $	.	$ & $	.	$ & $	.	$ & $	.	$ & $	.	$ \\
$	33	$ & $.	$ & $.	$ & $.	$ & $.	$ & $.	$ & $	{\red6}	$ & $	57	$ & $	72	$ & $	35	$ & $	{\blue5}	$ & $	.	$ & $	.	$ & $	.	$ & $	.	$ & $	.	$ \\
$	34	$ & $.	$ & $.	$ & $.	$ & $.	$ & $.	$ & ${\red	3}	$ & $	51	$ & $	81	$ & $	52	$ & $\textcolor{brown}{9}	$ & $	.	$ & $	.	$ & $	.	$ & $	.	$ & $	.	$ \\
$	35	$ & $.	$ & $.	$ & $.	$ & $.	$ & $.	$ & $	{\red1}	$ & $	43	$ & $	94	$ & $	71	$ & ${\green 25}	$ & $	.	$ & $	.	$ & $	.	$ & $	.	$ & $	.	$ \\
$	36	$ & $.	$ & $.	$ & $.	$ & $.	$ & $.	$ & $	{\red1}	$ & ${\red	31}	$ & $	95	$ & $	91	$ & $	40	$ & $	{\blue6}	$ & $	.	$ & $	.	$ & $	.	$ & $	.	$ \\
$	37	$ & $.	$ & $.	$ & $.	$ & $.	$ & $.	$ & $.	$ & $	{\red18}	$ & $	100	$ & $	109	$ & $	60	$ & $\textcolor{brown}{10}	$ & $	.	$ & $	.	$ & $	.	$ & $	.	$ \\
$	38	$ & $.	$ & $.	$ & $.	$ & $.	$ & $.	$ & $.	$ & ${\red	11}	$ & $	90	$ & $	130	$ & $	87	$ & ${\green 28}	$ & $	.	$ & $	.	$ & $	.	$ & $	.	$ \\
$	39	$ & $.	$ & $.	$ & $.	$ & $.	$ & $.	$ & $.	$ & $	{\red6}	$ & $	81	$ & $	140	$ & $	111	$ & $	47	$ & $	{\blue6}	$ & $	.	$ & $	.	$ & $	.	$  \\
$	40	$ & $.	$ & $.	$ & $.	$ & $.	$ & $.	$ & $.	$ & $	{\red3}	$ & $	67	$ & $	155	$ & $	141	$ & $	71	$ & $\textcolor{brown}{11}	$ & $	.	$ & $	.	$ & $	.	$  \\
$	41	$ & $	.	$ & $	.	$ & $	.	$ & $	.	$ & $	.	$ & $	.	$ & $	{\red1}	$ & $	{\red47}	$ & $	149	$ & $	163	$ & $	101	$ & ${\green 31}	$ & $.		$ & $	.	$ & $	.	$ \\
$	42	$ & $	.	$ & $	.	$ & $	.	$ & $	.	$ & $	.	$ & $	.	$ & $	{\red1}	$ & $	{\red31}	$ & $	153	$ & $	192	$ & $	132	$ & $	52	$ & $	{\blue7}	$ & $	.	$ & $	.	$ \\
$	43	$ & $	.	$ & $	.	$ & $	.	$ & $	.	$ & $	.	$ & $	.	$ & $	.	$ & $	{\red18}	$ & $	138	$ & $	209	$ & $	172	$ & $	81	$ & $	\textcolor{brown}{12}	$ & $	.	$ & $	.	$ \\
$	44	$ & $	.	$ & $	.	$ & $	.	$ & $	.	$ & $	.	$ & $	.	$ & $	.	$ & $	{\red11}	$ & $	121	$ & $	221	$ & $	200	$ & $	115	$ & $	{\green 34}	$ & $	.	$ & $	.	$ \\
$	45	$ & $	.	$ & $	.	$ & $	.	$ & $	.	$ & $	.	$ & $	.	$ & $	.	$ & $	{\red6}	$ & $	100	$ & $	231	$ & $	239	$ & $	158	$ & $	56	$ & $	{\blue7}	$ & $	.	$ \\
$	46	$ & $	.	$ & $	.	$ & $	.	$ & $	.	$ & $	.	$ & $	.	$ & $	.	$ & $	{\red3}	$ & $		{\red72}	$ & $	229	$ & $	271	$ & $	202	$ & $	94	$ & $\textcolor{brown}{13}	$ & $	.	$ \\
$	47	$ & $	.	$ & $	.	$ & $	.	$ & $	.	$ & $	.	$ & $	.	$ & $	.	$ & $	{\red1}	$ & $		{\red47}	$ & $	224	$ & $	296	$ & $	242	$ & $	130	$ & ${\green 37}	$ & $	.	$ \\
$	48	$ & $	.	$ & $	.	$ & $	.	$ & $	.	$ & $	.	$ & $	.	$ & $	.	$ & $	{\red1}	$ & $		{\red31}	$ & $	200	$ & $	320	$ & $	289	$ & $	180	$ & $	60	$ & $	{\blue8}	$ \\
$	49	$ & $	.	$ & $	.	$ & $	.	$ & $	.	$ & $	.	$ & $	.	$ & $	.	$ & $	.	$ & $		{\red18}	$ & $	175	$ & $	331	$ & $	342	$ & $	233	$ & $	104	$ & $\textcolor{brown}{14}	$ \\
$	50	$ & $	.	$ & $	.	$ & $	.	$ & $	.	$ & $	.	$ & $	.	$ & $	.	$ & $	.	$ & $		{\red11}	$ & $	143	$ & $	337	$ & $	374	$ & $	288	$ & $	145	$ & ${\green 40}	$ \\
\end{tabular}
}
\caption{Calculated values of $\tnd$.}
\label{t:tnd_50}
\end{center}
\end{table}

\begin{table}[ht]
\begin{center}
\includegraphics[width=\textwidth]{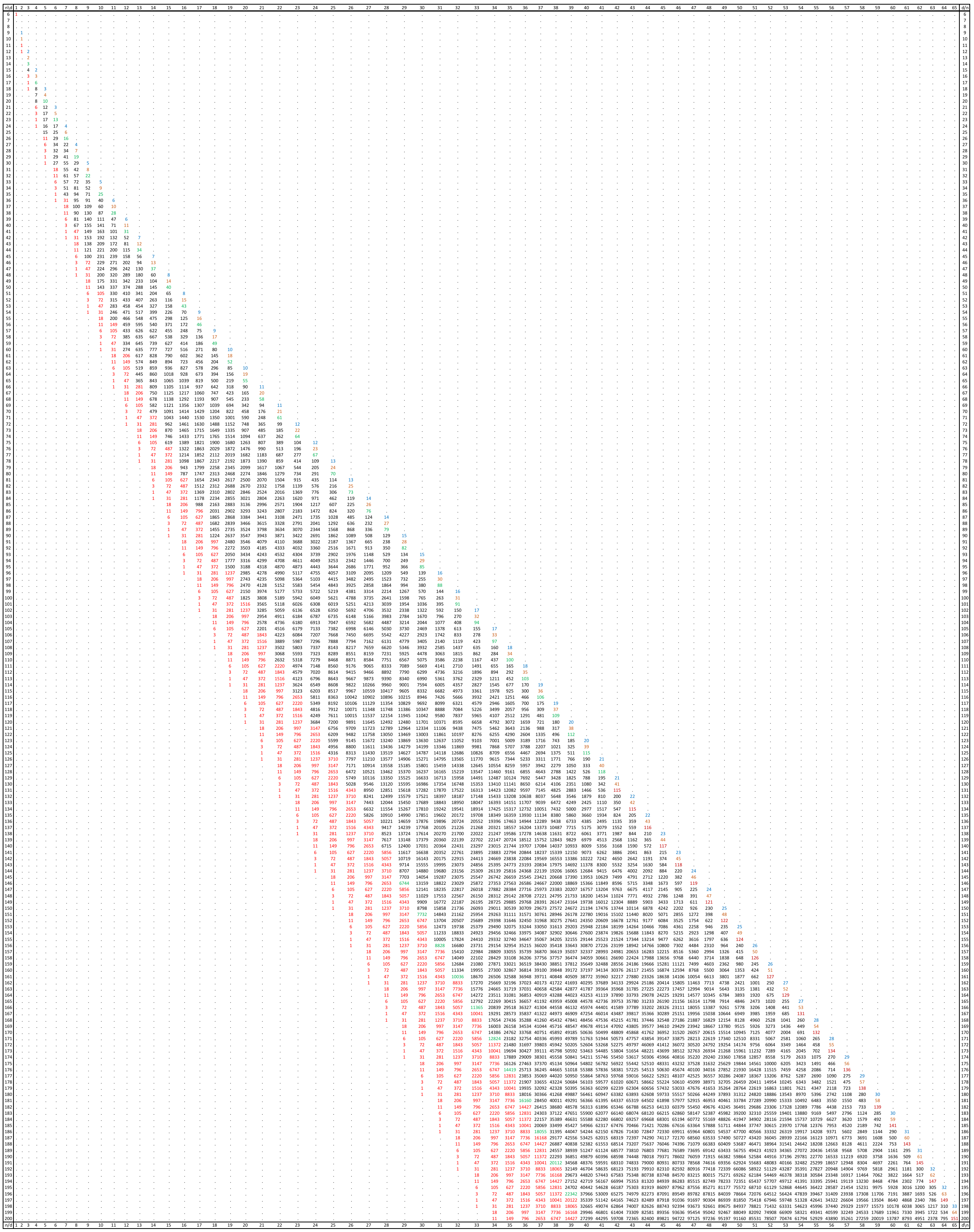}
\caption{Calculated values of $\tnd$.  The pdf can be zoomed for greater detail.}
\label{t:tnd_200}
\end{center}
\end{table}

\begin{table}[ht]
\begin{center}
{\footnotesize
\begin{tabular}{|c||c|c|c|c|c||c|c|}
\hline
$n$ & \ref{fix1} & \ref{fix2} & \ref{fix3} & \ref{fix4} & \ref{fix5} & $\tn$ & $\tn'$ \\
\hline\hline
$	1	$ & $	0	$ & $	0	$ & $	0	$ & $	0	$ & $	0	$ & $	0	$ & $	0	$ \\
$	2	$ & $	0	$ & $	0	$ & $	0	$ & $	0	$ & $	0	$ & $	0	$ & $	0	$ \\
$	3	$ & $	0	$ & $	0	$ & $	0	$ & $	0	$ & $	0	$ & $	0	$ & $	0	$ \\
$	4	$ & $	0	$ & $	0	$ & $	0	$ & $	0	$ & $	0	$ & $	0	$ & $	0	$ \\
$	5	$ & $	0	$ & $	0	$ & $	0	$ & $	0	$ & $	0	$ & $	0	$ & $	0	$ \\
$	6	$ & $	1	$ & $	1	$ & $	1	$ & $	1	$ & $	1	$ & $	1	$ & $	1	$ \\
$	7	$ & $	0	$ & $	0	$ & $	0	$ & $	0	$ & $	0	$ & $	0	$ & $	0	$ \\
$	8	$ & $	0	$ & $	0	$ & $	0	$ & $	0	$ & $	0	$ & $	0	$ & $	0	$ \\
$	9	$ & $	4	$ & $	1	$ & $	0	$ & $	2	$ & $	0	$ & $	1	$ & $	1	$ \\
$	10	$ & $	3	$ & $	0	$ & $	3	$ & $	1	$ & $	1	$ & $	1	$ & $	1	$ \\
\hline
$	11	$ & $	6	$ & $	0	$ & $	2	$ & $	2	$ & $	0	$ & $	1	$ & $	1	$ \\
$	12	$ & $	11	$ & $	2	$ & $	3	$ & $	5	$ & $	1	$ & $	3	$ & $	3	$ \\
$	13	$ & $	18	$ & $	0	$ & $	2	$ & $	4	$ & $	0	$ & $	2	$ & $	2	$ \\
$	14	$ & $	27	$ & $	0	$ & $	3	$ & $	5	$ & $	1	$ & $	3	$ & $	3	$ \\
$	15	$ & $	42	$ & $	3	$ & $	6	$ & $	10	$ & $	0	$ & $	6	$ & $	7	$ \\
$	16	$ & $	57	$ & $	0	$ & $	9	$ & $	9	$ & $	1	$ & $	6	$ & $	7	$ \\
$	17	$ & $	84	$ & $	0	$ & $	8	$ & $	10	$ & $	0	$ & $	7	$ & $	9	$ \\
$	18	$ & $	120	$ & $	3	$ & $	8	$ & $	18	$ & $	2	$ & $	12	$ & $	14	$ \\
$	19	$ & $	138	$ & $	0	$ & $	6	$ & $	18	$ & $	0	$ & $	11	$ & $	13	$ \\
$	20	$ & $	234	$ & $	0	$ & $	18	$ & $	22	$ & $	2	$ & $	18	$ & $	24	$ \\
\hline
$	21	$ & $	268	$ & $	4	$ & $	12	$ & $	28	$ & $	0	$ & $	21	$ & $	28	$ \\
$	22	$ & $	354	$ & $	0	$ & $	18	$ & $	30	$ & $	2	$ & $	25	$ & $	34	$ \\
$	23	$ & $	480	$ & $	0	$ & $	16	$ & $	36	$ & $	0	$ & $	31	$ & $	44	$ \\
$	24	$ & $	560	$ & $	5	$ & $	16	$ & $	42	$ & $	2	$ & $	38	$ & $	54	$ \\
$	25	$ & $	750	$ & $	0	$ & $	26	$ & $	46	$ & $	0	$ & $	46	$ & $	69	$ \\
$	26	$ & $	897	$ & $	0	$ & $	33	$ & $	55	$ & $	3	$ & $	56	$ & $	83	$ \\
$	27	$ & $	1082	$ & $	5	$ & $	30	$ & $	62	$ & $	0	$ & $	66	$ & $	101	$ \\
$	28	$ & $	1326	$ & $	0	$ & $	30	$ & $	66	$ & $	2	$ & $	76	$ & $	118	$ \\
$	29	$ & $	1584	$ & $	0	$ & $	36	$ & $	78	$ & $	0	$ & $	90	$ & $	141	$ \\
$	30	$ & $	2013	$ & $	6	$ & $	57	$ & $	93	$ & $	3	$ & $	117	$ & $	186	$ \\
\hline
$	31	$ & $	2256	$ & $	0	$ & $	48	$ & $	92	$ & $	0	$ & $	123	$ & $	200	$ \\
$	32	$ & $	2811	$ & $	0	$ & $	51	$ & $	107	$ & $	3	$ & $	151	$ & $	247	$ \\
$	33	$ & $	3258	$ & $	6	$ & $	50	$ & $	124	$ & $	0	$ & $	175	$ & $	288	$ \\
$	34	$ & $	3738	$ & $	0	$ & $	66	$ & $	124	$ & $	4	$ & $	196	$ & $	328	$ \\
$	35	$ & $	4554	$ & $	0	$ & $	70	$ & $	142	$ & $	0	$ & $	234	$ & $	397	$ \\
$	36	$ & $	5047	$ & $	7	$ & $	83	$ & $	161	$ & $	3	$ & $	264	$ & $	446	$ \\
$	37	$ & $	5886	$ & $	0	$ & $	78	$ & $	168	$ & $	0	$ & $	297	$ & $	510	$ \\
$	38	$ & $	6942	$ & $	0	$ & $	78	$ & $	184	$ & $	4	$ & $	346	$ & $	598	$ \\
$	39	$ & $	7778	$ & $	8	$ & $	98	$ & $	208	$ & $	0	$ & $	391	$ & $	678	$ \\
$	40	$ & $	9048	$ & $	0	$ & $	120	$ & $	220	$ & $	4	$ & $	448	$ & $	784	$ \\
\hline
$	41	$ & $	10104	$ & $	0	$ & $	108	$ & $	230	$ & $	0	$ & $	492	$ & $	869	$ \\
$	42	$ & $	11648	$ & $	8	$ & $	112	$ & $	260	$ & $	4	$ & $	568	$ & $	1004	$ \\
$	43	$ & $	13086	$ & $	0	$ & $	122	$ & $	278	$ & $	0	$ & $	630	$ & $	1121	$ \\
$	44	$ & $	14652	$ & $	0	$ & $	144	$ & $	290	$ & $	4	$ & $	702	$ & $	1257	$ \\
$	45	$ & $	16698	$ & $	9	$ & $	146	$ & $	320	$ & $	0	$ & $	797	$ & $	1434	$ \\
$	46	$ & $	18609	$ & $	0	$ & $	165	$ & $	347	$ & $	5	$ & $	884	$ & $	1592	$ \\
$	47	$ & $	20778	$ & $	0	$ & $	162	$ & $	364	$ & $	0	$ & $	977	$ & $	1772	$ \\
$	48	$ & $	23170	$ & $	10	$ & $	174	$ & $	390	$ & $	4	$ & $	1089	$ & $	1981	$ \\
$	49	$ & $	26118	$ & $	0	$ & $	190	$ & $	420	$ & $	0	$ & $	1217	$ & $	2224	$ \\
$	50	$ & $	28755	$ & $	0	$ & $	219	$ & $	445	$ & $	5	$ & $	1338	$ & $	2451	$ \\
\hline
\end{tabular}%
~\qquad~%
\begin{tabular}{|c||c|c|c|c|c||c|c|}
\hline
$n$ & \ref{fix1} & \ref{fix2} & \ref{fix3} & \ref{fix4} & \ref{fix5} & $\tn$ & $\tn'$ \\
\hline\hline
$	51	$ & $	31696	$ & $	10	$ & $	216	$ & $	472	$ & $	0	$ & $	1469	$ & $	2702	$ \\
$	52	$ & $	35265	$ & $	0	$ & $	225	$ & $	501	$ & $	5	$ & $	1624	$ & $	2995	$ \\
$	53	$ & $	38616	$ & $	0	$ & $	232	$ & $	532	$ & $	0	$ & $	1771	$ & $	3276	$ \\
$	54	$ & $	42950	$ & $	11	$ & $	266	$ & $	568	$ & $	6	$ & $	1970	$ & $	3653	$ \\
$	55	$ & $	47136	$ & $	0	$ & $	272	$ & $	592	$ & $	0	$ & $	2146	$ & $	3996	$ \\
$	56	$ & $	51519	$ & $	0	$ & $	291	$ & $	635	$ & $	5	$ & $	2343	$ & $	4366	$ \\
$	57	$ & $	56838	$ & $	12	$ & $	302	$ & $	676	$ & $	0	$ & $	2579	$ & $	4820	$ \\
$	58	$ & $	61590	$ & $	0	$ & $	330	$ & $	692	$ & $	6	$ & $	2782	$ & $	5215	$ \\
$	59	$ & $	67524	$ & $	0	$ & $	340	$ & $	744	$ & $	0	$ & $	3042	$ & $	5712	$ \\
$	60	$ & $	73716	$ & $	12	$ & $	364	$ & $	798	$ & $	6	$ & $	3322	$ & $	6242	$ \\
\hline
$	61	$ & $	80058	$ & $	0	$ & $	366	$ & $	818	$ & $	0	$ & $	3586	$ & $	6763	$ \\
$	62	$ & $	87486	$ & $	0	$ & $	390	$ & $	866	$ & $	6	$ & $	3912	$ & $	7388	$ \\
$	63	$ & $	94438	$ & $	13	$ & $	406	$ & $	924	$ & $	0	$ & $	4221	$ & $	7980	$ \\
$	64	$ & $	102546	$ & $	0	$ & $	438	$ & $	956	$ & $	6	$ & $	4568	$ & $	8655	$ \\
$	65	$ & $	111510	$ & $	0	$ & $	454	$ & $	1000	$ & $	0	$ & $	4953	$ & $	9406	$ \\
$	66	$ & $	120217	$ & $	13	$ & $	465	$ & $	1063	$ & $	7	$ & $	5339	$ & $	10143	$ \\
$	67	$ & $	129474	$ & $	0	$ & $	478	$ & $	1106	$ & $	0	$ & $	5731	$ & $	10909	$ \\
$	68	$ & $	140337	$ & $	0	$ & $	537	$ & $	1151	$ & $	7	$ & $	6204	$ & $	11829	$ \\
$	69	$ & $	150746	$ & $	14	$ & $	542	$ & $	1214	$ & $	0	$ & $	6657	$ & $	12707	$ \\
$	70	$ & $	162687	$ & $	0	$ & $	579	$ & $	1265	$ & $	7	$ & $	7169	$ & $	13702	$ \\
\hline
$	71	$ & $	174732	$ & $	0	$ & $	580	$ & $	1320	$ & $	0	$ & $	7683	$ & $	14706	$ \\
$	72	$ & $	187293	$ & $	15	$ & $	593	$ & $	1381	$ & $	7	$ & $	8230	$ & $	15766	$ \\
$	73	$ & $	201954	$ & $	0	$ & $	646	$ & $	1446	$ & $	0	$ & $	8857	$ & $	16991	$ \\
$	74	$ & $	215556	$ & $	0	$ & $	696	$ & $	1502	$ & $	8	$ & $	9446	$ & $	18137	$ \\
$	75	$ & $	230754	$ & $	15	$ & $	690	$ & $	1556	$ & $	0	$ & $	10095	$ & $	19412	$ \\
$	76	$ & $	248271	$ & $	0	$ & $	723	$ & $	1637	$ & $	7	$ & $	10846	$ & $	20870	$ \\
$	77	$ & $	263898	$ & $	0	$ & $	734	$ & $	1702	$ & $	0	$ & $	11513	$ & $	22175	$ \\
$	78	$ & $	283060	$ & $	16	$ & $	808	$ & $	1770	$ & $	8	$ & $	12345	$ & $	23801	$ \\
$	79	$ & $	301542	$ & $	0	$ & $	806	$ & $	1840	$ & $	0	$ & $	13125	$ & $	25330	$ \\
$	80	$ & $	321228	$ & $	0	$ & $	840	$ & $	1910	$ & $	8	$ & $	13969	$ & $	26979	$ \\
\hline
$	81	$ & $	342986	$ & $	17	$ & $	870	$ & $	1990	$ & $	0	$ & $	14903	$ & $	28811	$ \\
$	82	$ & $	363870	$ & $	0	$ & $	906	$ & $	2054	$ & $	8	$ & $	15790	$ & $	30549	$ \\
$	83	$ & $	387786	$ & $	0	$ & $	946	$ & $	2140	$ & $	0	$ & $	16811	$ & $	32552	$ \\
$	84	$ & $	411902	$ & $	17	$ & $	990	$ & $	2240	$ & $	8	$ & $	17854	$ & $	34584	$ \\
$	85	$ & $	437766	$ & $	0	$ & $	1010	$ & $	2294	$ & $	0	$ & $	18940	$ & $	36733	$ \\
$	86	$ & $	465111	$ & $	0	$ & $	1035	$ & $	2383	$ & $	9	$ & $	20107	$ & $	39018	$ \\
$	87	$ & $	491910	$ & $	18	$ & $	1086	$ & $	2492	$ & $	0	$ & $	21261	$ & $	41276	$ \\
$	88	$ & $	521439	$ & $	0	$ & $	1143	$ & $	2557	$ & $	9	$ & $	22511	$ & $	43739	$ \\
$	89	$ & $	552630	$ & $	0	$ & $	1146	$ & $	2646	$ & $	0	$ & $	23831	$ & $	46339	$ \\
$	90	$ & $	585595	$ & $	19	$ & $	1211	$ & $	2765	$ & $	9	$ & $	25251	$ & $	49115	$ \\
\hline
$	91	$ & $	618420	$ & $	0	$ & $	1232	$ & $	2838	$ & $	0	$ & $	26631	$ & $	51843	$ \\
$	92	$ & $	654549	$ & $	0	$ & $	1317	$ & $	2933	$ & $	9	$ & $	28173	$ & $	54875	$ \\
$	93	$ & $	691444	$ & $	19	$ & $	1328	$ & $	3046	$ & $	0	$ & $	29744	$ & $	57965	$ \\
$	94	$ & $	729144	$ & $	0	$ & $	1380	$ & $	3140	$ & $	10	$ & $	31341	$ & $	61107	$ \\
$	95	$ & $	771042	$ & $	0	$ & $	1422	$ & $	3246	$ & $	0	$ & $	33116	$ & $	64609	$ \\
$	96	$ & $	811697	$ & $	20	$ & $	1453	$ & $	3351	$ & $	9	$ & $	34849	$ & $	68018	$ \\
$	97	$ & $	855402	$ & $	0	$ & $	1522	$ & $	3456	$ & $	0	$ & $	36696	$ & $	71664	$ \\
$	98	$ & $	902544	$ & $	0	$ & $	1560	$ & $	3566	$ & $	10	$ & $	38695	$ & $	75602	$ \\
$	99	$ & $	947912	$ & $	20	$ & $	1592	$ & $	3668	$ & $	0	$ & $	40619	$ & $	79404	$ \\
$	100	$ & $	999738	$ & $	0	$ & $	1674	$ & $	3798	$ & $	10	$ & $	42817	$ & $	83730	$ \\
\hline
\end{tabular}
}
\caption{Calculated values of $\fix(\si)$, for $\si$ of types \ref{fix1}--\ref{fix5}, as well as $\tn$ and $\tn'$.}
\label{t:fix}
\end{center}
\end{table}

\begin{figure}[ht]
\begin{center}
\includegraphics[width=.8\textwidth]{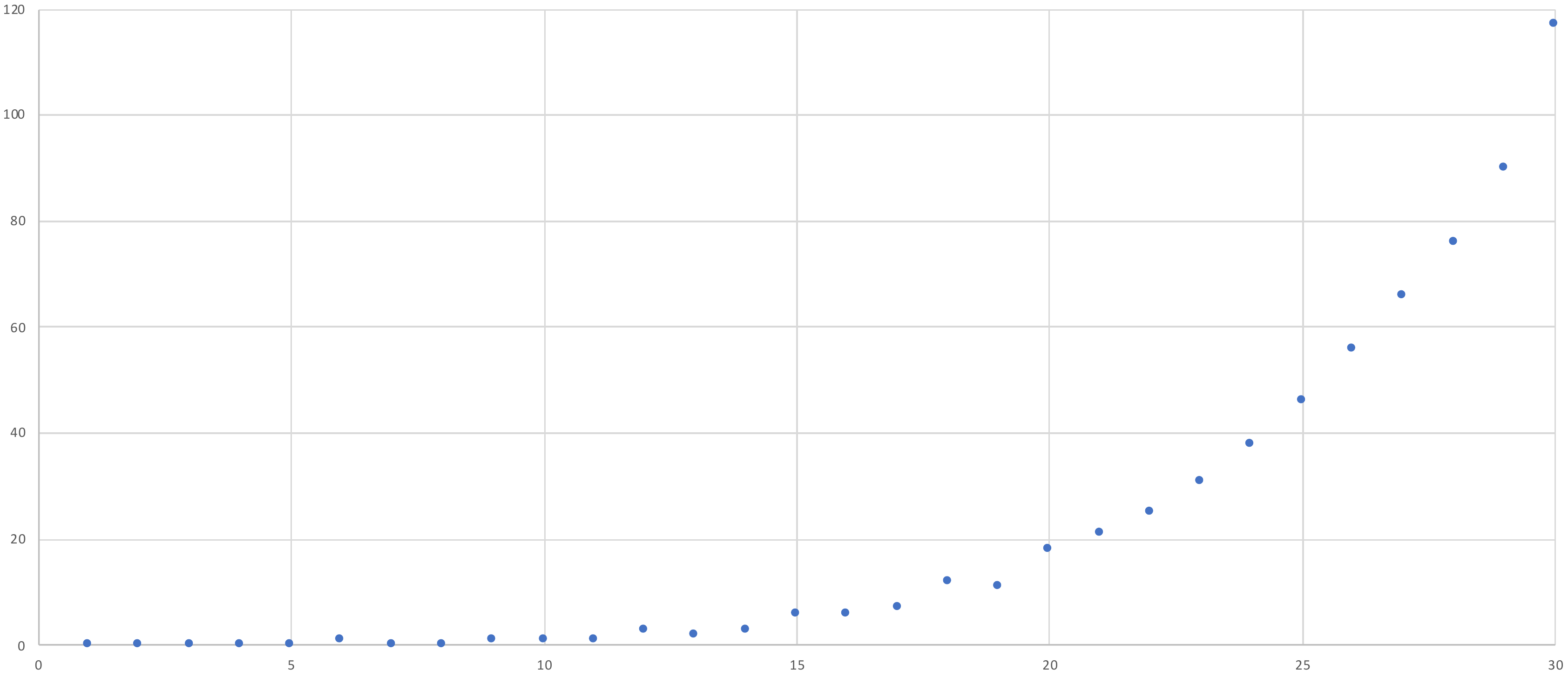}
\caption{Calculated values of $\tn$, $1\leq n\leq30$.}
\label{f:tn_30}
\end{center}
\end{figure}

\begin{figure}[ht]
\begin{center}
\includegraphics[width=.8\textwidth]{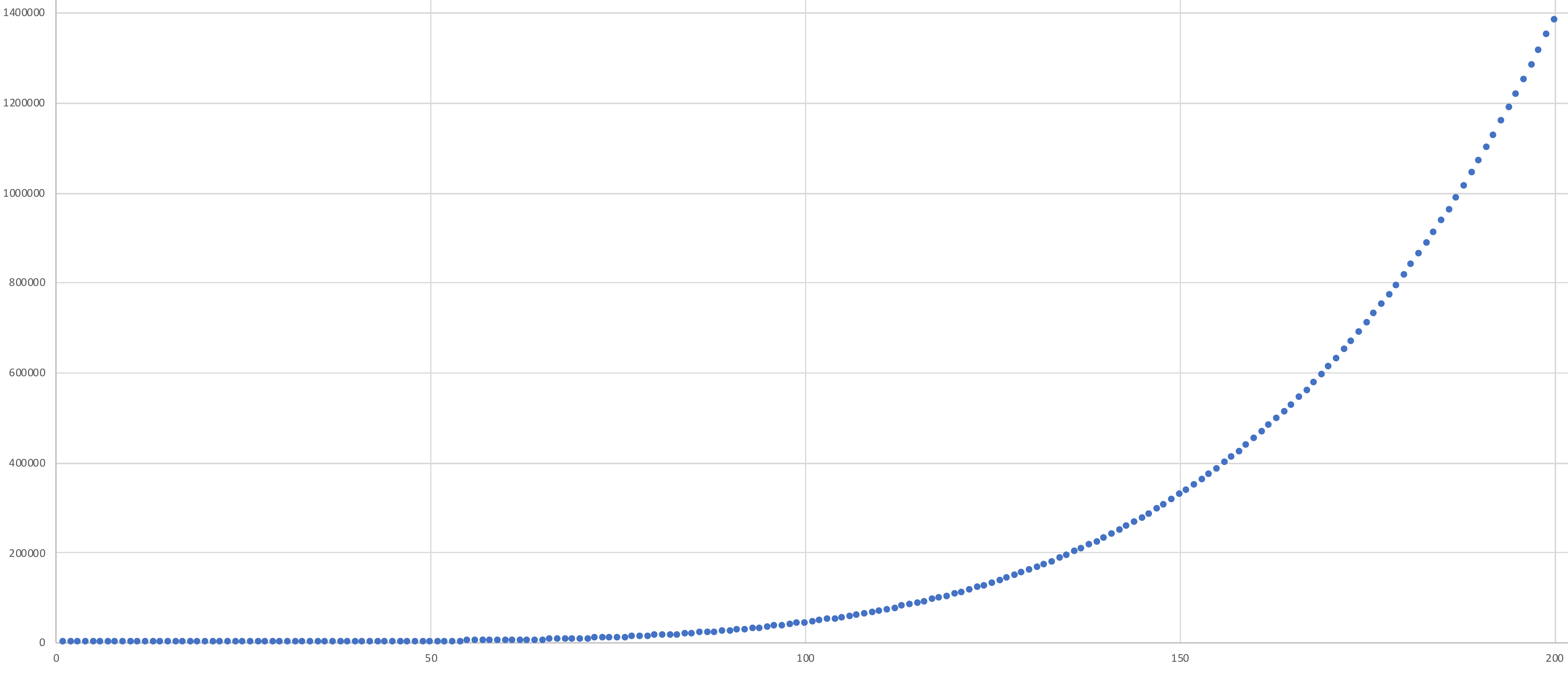}
\caption{Calculated values of $\tn$, $1\leq n\leq200$.}
\label{f:tn}
\end{center}
\end{figure}

\begin{figure}[ht]
\begin{center}
\includegraphics[width=.8\textwidth]{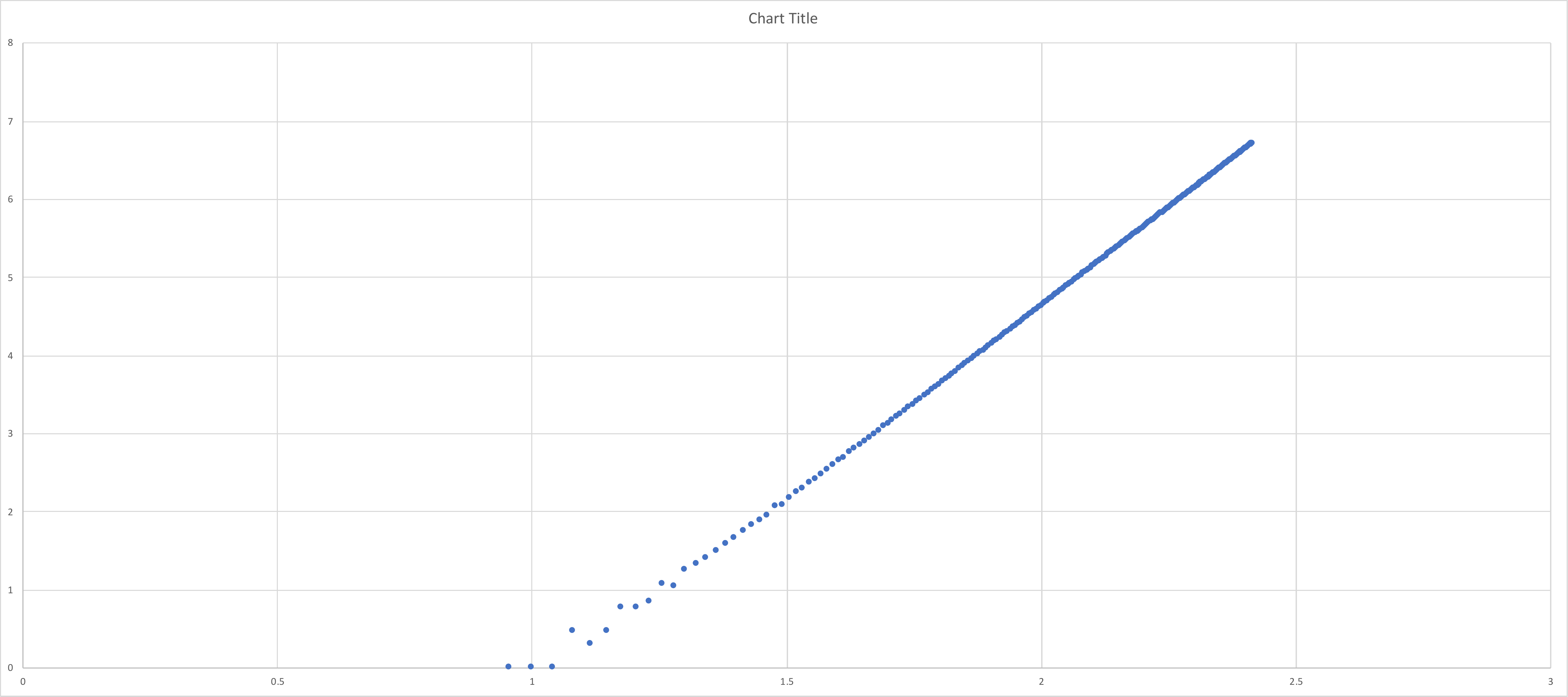}
\caption{Log-log plot of calculated values of $\tn$, $1\leq n\leq200$.}
\label{f:tn_log}
\end{center}
\end{figure}

\begin{figure}[ht]
\begin{center}
\includegraphics[width=\textwidth]{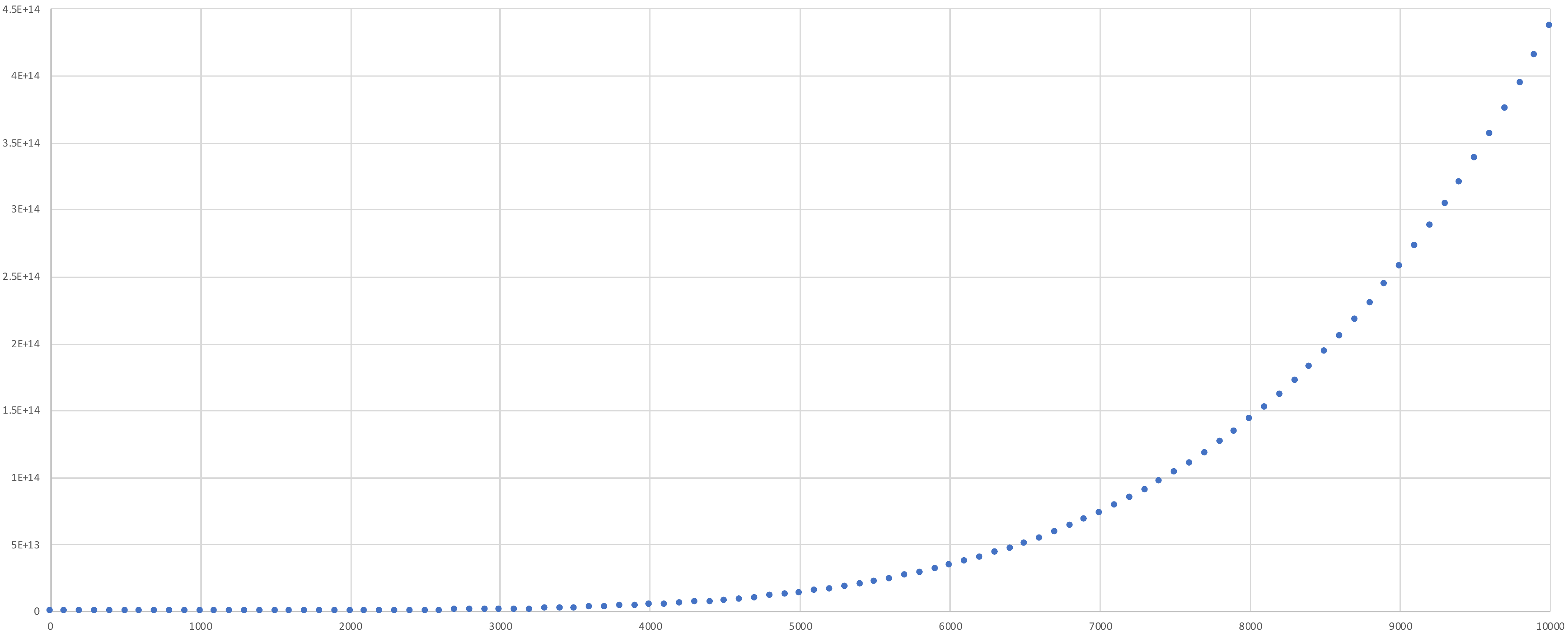}
\caption{Calculated values of $\tn$, $1\leq n\leq10000$ (by hundreds).}
\label{f:tn_100s}
\end{center}
\end{figure}

\begin{figure}[ht]
\begin{center}
\includegraphics[width=\textwidth]{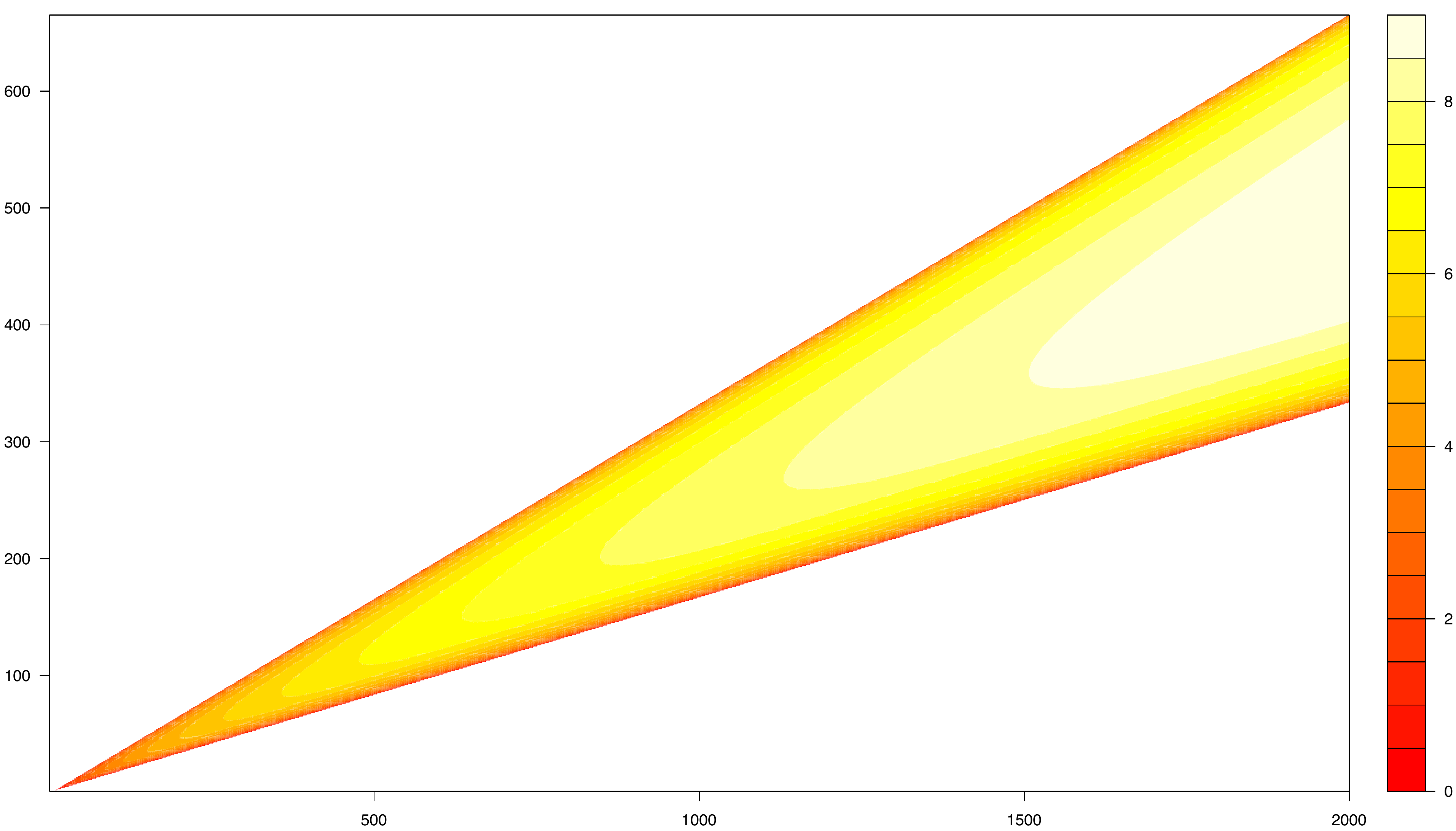}
\caption{Heat map of the numbers $\log_{10}(\tnd)$, with $1\leq n\leq2000$ and $1\leq d\leq666$ on the horizontal and vertical axes, respectively.}
\label{f:HM}
\end{center}
\end{figure}

\begin{figure}[ht]
\begin{center}
\includegraphics[height=5cm]{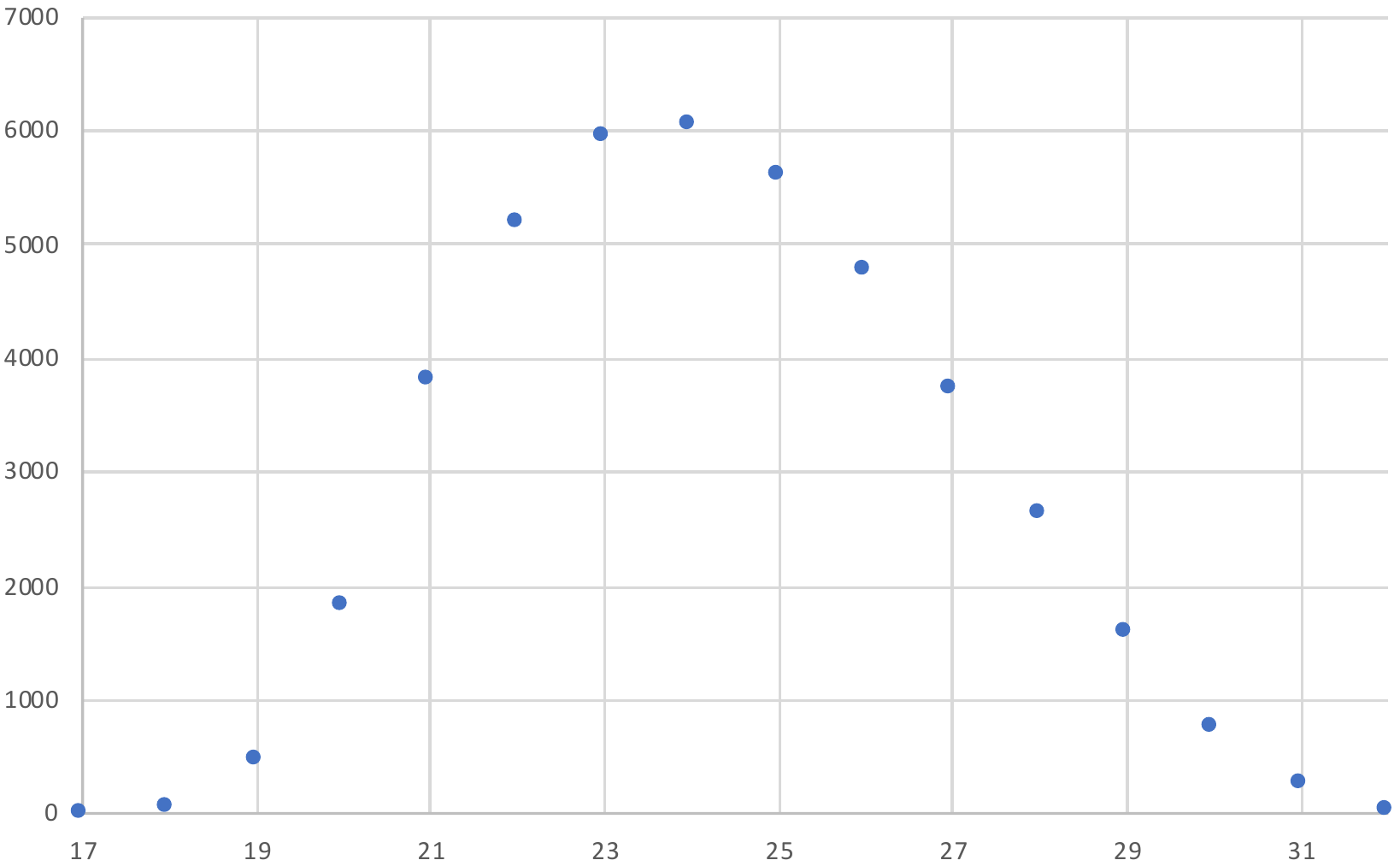}
\quad
\includegraphics[height=5cm]{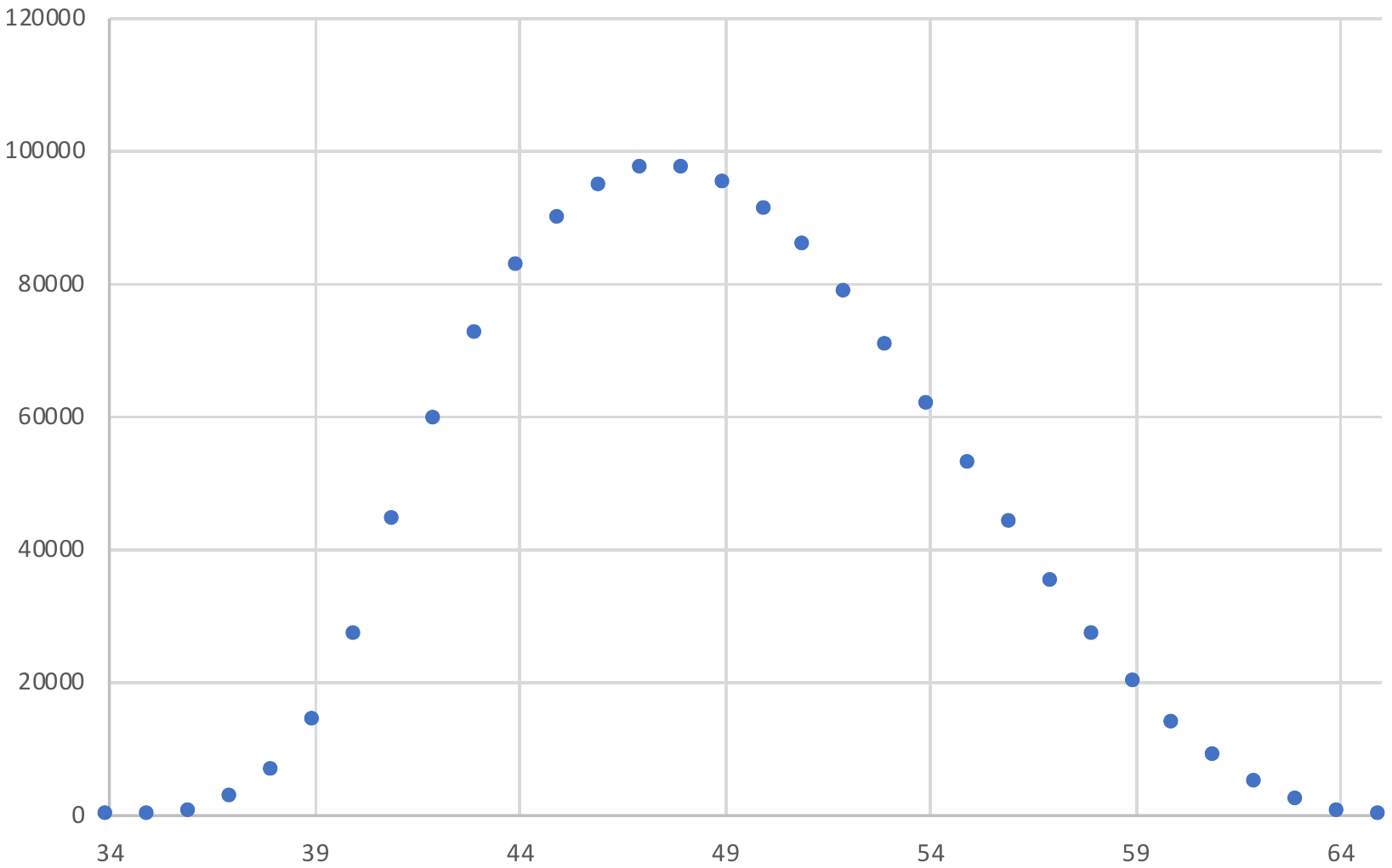}
\caption{Calculated values of $\td_{100}$ (left) and $\td_{200}$ (right).}
\label{f:tnd_100_200}
\end{center}
\end{figure}

\begin{figure}[ht]
\begin{center}
\includegraphics[width=.8\textwidth]{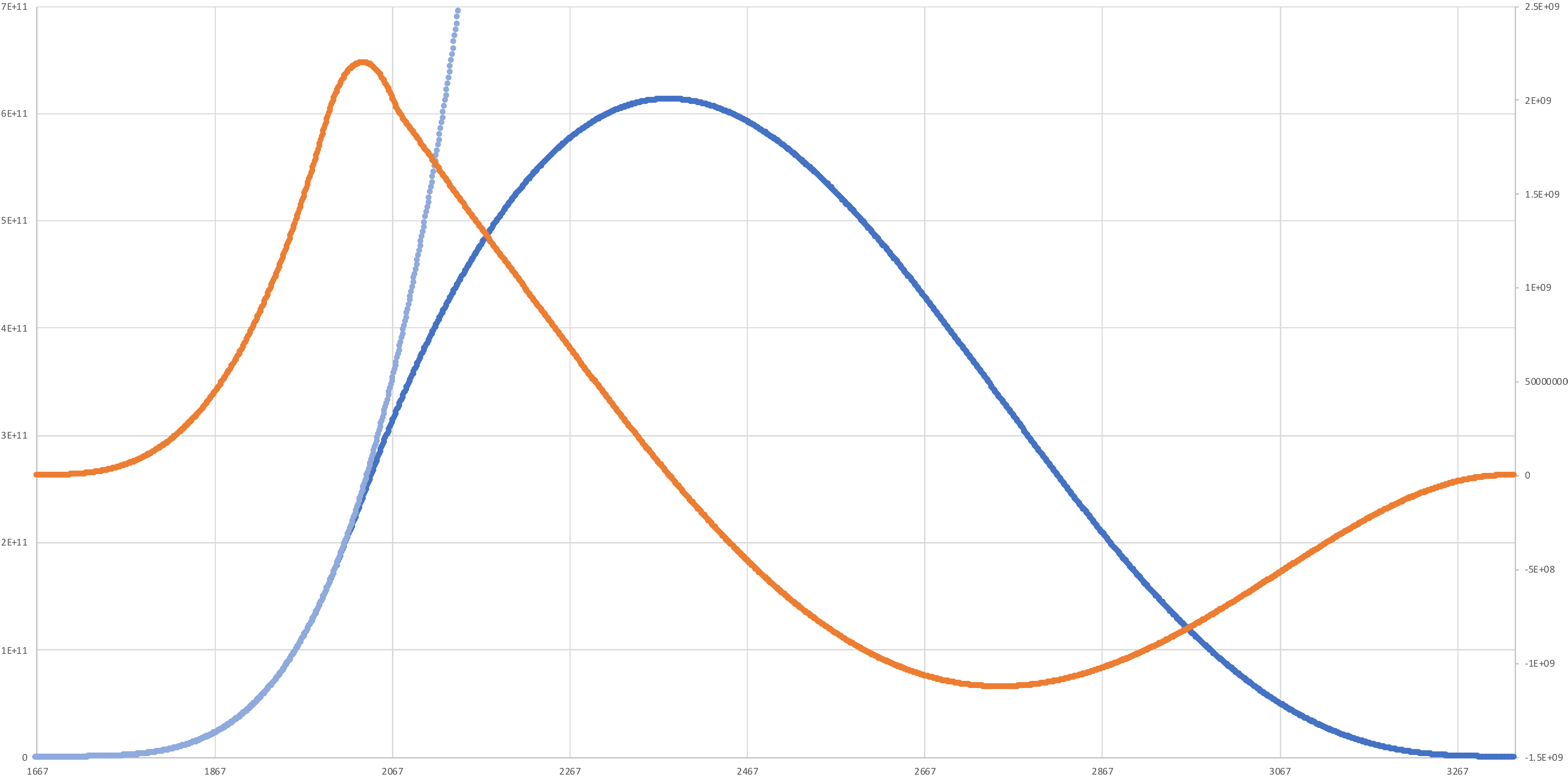}
\caption{Calculated values of $\td_{10000}$ (dark blue) and the discrete derivative (orange).  The light blue curve is explained in Section \ref{ss:ak}.}
\label{f:tnd_10000}
\end{center}
\end{figure}

\begin{figure}[ht]
\begin{center}
\includegraphics[height=4.5cm]{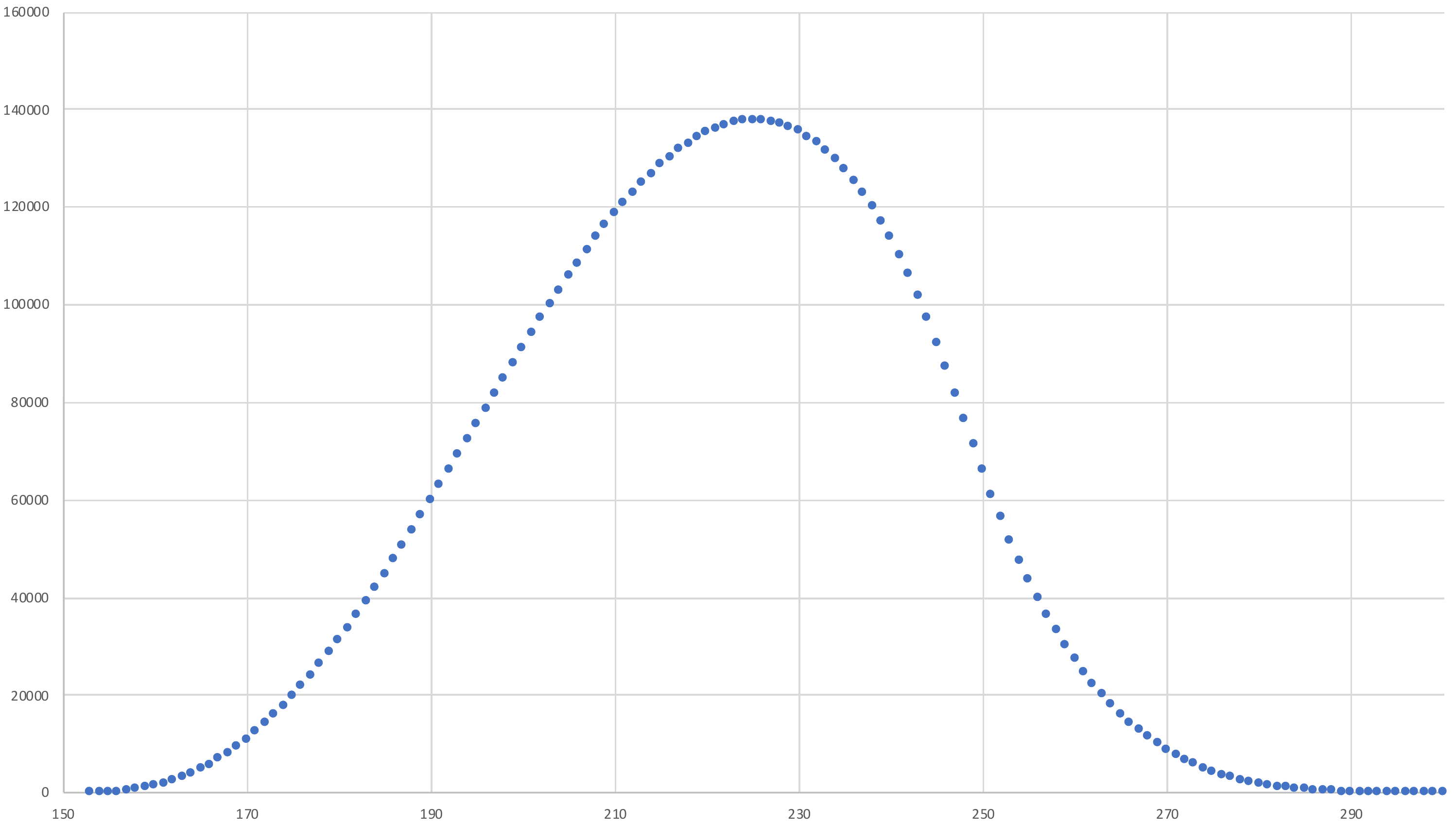}
\quad
\includegraphics[height=4.5cm]{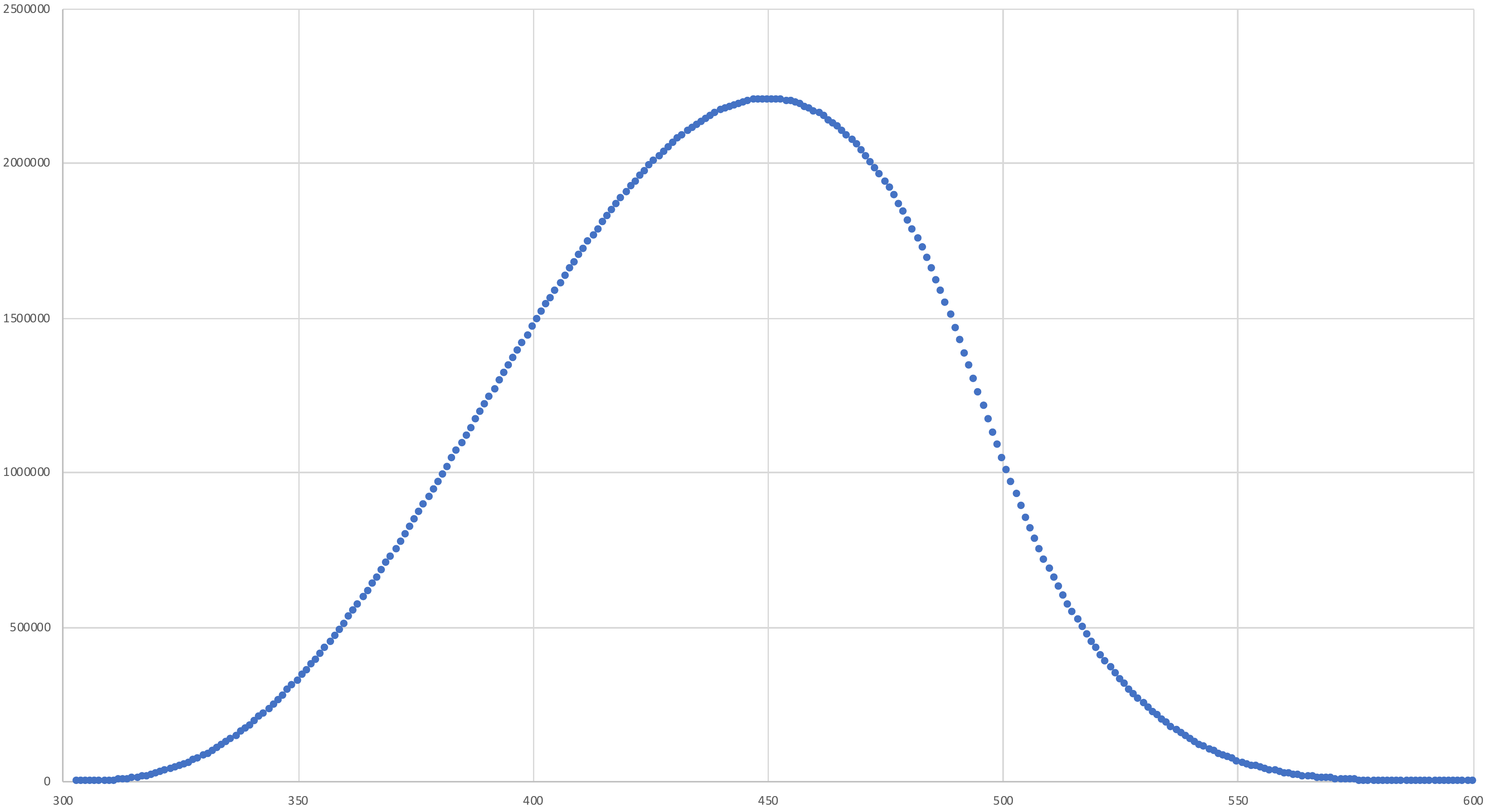}
\caption{Calculated values of ${}^{50}\tn$ (left) and ${}^{100}\tn$ (right).}
\label{f:tnd_50_100}
\end{center}
\end{figure}

\begin{figure}[ht]
\begin{center}
\includegraphics[width=.8\textwidth]{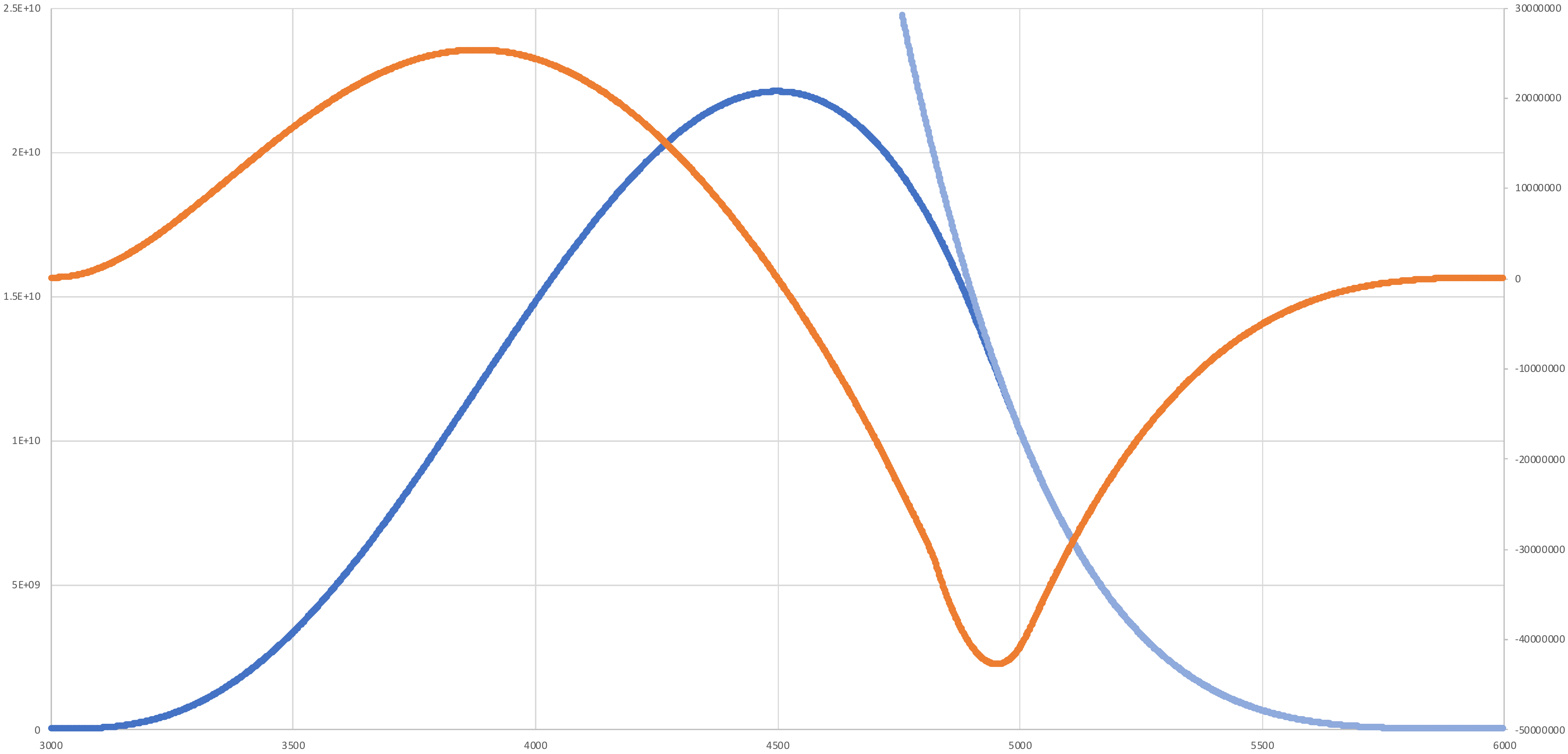}
\caption{Calculated values of ${}^{1000}\tn$ (dark blue) and the discrete derivative (orange).  The light blue curve is explained in Section \ref{ss:ak}.}
\label{f:tnd_1000}
\end{center}
\end{figure}


%
%
%
%
%
%
%
%
%
%

\section{Observations from the data}\label{s:obs}

The computational data presented in Section \ref{s:computed} has allowed us to discover some interesting patterns, and has led to a number of natural conjectures, which we discuss in the current section.  We do not mean to imply that the following discussion is exhaustive.  We hope that other researchers may be able to shed light on some of our conjectures, and make further discoveries.

\subsection{Asymptotics}\label{ss:as}

First and foremost, we had hoped that the computed values of $\tn$ might have suggested an obvious formula, perhaps after consultation with the OEIS \cite{OEIS}.  However, this regrettably was not the case.  The ``bumpy'' nature of the $\tn$ sequence, which can be seen in Figure \ref{f:tn_30}, suggests there is no simple formula.  However, the sequence does appear to become smoother at a larger scale, as can be seen in Figure~\ref{f:tn}.  This is (apparently) confirmed further in Figure \ref{f:tn_log}, which gives a log-log plot, and suggests a possible power law.

By analogy with the case of triangles, where the answer is asymptotic to $\frac{n^2}{48}$ (cf.~Honsberger's Theorem), we might wonder if $\tn\sim \frac{n^k}C$ for some integers $k$ and $C$.  
If this were so, then 
\[
\frac{t_{n+1}}{\tn} \sim \frac{(n+1)^k}{n^k} = \frac{n^k+kn^{k-1}+\cdots}{n^k} \sim 1+\frac kn .
\]
Computed values of $\left(\frac{t_{n+1}}{\tn}-1\right)\times n$ strongly suggest a likely value of $k=5$.  Computed values of~$C_n=\frac{n^5}{\tn}$ suggest a possible value in the order of $C\approx 229000$.  Table \ref{t:C} gives calculated values of $C_n$ for $n$ up to $10000$ (by hundreds); Table \ref{t:rhon} goes up to $n=20000$ (by thousands).  The third column of Table \ref{t:C} gives the difference $d_n = C_n-C_{n-100}$, which shows that the fall of $C_n$ becomes slower as $n$ increases.  The values $d_n$ are themselves decreasing, as governed by the ratio $r_n=d_n/d_{n-100}$.  This ratio seems to be tending towards $0.98$ or $0.99$, or thereabouts.  If so, then continuing the trend, $C_n$ seems to approach $229028$ or $229024$, respectively.  So we have the following tentative conjecture; we have made an explicit choice for the denominator so as to have a concrete statement, but there is some level of uncertainty, as just discussed.

\begin{con}\label{c:as}
The number $\tn$ of integer tetrahedra of perimeter $n$, up to congruence, satisfies
\[
\tn \sim \frac{n^5}{229024} \qquad\text{as $n\to\infty$.}
\]
\end{con}

Taking Conjecture \ref{c:as} as given, one could then aim to analyse the sequence $s_n = \tn-\frac{n^5}{229024}$.  For example, one might hope that this is a quartic polynomial.  As above, evidence for this could be gathered by calculating values of $\left(\frac{s_{n+1}}{s_n}-1\right)\times n$.  However, these values appear to be quite chaotic, even when we look at different values of $C$ around the conjectured value of 229024.  However, we have more success by splitting up the $s_n$ sequence, and looking at subsequences $s_i, s_{i+k},s_{i+2k},\dots$, for various integers $k\geq2$ and for each $0\leq i\leq k-1$, especially when $k$ is a multiple of $12$.  Here the relevant value to calculate is $\left(\frac{s_{n+k}}{s_n}-1\right)\times \frac nk$, and this does seem to approach $4$ as $n\to\infty$ (for $k\equiv0\Mod{12}$), for various $C$.
Might there be a multi-case formula for $\tn$?  Attempts to fit polynomials to subsequences (of the kind described above), using various statistical packages, have not led to any success.
In any case, the uncertainty as to the value of $C\approx 229024$ leads to even greater uncertainty~here.

Here is an additional note about Conjecture \ref{c:as} (we will return to it again in Section \ref{ss:max}).  Recall that $\tn=\frac1{24}\sum_{\si\in\S_4}\fix(\si)$; cf.~\eqref{e:BS1}.  For $\si\in\S_4\sm\{\id_4\}$, it is easy to see that $\fix(\si)$ is bounded above by~$n^4$; this can be seen by examining the constraints on the labels of a graph $G$ fixed by $\si$ in cases \ref{fix2}--\ref{fix5}, as listed in Section \ref{s:BL}.  Thus, if~$\tn$ is indeed (asymptotically) a quintic, then so too is $\fix(\id_4)=|\T_n|$, and so 
\[
\tn \sim \frac{|\T_n|}{24} \qquad\text{as $n\to\infty$.}
\]
This can be interpreted as saying that virtually all tetrahedra have no symmetry at all.  It also leads to an equivalent formulation of Conjecture \ref{c:as}, again subject to some uncertainty in the value of the stated numerical coefficient:

\begin{con}\label{c:as2}
The set $\T_n$ of graphs corresponding to integer tetrahedra of perimeter $n$ satisfies
\[
|\T_n| \sim \frac{3n^5}{28628} \qquad\text{as $n\to\infty$.}
\]
\end{con}

\begin{rem}\label{r:Kurz_as}
Asymptotics of the numbers $\td$ were not discussed explicitly by Kurz \cite{Kurz2007}.  However, it was observed on \cite[p.~7]{Kurz2007} that $\td$ seems to be approximately $0.103$ times the number of graphs from $\G$ satisfying \ref{T1'} and \ref{T2}--\ref{T5}; see Proposition \ref{p:Tn} and Remark \ref{r:td_tnd}.  The language used in \cite{Kurz2007} was different, and instead spoke of symmetric matrices, but it follows from \cite[Lemma 4]{Kurz2007} that the number of such graphs/matrices is asymptotic to~$\frac{d^5}4$.  Thus, a tentative conjecture for an asymptotic expression for $\td$ is
\begin{equation}\label{e:td_con}
\td \sim 0.103\times\frac{d^5}4 \approx \frac{d^5}{38.835} \qquad\text{as $d\to\infty$.}
\end{equation}
This fares quite well with computed values; see Table~\ref{t:rhod}.  The non-integrality of the denominator $38.835$ in \eqref{e:td_con} suggests that the denominator in Conjecture \ref{c:as} for the $t_n$ sequence might not be an integer either.
Again, we will say more about this in Section~\ref{ss:max}.
\end{rem}

\begin{table}[ht]
\begin{center}
{\footnotesize
\begin{tabular}{|c||c|c|c|}
\hline
$n$ & $C_n$ & $d_n$ & $r_n$  \\
\hline\hline
$	100	$ & $	233552.094	$ & $		$ & $		$ \\
$	200	$ & $	231442.975	$ & $	-2109.119	$ & $		$ \\
$	300	$ & $	230489.251	$ & $	-953.724	$ & $	0.452	$ \\
$	400	$ & $	230045.517	$ & $	-443.734	$ & $	0.465	$ \\
$	500	$ & $	229786.224	$ & $	-259.293	$ & $	0.584	$ \\
$	600	$ & $	229619.532	$ & $	-166.692	$ & $	0.643	$ \\
$	700	$ & $	229509.307	$ & $	-110.225	$ & $	0.661	$ \\
$	800	$ & $	229428.544	$ & $	-80.763	$ & $	0.733	$ \\
$	900	$ & $	229367.780	$ & $	-60.765	$ & $	0.752	$ \\
$	1000	$ & $	229321.619	$ & $	-46.161	$ & $	0.760	$ \\
\hline
$	1100	$ & $	229284.562	$ & $	-37.057	$ & $	0.803	$ \\
$	1200	$ & $	229255.271	$ & $	-29.291	$ & $	0.790	$ \\
$	1300	$ & $	229230.648	$ & $	-24.623	$ & $	0.841	$ \\
$	1400	$ & $	229210.453	$ & $	-20.195	$ & $	0.820	$ \\
$	1500	$ & $	229193.275	$ & $	-17.179	$ & $	0.851	$ \\
$	1600	$ & $	229178.649	$ & $	-14.625	$ & $	0.851	$ \\
$	1700	$ & $	229166.106	$ & $	-12.543	$ & $	0.858	$ \\
$	1800	$ & $	229155.133	$ & $	-10.973	$ & $	0.875	$ \\
$	1900	$ & $	229145.626	$ & $	-9.507	$ & $	0.866	$ \\
$	2000	$ & $	229137.158	$ & $	-8.468	$ & $	0.891	$ \\
\hline
$	2100	$ & $	229129.647	$ & $	-7.511	$ & $	0.887	$ \\
$	2200	$ & $	229122.995	$ & $	-6.652	$ & $	0.886	$ \\
$	2300	$ & $	229116.976	$ & $	-6.019	$ & $	0.905	$ \\
$	2400	$ & $	229111.617	$ & $	-5.359	$ & $	0.890	$ \\
$	2500	$ & $	229106.724	$ & $	-4.893	$ & $	0.913	$ \\
$	2600	$ & $	229102.310	$ & $	-4.414	$ & $	0.902	$ \\
$	2700	$ & $	229098.275	$ & $	-4.035	$ & $	0.914	$ \\
$	2800	$ & $	229094.575	$ & $	-3.699	$ & $	0.917	$ \\
$	2900	$ & $	229091.195	$ & $	-3.381	$ & $	0.914	$ \\
$	3000	$ & $	229088.078	$ & $	-3.116	$ & $	0.922	$ \\
\hline
$	3100	$ & $	229085.208	$ & $	-2.870	$ & $	0.921	$ \\
$	3200	$ & $	229082.558	$ & $	-2.650	$ & $	0.923	$ \\
$	3300	$ & $	229080.095	$ & $	-2.463	$ & $	0.929	$ \\
$	3400	$ & $	229077.802	$ & $	-2.292	$ & $	0.931	$ \\
$	3500	$ & $	229075.679	$ & $	-2.123	$ & $	0.926	$ \\
$	3600	$ & $	229073.696	$ & $	-1.982	$ & $	0.934	$ \\
$	3700	$ & $	229071.831	$ & $	-1.865	$ & $	0.941	$ \\
$	3800	$ & $	229070.095	$ & $	-1.736	$ & $	0.931	$ \\
$	3900	$ & $	229068.457	$ & $	-1.637	$ & $	0.943	$ \\
$	4000	$ & $	229066.931	$ & $	-1.526	$ & $	0.932	$ \\
\hline
$	4100	$ & $	229065.482	$ & $	-1.449	$ & $	0.949	$ \\
$	4200	$ & $	229064.122	$ & $	-1.360	$ & $	0.939	$ \\
$	4300	$ & $	229062.837	$ & $	-1.285	$ & $	0.945	$ \\
$	4400	$ & $	229061.628	$ & $	-1.209	$ & $	0.941	$ \\
$	4500	$ & $	229060.484	$ & $	-1.144	$ & $	0.946	$ \\
$	4600	$ & $	229059.392	$ & $	-1.092	$ & $	0.955	$ \\
$	4700	$ & $	229058.366	$ & $	-1.026	$ & $	0.939	$ \\
$	4800	$ & $	229057.387	$ & $	-0.979	$ & $	0.955	$ \\
$	4900	$ & $	229056.459	$ & $	-0.928	$ & $	0.947	$ \\
$	5000	$ & $	229055.573	$ & $	-0.887	$ & $	0.956	$ \\
\hline
\end{tabular}%
~\qquad\qquad\qquad~%
\begin{tabular}{|c||c|c|c|}
\hline
$n$ & $C_n$ & $d_n$ & $r_n$  \\
\hline\hline
$	5100	$ & $	229054.729	$ & $	-0.843	$ & $	0.951	$ \\
$	5200	$ & $	229053.929	$ & $	-0.801	$ & $	0.949	$ \\
$	5300	$ & $	229053.160	$ & $	-0.768	$ & $	0.960	$ \\
$	5400	$ & $	229052.428	$ & $	-0.732	$ & $	0.953	$ \\
$	5500	$ & $	229051.730	$ & $	-0.698	$ & $	0.952	$ \\
$	5600	$ & $	229051.059	$ & $	-0.672	$ & $	0.963	$ \\
$	5700	$ & $	229050.420	$ & $	-0.638	$ & $	0.950	$ \\
$	5800	$ & $	229049.805	$ & $	-0.615	$ & $	0.964	$ \\
$	5900	$ & $	229049.217	$ & $	-0.588	$ & $	0.956	$ \\
$	6000	$ & $	229048.654	$ & $	-0.563	$ & $	0.958	$ \\
\hline
$	6100	$ & $	229048.111	$ & $	-0.543	$ & $	0.965	$ \\
$	6200	$ & $	229047.591	$ & $	-0.520	$ & $	0.956	$ \\
$	6300	$ & $	229047.090	$ & $	-0.501	$ & $	0.965	$ \\
$	6400	$ & $	229046.610	$ & $	-0.480	$ & $	0.958	$ \\
$	6500	$ & $	229046.146	$ & $	-0.464	$ & $	0.965	$ \\
$	6600	$ & $	229045.700	$ & $	-0.446	$ & $	0.963	$ \\
$	6700	$ & $	229045.270	$ & $	-0.430	$ & $	0.963	$ \\
$	6800	$ & $	229044.855	$ & $	-0.415	$ & $	0.965	$ \\
$	6900	$ & $	229044.455	$ & $	-0.400	$ & $	0.965	$ \\
$	7000	$ & $	229044.070	$ & $	-0.385	$ & $	0.962	$ \\
\hline
$	7100	$ & $	229043.698	$ & $	-0.373	$ & $	0.968	$ \\
$	7200	$ & $	229043.339	$ & $	-0.359	$ & $	0.963	$ \\
$	7300	$ & $	229042.990	$ & $	-0.349	$ & $	0.972	$ \\
$	7400	$ & $	229042.654	$ & $	-0.337	$ & $	0.966	$ \\
$	7500	$ & $	229042.329	$ & $	-0.325	$ & $	0.966	$ \\
$	7600	$ & $	229042.014	$ & $	-0.314	$ & $	0.967	$ \\
$	7700	$ & $	229041.709	$ & $	-0.306	$ & $	0.972	$ \\
$	7800	$ & $	229041.413	$ & $	-0.296	$ & $	0.968	$ \\
$	7900	$ & $	229041.127	$ & $	-0.286	$ & $	0.966	$ \\
$	8000	$ & $	229040.850	$ & $	-0.277	$ & $	0.969	$ \\
\hline
$	8100	$ & $	229040.580	$ & $	-0.270	$ & $	0.975	$ \\
$	8200	$ & $	229040.320	$ & $	-0.261	$ & $	0.965	$ \\
$	8300	$ & $	229040.067	$ & $	-0.253	$ & $	0.971	$ \\
$	8400	$ & $	229039.821	$ & $	-0.246	$ & $	0.973	$ \\
$	8500	$ & $	229039.582	$ & $	-0.239	$ & $	0.969	$ \\
$	8600	$ & $	229039.349	$ & $	-0.233	$ & $	0.978	$ \\
$	8700	$ & $	229039.124	$ & $	-0.225	$ & $	0.966	$ \\
$	8800	$ & $	229038.905	$ & $	-0.219	$ & $	0.972	$ \\
$	8900	$ & $	229038.692	$ & $	-0.213	$ & $	0.974	$ \\
$	9000	$ & $	229038.485	$ & $	-0.207	$ & $	0.970	$ \\
\hline
$	9100	$ & $	229038.282	$ & $	-0.203	$ & $	0.981	$ \\
$	9200	$ & $	229038.086	$ & $	-0.196	$ & $	0.965	$ \\
$	9300	$ & $	229037.895	$ & $	-0.191	$ & $	0.977	$ \\
$	9400	$ & $	229037.708	$ & $	-0.187	$ & $	0.975	$ \\
$	9500	$ & $	229037.527	$ & $	-0.181	$ & $	0.970	$ \\
$	9600	$ & $	229037.350	$ & $	-0.177	$ & $	0.977	$ \\
$	9700	$ & $	229037.178	$ & $	-0.172	$ & $	0.971	$ \\
$	9800	$ & $	229037.010	$ & $	-0.168	$ & $	0.979	$ \\
$	9900	$ & $	229036.846	$ & $	-0.164	$ & $	0.975	$ \\
$	10000	$ & $	229036.686	$ & $	-0.160	$ & $	0.974	$ \\
\hline
\end{tabular}
}
\caption{Calculated values of $C_n=n^5/\tn$, and associated differences $d_n=C_n-C_{n-100}$ and ratios $r_n=d_n/d_{n-100}$; see Section \ref{ss:as} for more details.}
\label{t:C}
\end{center}
\end{table}

\subsection{Distributions}\label{ss:tnd}

Looking at Figures \ref{f:tnd_100_200} and \ref{f:tnd_10000}, which graph the numbers $\td_{100}$, $\td_{200}$ and $\td_{10000}$ (for all allowable values of $d$), 
one sees a general shape emerging.  
In fact, by plotting scaled graphs of $\tnd$, one sees that the shapes are essentially identical, up to scale; several such graphs can be seen at \cite{web}, including some animations.  This means that in principle one might be able to estimate the value of (say) $t_{30000}$ by interpolating the shape of the $\td_{30000}$ curve from the $\td_{10000}$ or $\td_{20000}$ curve (for which we have full data), and calculating a few of the maximum values of $\td_{30000}$ to obtain the scaling factor.  
Table \ref{t:rhon} gives the maximum values of $\td_{30000}$ and $\td_{40000}$, which occur at $d=7140$ and $d=9520$, respectively.
The maximum value of $\tnd$ seems to regularly occur when $d\approx0.238n\approx\frac n{4.2}$.  We do not currently know the significance of this number.  We will come back to this point in Section~\ref{ss:max}.

Since there are so many points, Figure \ref{f:tnd_10000} appears to show a number of continuous curves.  The dark blue ``curve'' plots the sequence $\td_{10000}$ ($1667\leq d\leq 3332$), while the orange ``curve'' is the discrete derivative of this sequence: i.e., the values of $\td_{10000}-{}^{d-1}t_{10000}$.  (The meaning of the light blue curve will be explained in Section \ref{ss:ak}.)  Thus, the orange curve is (an approximation to) the derivative of the blue curve.  One may see that although the orange curve appears to be continuous, there is a sharp corner just after the maximum slope (which occurs at around $d=2073$).  

Similarly, Figures \ref{f:tnd_50_100} and \ref{f:tnd_1000} graph the numbers ${}^{50}\tn$, ${}^{100}\tn$ and ${}^{1000}\tn$ (for all allowable values of~$n$), 
and again one sees a general shape emerging.  Again Figure \ref{f:tnd_1000} shows the discrete derivative of the ${}^{1000}\tn$ sequence (in orange); although it is not as easy to see, the orange curve also has a sharp corner, this time just before the minimum slope (which occurs at around $n\approx4830$).

The data used to create Figure \ref{f:tnd_1000} also allows one to calculate
\[
{}^{1000}t = \sum_{n=3003}^{6000}{}^{1000}\tn = 25728695195597,
\]
which agrees with the largest computed value given by Kurz in \cite[Table 2]{Kurz2007}.  Some further values of~$\td$ (up to $d=2000$) are given in Table \ref{t:rhod}.


\subsection{Scratching the surface: maximal diameter}\label{ss:surface}

Tables \ref{t:tnd_50} and \ref{t:tnd_200} give values of $\tnd$, and a couple of simple patterns seem to emerge when looking at the right-most entry of each row, corresponding to the maximum $d$ for a fixed $n$.  Specifically, we see the three sequences
\begin{align}\label{e:3seq}
&{\blue1},{\blue1},{\blue2},{\blue2},{\blue3},{\blue3},\ldots , 
&&\textcolor{brown}{0},\textcolor{brown}{1},\textcolor{brown}{2},\textcolor{brown}{3},\textcolor{brown}{4},\textcolor{brown}{5},\ldots ,
&& {\green0},{\green1},{\green3},{\green6},{\green10},{\green13},{\green16},{\green19},{\green22},{\green25},{\green28},{\green31},{\green34},\ldots,
\intertext{corresponding to $n\equiv0$, $n\equiv1$ and $n\equiv2\Mod3$, respectively.  These are the numbers}
\nonumber&{}^dt_{3d+3}\ (d=1,2,3,\ldots), && {}^dt_{3d+4} \ (d=1,2,3,\ldots), && {}^dt_{3d+5} \ (d=1,2,3,\ldots),
\end{align}
respectively.
%
The apparent patterns in the first two sequences are obvious; it appears that
\[
{}^dt_{3d+3} = \lceil \tfrac d2 \rceil \AND {}^dt_{3d+4} = d-1 \qquad\text{for $d\geq1$.}
\]
There are simple explanations for these.  We give the details for the first, and sketch them for the second.  (We will consider the third sequence later.)  We begin with a simple observation:

\begin{lemma}\label{l:h}
Any altitude of an integer triangle is greater than $1/\sqrt2$.
\end{lemma}

\pf
Let the side-lengths of the triangle be $a\leq b\leq d$, and let $x$ be any of $a,b,d$.  Let $h$ be the altitude measured from a side of length $x$, and denote the area of the triangle by $A=xh/2$.  By Heron's formula, and keeping $a\leq b\leq d$ and $a+b\geq d+1$ in mind, we have
\begin{align*}
h^2 = \frac{4A^2}{x^2} \geq \frac{4A^2}{d^2} &= \frac4{d^2}\cdot\frac1{16}\cdot (a+b+d)\cdot(a+b-d)\cdot(a-b+d)\cdot(-a+b+d) \\
&\geq \frac1{4d^2}\cdot(2d+1)\cdot1\cdot a\cdot d > \frac{2ad^2}{4d^2} = \frac a2 \geq \frac12.
\end{align*}
The result follows.
\epf

\begin{rem}\label{r:h}
If the shortest side of an integer triangle is at least $2$, then it follows from the above proof that any altitude is greater than $1$.  Of course this is not true if the shortest side has length $1$, and if the altitude is measured from one of the two sides of length $d$; in fact it is easy to show that $h^2=1-1/4d^2$ in this case.  Since the latter is decreasing in $d$, and equal to $3/4$ when $d=1$, the lower bound of $1/\sqrt2\approx0.707$ in Lemma \ref{l:h} could be replaced by $\sqrt3/2\approx0.866$, though the latter is not a \emph{strict} lower bound.
\end{rem}

\begin{lemma}\label{l:3d+3}
For any $d\geq1$ we have ${}^dt_{3d+3} = \lceil \tfrac d2 \rceil$.
\end{lemma}

\pf
First consider an integer triangle with side-lengths $B\leq C\leq d$ satisfying ${B+C=d+1}$, and join two copies of this triangle together in the way shown in Figure \ref{f:3d+3}.  By Lemma \ref{l:h}, we can fold these towards each other (as in Figure \ref{f:T3}) until the tips are $1$ unit apart, thus obtaining an integer tetrahedron of diameter $d$ and perimeter $3d+3$.  There are $\lceil \tfrac d2 \rceil$ such triangles, and they give rise to pairwise-noncongruent tetrahedra.

Conversely, let $T$ be an arbitrary tetrahedron with diameter $d$ and perimeter $n=3d+3$.  We can construct $T$ by folding, as in Figure \ref{f:T3}, assuming that $A=d$ and $B\leq C$.  To ensure that $n=3d+3$ we must of course have $B+C=b+c=d+1$ and $a=1$.  To complete the proof, we must show that~$B=c$, as then also $C=b$.  Aiming for a contradiction, suppose instead that $B\not=c$, and consider the points~$P$ and~$Q$ as shown in Figure \ref{f:T3}.  These are on the circles $x^2+y^2=B^2$ and $x^2+y^2=c^2$, respectively.  Since $B$ and $c$ differ by at least $1$ (as they are distinct positive integers), it follows that $|PQ|\geq1=a$.  But this contradicts $|PQ|<a<|PR|$ from the proof of Proposition \ref{p:Tn}.
\epf

\begin{figure}[ht]
\begin{center}
\begin{tikzpicture}[scale=.6,inner sep=1.0]
\nc\AAA{10}
\nc\BBB4
\nc\CCC7
\nc\bbb7
\nc\ccc4
\begin{scope}[shift={(0,0)}]
\coordinate (P) at (0,0);
\coordinate (R) at (\AAA,0);
\coordinate (S) at ({(\AAA^2+\BBB^2-\CCC^2)/(2*\AAA)},{sqrt(\BBB^2-((\AAA^2+\BBB^2-\CCC^2)/(2*\AAA))^2)});
\coordinate (Q) at ({(\AAA^2+\ccc^2-\bbb^2)/(2*\AAA)},{-sqrt(\ccc^2-((\AAA^2+\ccc^2-\bbb^2)/(2*\AAA))^2)});
\coordinate (Q') at ({(\AAA^2+\ccc^2-\bbb^2)/(2*\AAA)},{sqrt(\ccc^2-((\AAA^2+\ccc^2-\bbb^2)/(2*\AAA))^2)});
\coordinate (Q'') at (4.5,3.5);
\fill[blue!20] (P) -- (Q) --(R); \fill[red!20] (R) -- (S) -- (P); 
\draw[blue,ultra thick] (P) --node[circle,fill=white]{$B$} (Q) --node[circle,fill=white]{$C$} (R); \draw[red,ultra thick] (R) --node[circle,fill=white]{$C$} (S) --node[circle,fill=white]{$B$} (P); 
\draw[->] (-1,0)--(11,0); \node () at (11,-.4) {\footnotesize $x$};
\draw[->] (0,-1)--(0,3.5); \node () at (-.4,3.5) {\footnotesize $y$};
\draw[ultra thick] (P) --node[circle,fill=white]{$d$} (R);
\end{scope}
\end{tikzpicture}
\caption{Creating an integer tetrahedron of diameter $d$ and perimeter $n=3d+3$ by folding congruent triangles (with $B+C=d+1$); see the proof of Lemma \ref{l:3d+3} for more details, and cf.~Figure \ref{f:T3}.}
\label{f:3d+3}
\end{center}
\end{figure}

\begin{lemma}\label{l:3d+4}
For any $d\geq1$ we have ${}^dt_{3d+4} = d-1$.
\end{lemma}

\pf[\bf Sketch of proof.]
Again we must analyse the pairs of triangles $(B,C,d)$ and $(b,c,d)$ that can be folded to create appropriate tetrahedra, as in Figure \ref{f:T3} (with $A=d$).  Up to symmetry, this time we have either
\bena
\item \label{3d+4a} $B+C=b+c=d+1$ and $a=2$, or else
\item \label{3d+4b} $B+C=d+1$, $b+c=d+2$ and $a=1$.
\een
We first note that case \ref{3d+4b} never actually occurs.  Indeed, here, since $B+C\not=b+c$, we must either have $B\not=c$ or $C\not=b$ (or both).  Considering circles, as in the proof of Lemma \ref{l:3d+3}, we see that $|PQ|\geq 1=a$, so the folding procedure cannot be carried out.

This leaves us to consider case \ref{3d+4a}, and we assume without loss of generality that $B\leq c$.  Again considering circles, $B$ and $c$ cannot differ by more than $1$, so we must have either
\begin{enumerate}[label=\textup{(a\arabic*)},leftmargin=9mm]
\item \label{a1} $c=B$ (and $b=C$), or else
\item \label{a2} $c=B+1$ (and $b=C-1$).
\een
We consider these in turn, and show that there are $\lceil \tfrac d2 \rceil-1$ and $\lfloor \tfrac d2 \rfloor$ tetrahedra in each case.  Since these sum to $d-1$, this will complete the proof.

\pfitem{\ref{a1}}  By symmetry, we may also assume that $B\leq C(=d+1-B)$.  So $1\leq B\leq \lceil \tfrac d2 \rceil$.  Using Remark~\ref{r:h}, we see that the folding procedure can be carried out in every case except for $B=1$.

\pfitem{\ref{a2}}  We again have $1\leq B\leq \lceil \tfrac d2 \rceil$, and this time the folding procedure can be carried out in each case.  However, when $d$ is odd, the $B= \lceil \tfrac d2 \rceil-1$ and $B= \lceil \tfrac d2 \rceil$ cases produce congruent tetrahedra.  So we obtain $\tfrac d2$ tetrahedra when $d$ is even, and $ \lceil \tfrac d2 \rceil-1$ when $d$ is odd.  In both cases, this is equal to $\lfloor \tfrac d2 \rfloor$.
\epf

The behaviour of the third sequence in \eqref{e:3seq} appears to be rather more complex.  For convenience in the following discussion, we will write ${u_d = {}^dt_{3d+5}}$ for $d\geq1$.  While these numbers seem to quickly stabilise into an arithmetic progression, $u_d = 3d - 5$ ($d\geq5$), this only persists until $d=41$, where we see an interesting change in the sequence:
\[
\ldots,
{\green100},
{\green103},
{\green106},
{\green109},
{\green112},
{\green115},
{\green118},
{\red\bf115},
{\red\bf115},
{\red\bf116},
{\red\bf117},
{\red\bf118},
{\red\bf119},
{\red\bf121},
{\red\bf122},
{\red\bf124},
{\red\bf126},
{\red\bf127},
{\red\bf129},
\ldots
\]
The colours in the above lists are reflected in Tables \ref{t:tnd_50} and \ref{t:tnd_200}, and we note that $u_{41}={\green118}$ and $u_{42}={\red\bf115}$.
The sequence of differences $u_{d+1}-u_d$ ($d\geq42$) begins
\[
0, 1, 1, 1, 1, 2, 1, 2, 2, 1, 2, 2, 1, 2, 2, 2, 1, 2, 2, 2, 2, 2, 2, 1, 2, 2, 2, 2, 2, 2, 2, 2, 1, \ldots.
\]
Although this appears somewhat chaotic, further calculations show that the last difference of $1$ seems to occur at $d=1110$, with all subsequent differences being $2$.  See Figure \ref{f:3d+5}, which plots the differences $u_{d+1}-u_d$ (right), and also the ratios $u_d/d$ (left).  The above discussion suggests the possible formula
\begin{equation}\label{e:3d+5}
u_d = {}^dt_{3d+5} = 2d+15 \qquad\text{for $d\geq1111$.}
\end{equation}
We have verified this computationally up to $d=50,000,000$.  Interestingly, we obtained much better performance using the ``na\"ive'' algorithm described in Section \ref{ss:code} (see steps \ref{a:I}--\ref{a:IV}), as compared to our modification of Kurz's algorithm from \cite{Kurz2007}.  (For example, the latter calculated $u_{3000}=6015$ in around 2 hours, while the former took a fraction of a second; it calculated the $50,000,000$th term in around 80 minutes.)  The reason for this appears to be that our original algorithm directly creates and stores all tetrahedra with given dimensions, and then counts them; this led to serious memory issues in general, but for parameters $n,d$ for which $\tnd$ is relatively small (such as $n=3d+5$) this is actually an advantage.  As discussed in Section \ref{ss:code}, the Kurz-based algorithm moves through tuple space, counting or rejecting tuples as appropriate; it seems that when $n$ is small relative to $d$, many more tuples are rejected than counted, leading to longer running times.

We have not attempted to prove \eqref{e:3d+5}, but we expect this could be done (if it is true) in a similar way to Lemmas \ref{l:3d+3} and \ref{l:3d+4} above.  
Computations suggest that for $d\geq1111$, there are:
\bit
\item $\lfloor\tfrac d2\rfloor-1$ tetrahedra with $a=1$, 
\item $\lceil\tfrac {3d}2\rceil-3$ tetrahedra with $a=3$, and
\item $19$ ``sporadic'' tetrahedra with $a=2$,
\eit
and these sum to $2d+15$.
For smaller $d$ the number of sporadic tetrahedra is different; for example, it is $20$ for $205\leq d\leq1110$.  Figure \ref{f:sporadic} plots the number of such sporadic tetrahedra for $1\leq d\leq 2000$.

\begin{figure}[ht]
\begin{center}
\includegraphics[height=4.8cm]{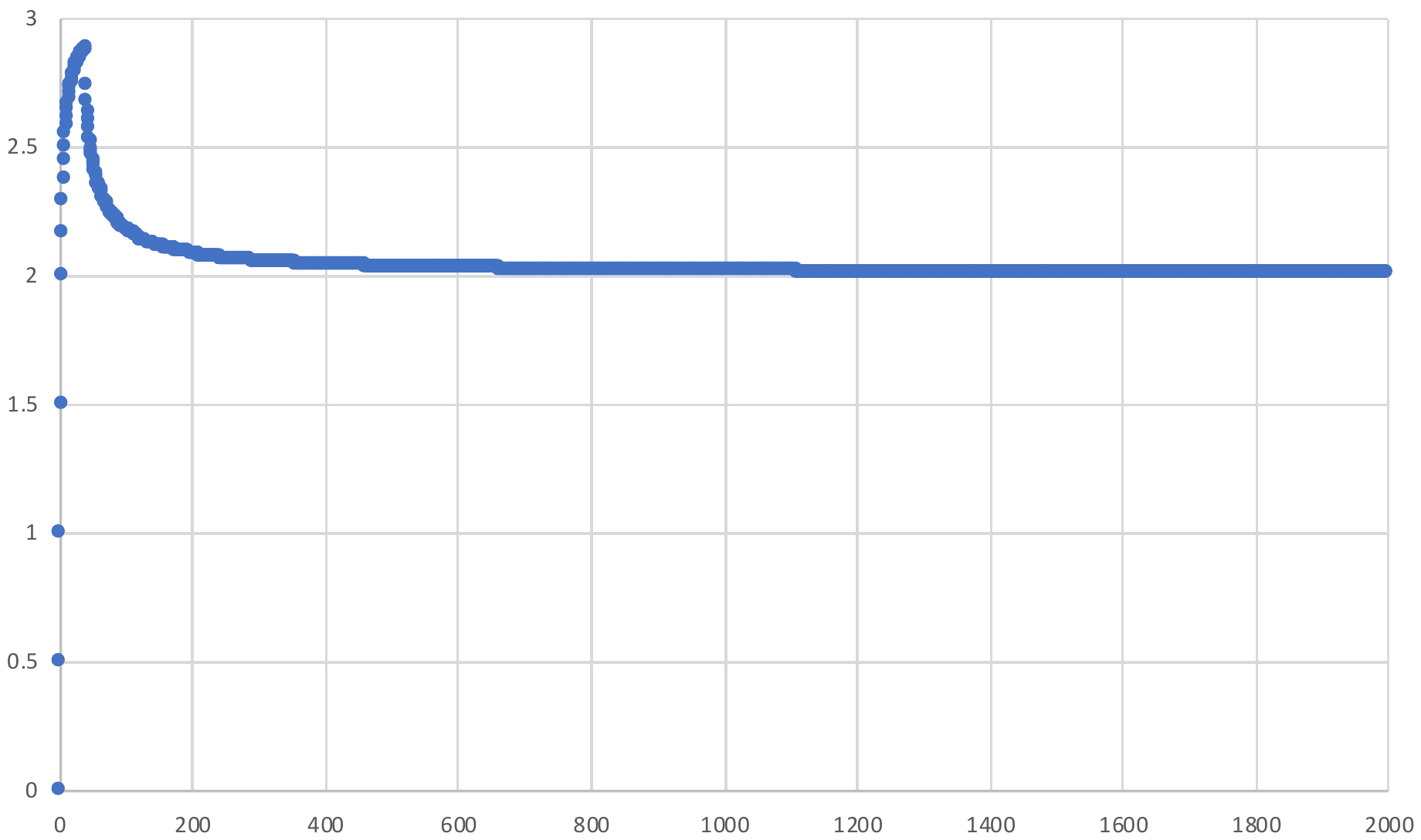}
\quad
\includegraphics[height=4.8cm]{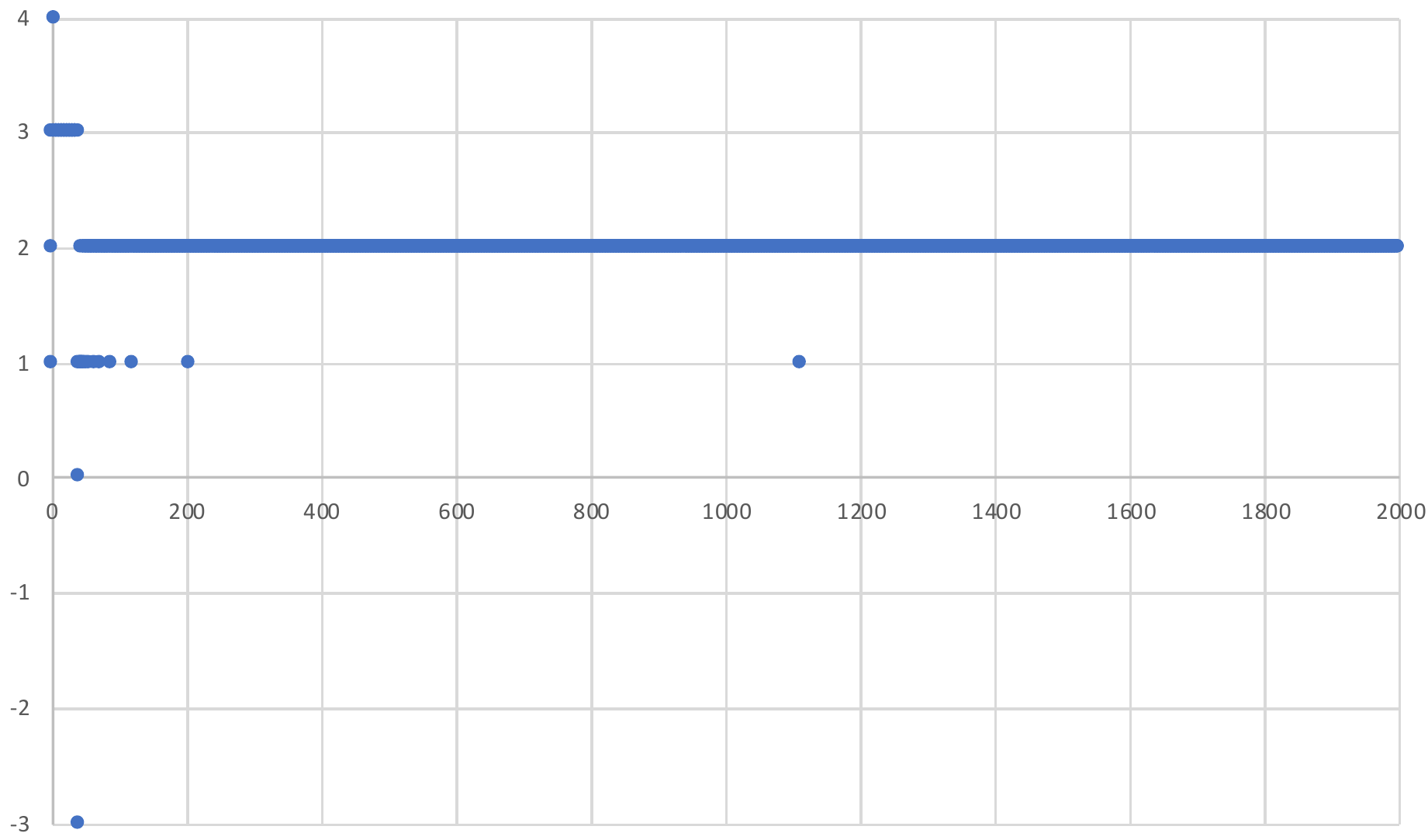}
\caption{Left and right: calculated values of $u_d/d$ and $u_{d+1}-u_d$, $1\leq d\leq2000$, where $u_d={}^dt_{3d+5}$.}
\label{f:3d+5}
\end{center}
\end{figure}

\begin{figure}[ht]
\begin{center}
\includegraphics[width=\textwidth]{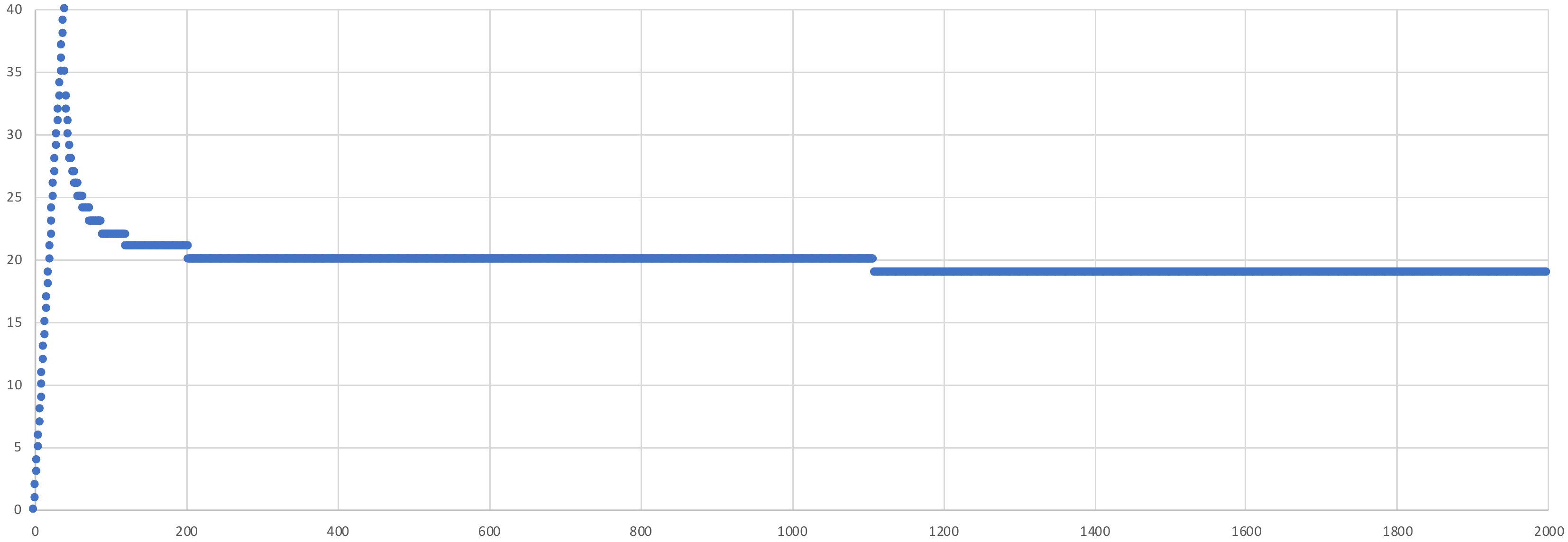}
\caption{The number of ``sporadic'' tetrahedra with diameter $d$ and perimeter $3d+5$, $1\leq d\leq2000$; see Section \ref{ss:surface} for more details.}
\label{f:sporadic}
\end{center}
\end{figure}

We have not attempted to systematically study the numbers ${}^dt_{3d+k}$ for (fixed) $k\geq6$, but Figure~\ref{f:3dk} gives some graphs of computed values for $k=3,4,\ldots,10$.  It appears that for fixed $k$, the sequence~${}^dt_{3d+k}$~($d=1,2,3,\ldots$) is eventually linear in $d$, but the ``pre-linear'' behaviour becomes more complex, and lasts longer as $k$ increases.  Taking $k=10$, for example, and writing $v_d={}^dt_{3d+10}$, Figure~\ref{f:3d+10} shows the ratios~$v_d/d$ and differences $v_{d+1}-v_d$ (cf.~Figure \ref{f:3d+5}, which does the same for $k=5$).  Both appear to approach a limit of $10$, strongly suggesting that eventually $v_d=10d+l$ for some $l$.  Calculations suggest that $l=2470$; this value of $l$ holds for~$2,000,000\leq d\leq20,000,000$.  But we note that $l=2471$ for $d=1,000,000$, so that the ``pre-linear'' behaviour persists beyond the $1,000,000$th term.


\begin{figure}[ht]
\begin{center}
\includegraphics[width=\textwidth]{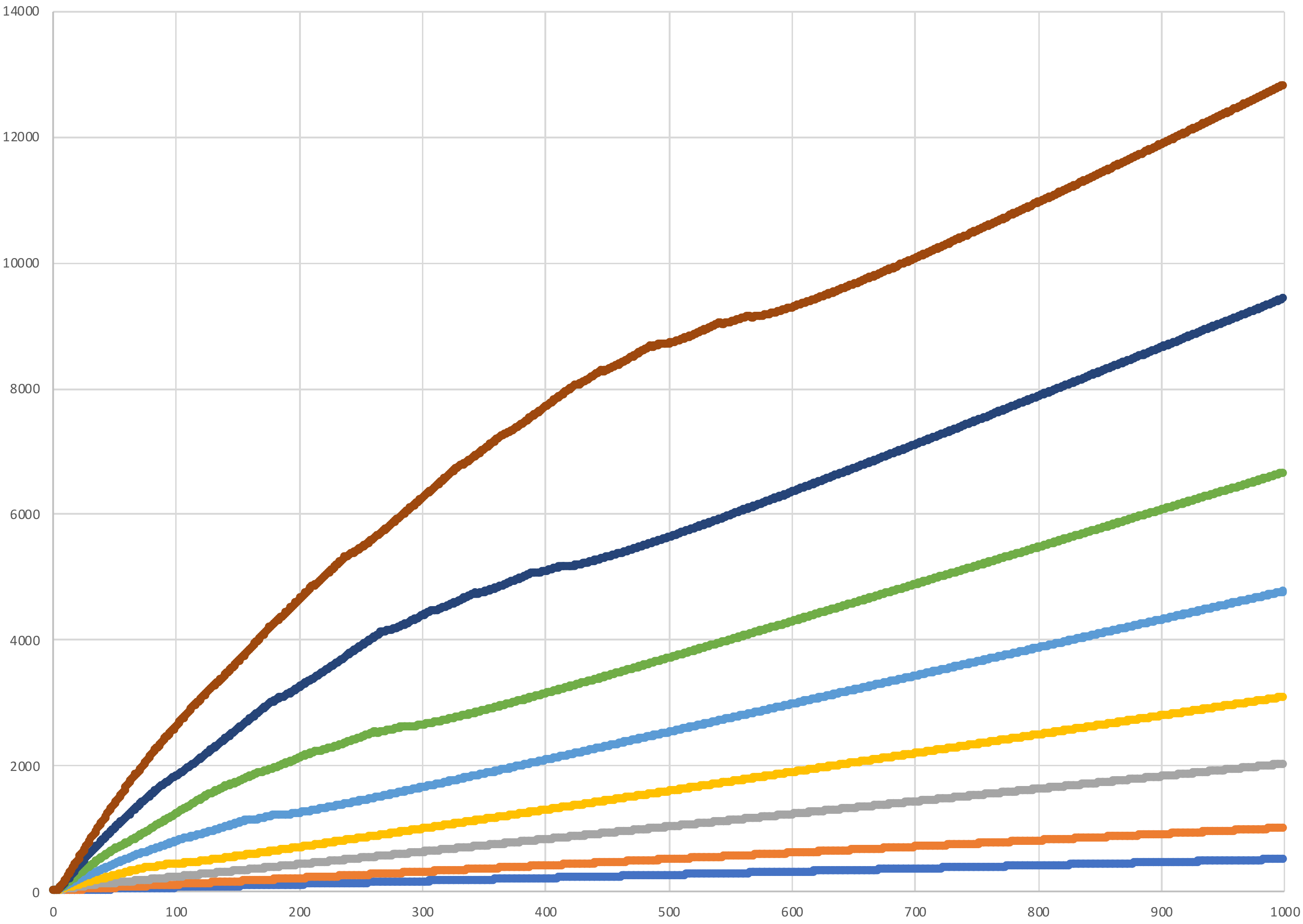}
\caption{Calculated values of ${}^dt_{3d+k}$, $1\leq d\leq1000$, for $k=3,4,\ldots,10$ (bottom to top).}
\label{f:3dk}
\end{center}
\end{figure}

\begin{figure}[ht]
\begin{center}
\includegraphics[height=5.1cm]{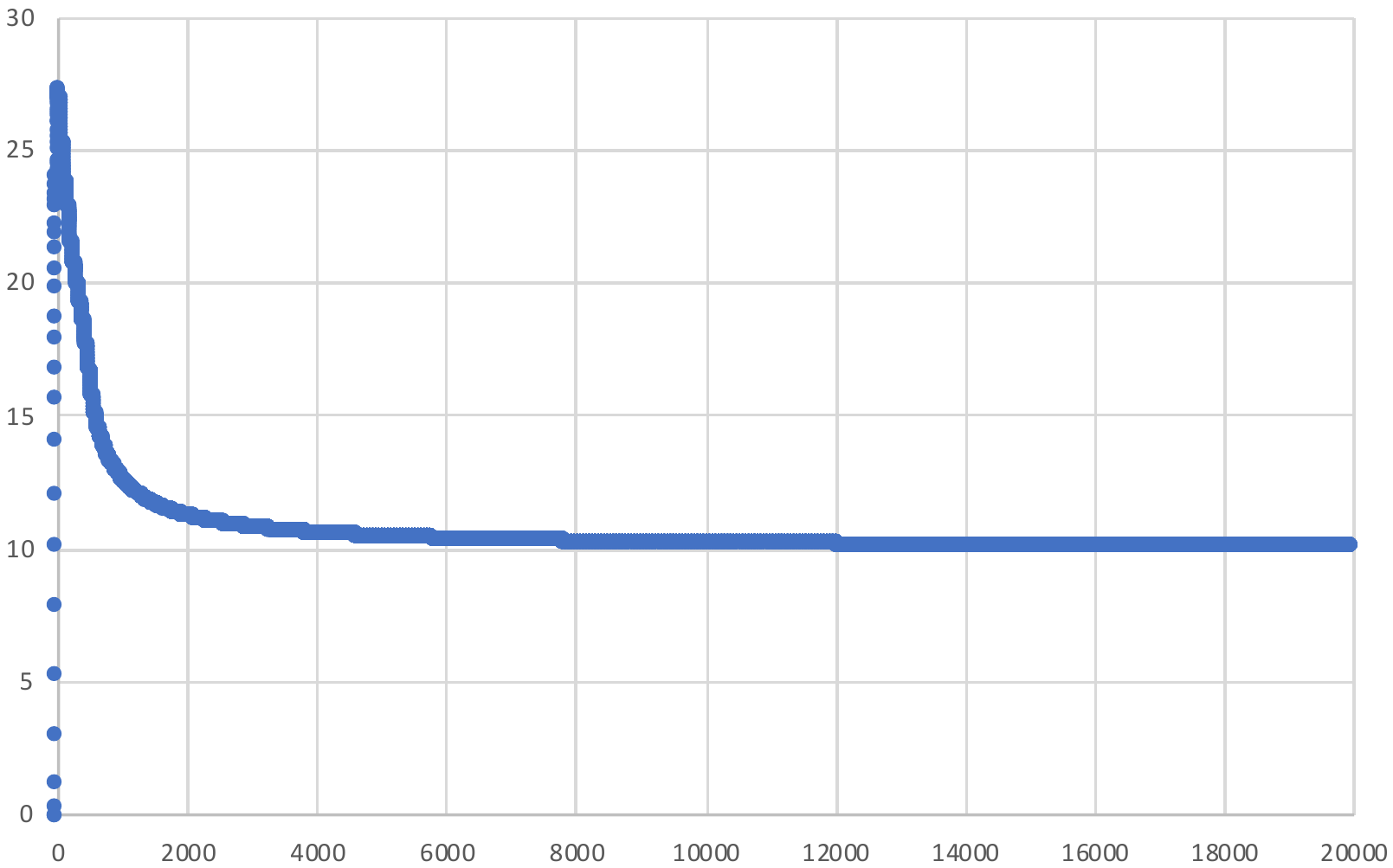}
\quad
\includegraphics[height=5.1cm]{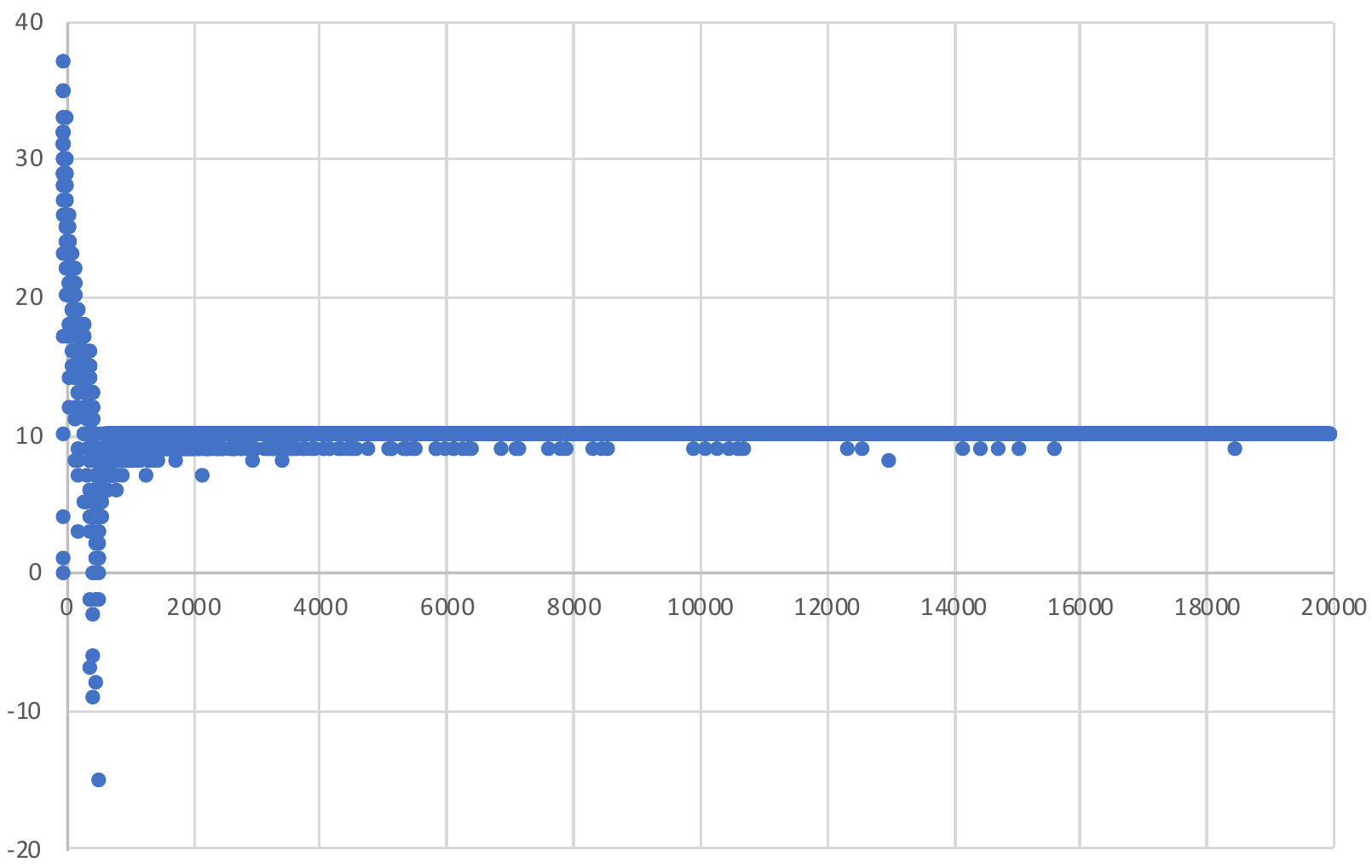}
\caption{Left and right: calculated values of $v_d/d$ and $v_{d+1}-v_d$, $1\leq d\leq20000$, where $v_d={}^dt_{3d+10}$.}
\label{f:3d+10}
\end{center}
\end{figure}

\subsection{Initial segments: a glimpse of order and hope?}\label{ss:ak}

In the previous section we looked at the very top parts of the columns in the~$\tnd$ data, as shown in Tables \ref{t:tnd_50} and \ref{t:tnd_200}.  One of the most interesting/promising observations arises when one looks at the \emph{bottom} parts of these columns.  

Examining consecutive columns, one sees that the first few values at the bottom of one column are present in the next.  These are the red entries in the lower parts of Tables \ref{t:tnd_50} and \ref{t:tnd_200}, specifically the values $\tnd$ with $n$ approximately ranging from $5d$ to $6d$.  Note that each column seems to add an extra number to the ``Stable Column Sequence'' (as we will call it), but that no extra number is added from column $d=28$ to $d=29$; this is indicated by green in Table \ref{t:tnd_200}.  Although we do not have a complete explanation for why this happens, it appears that for fixed $d$, the last value of $n$ for which~$\tnd$ is not this ``stable'' value is $\lceil 5.035d\rceil$; note that
\[
\lceil5.035\cdot28\rceil = \lceil140.98\rceil = 141 \AND \lceil5.035\cdot29\rceil=\lceil146.015\rceil = 147,
\]
so that column $d=29$ adds an additional $6$ ``unstable'' $\tnd$ values from column $d=28$.

In any case, the resulting Stable Column Sequence begins:
\begin{equation}
\label{e:ak} 1,1,3,6,11,18,31,47,72,105,149,206,281,372,487,627,796,997,1237,1516,1843,2220,2653,3147,\ldots
\end{equation}
\newpage
\noindent This sequence does not appear on the OEIS.  The first few values are easy to understand:
\bit
\item for $d\geq1$, ${}^dt_{6d}=1$ counts only $\tup dddddd$, the equilateral tetrahedron, 
\item for $d\geq2$, ${}^dt_{6d-1}=1$ counts only $\tup ddddd{d-1}$, 
\item for $d\geq3$, ${}^dt_{6d-2}=3$ counts only $\tup ddddd{d-2}$, $\tup dddd{d-1}{d-1}$ and $\tup ddd{d-1}d{d-1}$.
\eit
One could similarly explain other values, with ad hoc arguments: e.g., ${}^dt_{6d-3}=6$ for $d\geq4$.

If we denote the sequence \eqref{e:ak} by $a_k$ ($k=0,1,2,\ldots$), then calculating $\left(\frac{a_{k+1}}{a_k}-1\right)\times k$ suggests that~$a_k$ is quartic in $k$.  Further experimentation quickly suggests the leading coefficient is $\frac1{96}$.  Additional analysis suggests that $a_k-\frac{k^2}{96}$ is quadratic with leading coefficient $\frac7{16}=\frac{42}{96}$.  After examining $a_k-\frac{k^4}{96}-\frac{7k^2}{16}$, it eventually appears that we have the exact formula
\begin{equation}\label{e:ak2}
a_k = \frac{k^4+42k^2+b_k}{96} \qquad\text{where $b_0=96$, and for $k\geq1$:} \qquad b_k = 
\begin{cases}
192	&\text{if $k\equiv0\Mod{12}$}\\
53	&\text{if $k\equiv1,5,7,11\Mod{12}$}\\
104	&\text{if $k\equiv2,10\Mod{12}$}\\
117	&\text{if $k\equiv3,9\Mod{12}$}\\
128	&\text{if $k\equiv4,8\Mod{12}$}\\
168	&\text{if $k\equiv6\Mod{12}$.}
\end{cases}
\end{equation}
Equivalently, for $k\geq1$, $a_k$ is the nearest integer to
\[
\frac{k^4+42k^2+148}{96} \qquad\text{or}\qquad \frac{k^4+42k^2+85}{96},
\]
for even and odd $k$, respectively.  (The values of $148$ and $85$ are simply the average of the maximum and minimum values of $b_k$ for even and odd $k$, respectively.  A single ``nearest integer formula'' cannot be given, unfortunately, because the distance between the maximum and minimum values of $b_k$ for arbitrary $k$ (i.e., $192-53=139$) is greater than the denominator, $96$.)  We do not currently know if there is any significance in the exact values of the numbers $b_k\in\{53,104,117,128,168,192\}$.

\begin{con}\label{c:ak}
For suitably small $k$ (approximately $0\leq k<d$), we have ${}^dt_{6d-k} = a_k$, where $a_k$ is defined in \eqref{e:ak2}.  (Note then that for such $k$, ${}^dt_{6d-k}$ would depend only on $k$, and not on $d$.)
\end{con}

The formula \eqref{e:ak2} can be used to calculate many more values of $a_k$.  Figure \ref{f:abk} plots these, together with the associated values of $b_k$; the periodic nature of the latter can be readily seen in the graph.

\newpage

The existence of the Stable Column Sequence implies that there are six ``Stable Row Sequences'', one for each residue of $n$, modulo $6$.  So the $n$th row of the $\tnd$ table, for suitably large $n$, begins:
\bit
\item $18   ,   206 ,     997  ,   3147  ,   7736 ,   16168 ,   30171 ,   51797,    83422,\ldots$ when $n\equiv0\Mod6$,
\item $11  ,    149  ,    796 ,    2653  ,   6747,    14427  ,  27368,    47567,    77347,\ldots$ when $n\equiv1\Mod6$,
\item $6   ,   105     , 627  ,   2220    , 5856 ,   12831    ,24765 ,   43602  ,  71610,\ldots$ when $n\equiv2\Mod6$,
\item $3  ,     72      ,487   ,  1843   ,  5057  ,  11372   , 22353  ,  39891   , 66199,\ldots$ when $n\equiv3\Mod6$,
\item $1 ,      47      ,372   ,  1516   ,  4343  ,  10041   , 20122  ,  36422   , 61101,\ldots$ when $n\equiv4\Mod6$,
\item $1,       31      ,281   ,  1237   ,  3710  ,   8833    ,18065   , 33187   , 56306,\ldots$ when $n\equiv5\Mod6$.
\eit
These can all be seen in Table \ref{t:tnd_200}, and they are all of course subsequences of the Stable Column Sequence.
None of these six (sub)sequences appear on the OEIS.
Note that the (initial segments of the) Stable Row Sequences account for the left-most part of the graphs in Figures \ref{f:tnd_100_200} and~\ref{f:tnd_10000}; see especially Figure~\ref{f:tnd_10000}, which shows values of $a_k$ in light blue.  Similarly, the Stable Column Sequence accounts for the right-most part of the graphs in Figures \ref{f:tnd_50_100} and \ref{f:tnd_1000}; see especially Figure~\ref{f:tnd_1000}, which shows values of $a_k$ in light blue.

So we have seen that the first portion of a row, and the last portion of a column, in the $\tnd$ table appears to be governed by a certain quartic (indeed, biquadratic) polynomial.  Unfortunately, similar experimentation shows that the remainder of these rows/columns do not look at all like any polynomial.  Indeed, even considering small portions at the opposite end of each row/column does not seem to yield a polynomial pattern.

A possible approach to Conjecture \ref{c:ak} might be to show that any tuple $\tup AaBbCc$ with $d=A=\max\{A,B,C,a,b,c\}$, and with perimeter greater than (around) $5.035d$, automatically satisfies conditions \ref{T2}--\ref{T6} from Proposition \ref{p:Tn}.  
The hope is that such a tuple is ``close enough'' to the equilateral $\tup dddddd$ that it corresponds to a tetrahedron by default.  This would not immediately give the desired formula for $a_k$, but it would at least show that some kind of Stable Column Sequence does exist.

\begin{figure}[ht]
\begin{center}
\includegraphics[height=3.6cm]{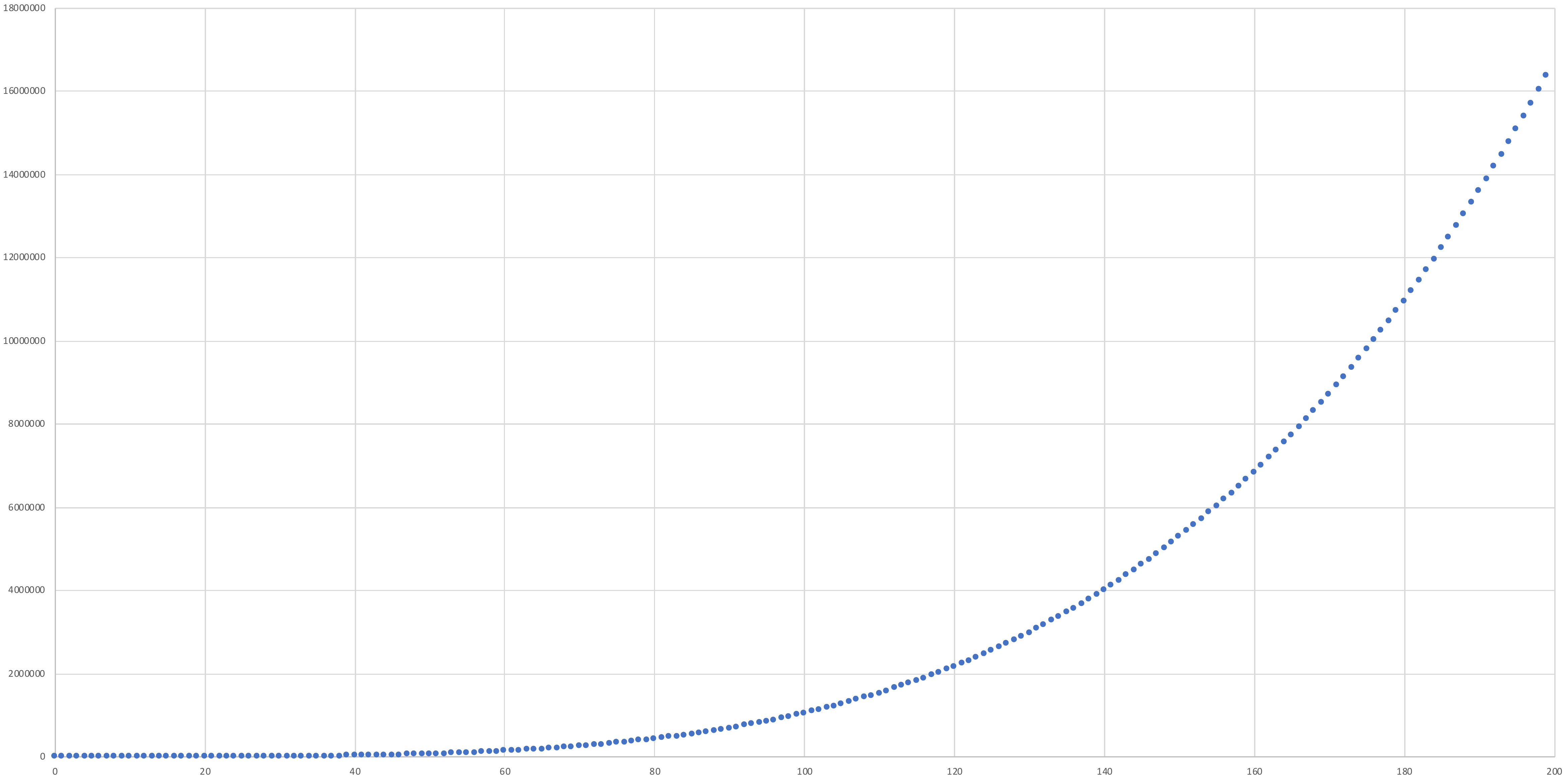}
\qquad
\includegraphics[height=3.6cm]{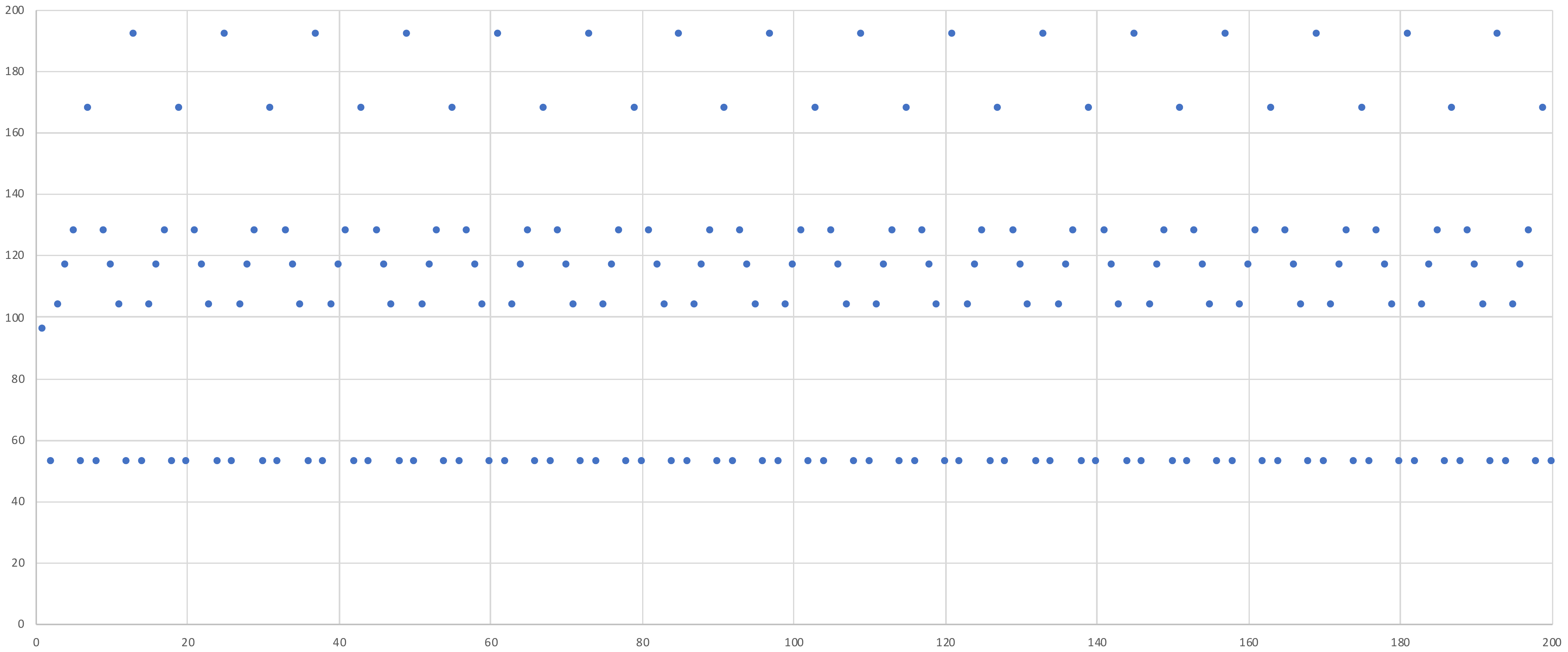}
\caption{Calculated values of $a_k$ (left) and $b_k$ (right), $1\leq k\leq200$, as defined in \eqref{e:ak2}.}
\label{f:abk}
\end{center}
\end{figure}

\subsection{Sums vs.~maximums: more on asymptotics}\label{ss:max}

We noted in Section \ref{ss:tnd} that one might hope to approximate values of $\tn$ (for large $n$) by interpolating the ${}^dt_m$ curve (for smaller $m$), and using the value of $\max_d\tnd$ as the scaling factor.  In this section we discuss a simpler attempt to estimate $\tn$ (and $\td$) using single values of $\tnd$.

To keep the following discussion manageable, it will be convenient to define a number of additional parameters:
\bit
\item For fixed $n$ we define $\mun = \max_d \tnd$, and the ratio $\rhon = \tn/\mun$.  \\  
Further, let $d^*(n)$ be the diameter corresponding to the maximum value of $\tnd$: i.e., $\mun={}^{d^*(n)}t_n$.
\item For fixed $d$ we define $\mud = \max_n \tnd$, and the ratio $\rhod = \td/\mud$.  \\
Further, let $n^*(d)$ be the perimeter corresponding to the maximum value of $\tnd$: i.e., $\mud={}^dt_{n^*(d)}$.
\eit
Computationally, it seems that 
\begin{equation}\label{e:d*n*}
d^*(n) \approx 0.238n\approx\frac{n}{4.2} \AND n^*(d) \approx 4.5d.
\end{equation}
We do not currently understand the significance of the number $0.238\approx\frac1{4.2}$.  However, $4.5d$ is of course almost exactly mid-way between the extreme values of $3d+3\leq n\leq6d$; cf.~Lemma \ref{l:d}.  These values of $d^*(n)$ and $n^*(d)$ can be seen by locating the peaks in Figures \ref{f:tnd_100_200}--\ref{f:tnd_1000}.

Figure \ref{f:rho} shows graphs of calculated values of $\rhon$ ($n\leq2000$) and $\rhod$ ($d\leq300$).  These graphs strongly suggest that (at least approximately) $\rhon$ is linear in $n$, and $\rhod$ is linear in $d$, so it seems worthwhile to look at the values of $\rhon/n$ and $\rhod/d$.  
Accordingly, 
\bit
\item Table \ref{t:rhon} shows calculated values of $\tn$, $\mun$, $\rhon$ and $\rhon/n$ (as well as $d^*(n)$), and
\item Table \ref{t:rhod} shows calculated values of $\td$, $\mud$, $\rhod$ and $\rhod/d$ (as well as $n^*(d)$).
\eit
Note that we have been able to calculate $\mun$ and $\mud$ for higher values of $n$ and $d$ than for $\tn$ and $\td$.  To obtain these 
we used the approximations for $d^*(n)$ and $n^*(d)$ from \eqref{e:d*n*}, and calculated enough values of $\tnd$ around the critical values, to ensure that we had found the peak of the curve.  For the round values of $n$ and $d$ we used, the approximations in \eqref{e:d*n*} were in fact exact, with a single exception: $n^*(1900)=8551$, whereas $4.5\cdot1900=8550$.  But here we note that
\bit
\item ${}^{1900}t_{8549} = 287666266084$,
\item ${}^{1900}t_{8550} = 287666506644 = {}^{1900}t_{8549} + 240560$,
\item ${}^{1900}t_{8551} = 287666521356 = {}^{1900}t_{8550} + 14712$,
\item ${}^{1900}t_{8552} = 287666235410 = {}^{1900}t_{8551} + 271234$,
\eit
meaning that ${}^{1900}t_{8550}$ is comparatively very close to the maximum value of ${}^{1900}\mu = {}^{1900}t_{8551}$.

The values displayed in Tables \ref{t:rhon} and \ref{t:rhod} strongly suggest that $\rhon/n$ and $\rhod/d$ tend to limits of around $0.07122$ and $1.1657$, respectively, leading to asymptotic expressions:
\begin{equation}\label{e:as_rho}
\rhon \sim 0.07122n \text{ \ as $n\to\infty$} \AND \rhod \sim 1.1657d \text{ \ as $d\to\infty$.}
\end{equation}
(Interestingly, $0.07122$ is quite close to $\frac1{14}$, although computational evidence suggests that the limit $\lim_{n\to\infty}\frac{\rhon}n$ is around $\frac1{14.041}$, and this can make quite a big difference in the calculations that follow.)

By definition, we have $\tn=\rhon\times\mun$ and $\td=\rhod\times\mud$.  Thus, if we could also obtain asymptotic expressions for $\mun$ and $\mud$, then these could be combined with \eqref{e:as_rho} to yield asymptotic expressions for~$\tn$ and $\td$ themselves (cf.~Conjecture \ref{c:as} and Remark \ref{r:Kurz_as}).

Calculations show that the ratios $\left(\frac{\mu_{n+1}}{\mun}-1\right)\times n$ and $\left(\frac{{}^{d+1}\mu}{\mud}-1\right)\times d$ both seem to approach $4$, suggesting that $\mun\sim \frac{n^4}D$ and $\mud\sim \frac{d^4}E$ for some constants $D$ and $E$.  Tables \ref{t:rhon} and \ref{t:rhod} also show values of $D_n=\frac{n^4}{\mun}$ and ${}^dD=\frac{d^4}{\mud}$, respectively.  (These tables also give values of $C_n = \frac{n^5}{\tn}$ and ${}^dC = \frac{d^5}{\td}$.)  These seem to approach limits of around $D_n\to16310.5$ and ${}^dD\to45.301$, respectively.


Putting all of the above together, we have (very approximate) candidates for asymptotic formulas:
\begin{align}
\label{e:as_n}
\tn &= \rhon\times\mun \sim 0.07122n \times \frac{n^4}{16310.5} \approx \frac{n^5}{229016} &&\text{as $n\to\infty$} \\[5mm]
\text{and} \qquad\qquad 
\label{e:as_d}
\td &= \rhod\times\mud \sim 1.1657d \times \frac{d^4}{45.301} \approx \frac{d^5}{38.86} &&\text{as $d\to\infty$.}
\end{align}
Note that \eqref{e:as_n} and \eqref{e:as_d} are reasonably close to Conjecture \ref{c:as} and \eqref{e:td_con}.  
We should stress, however, that in the above discussion we often rounded off certain numbers to seemingly-arbitrary precision.  Thus, the asymptotic expressions in \eqref{e:as_n} and \eqref{e:as_d} are very speculative, and are not meant to be exact predictions.
%

\begin{figure}[ht]
\begin{center}
\includegraphics[height=4.8cm]{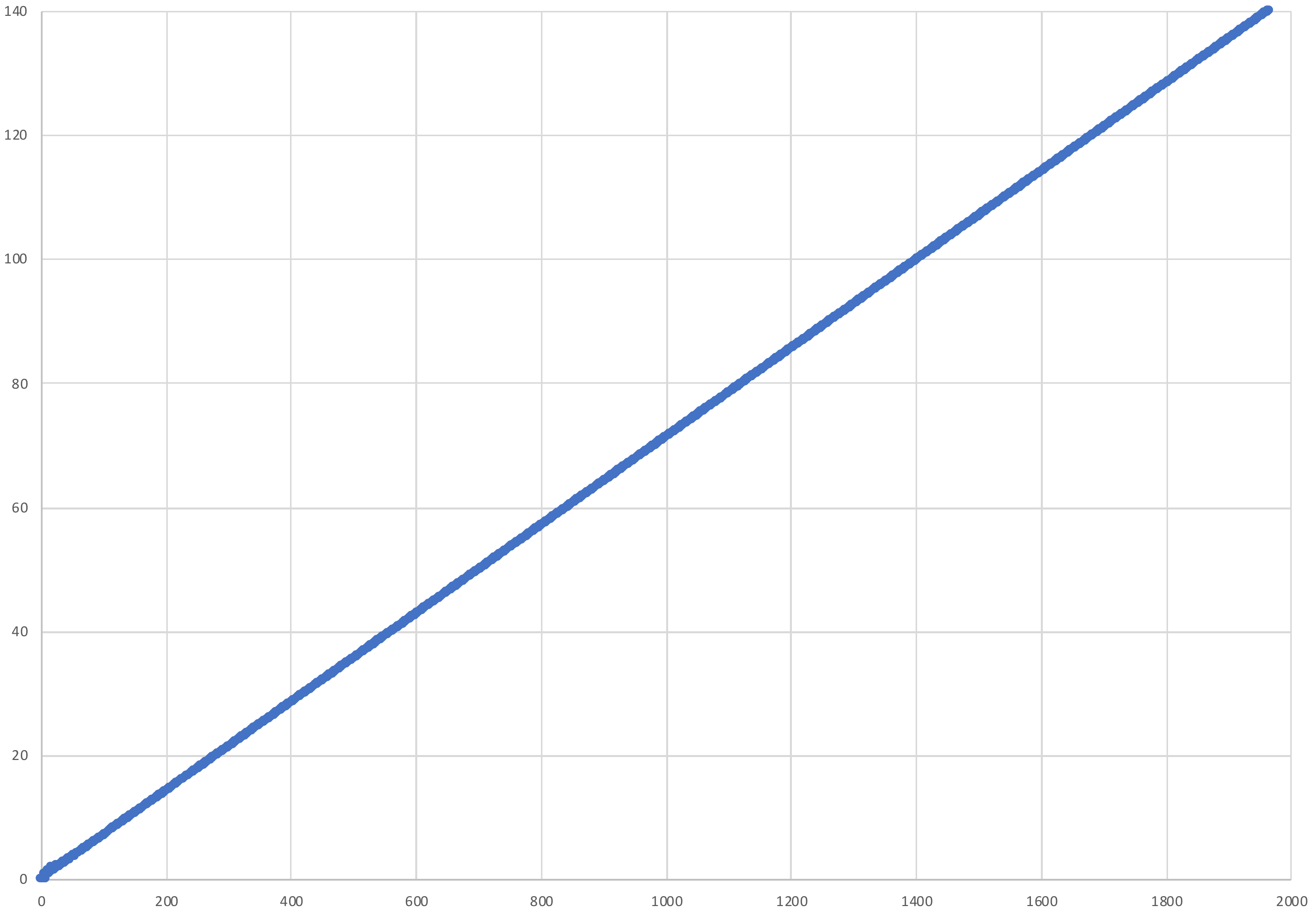}
\qquad\qquad
\includegraphics[height=4.8cm]{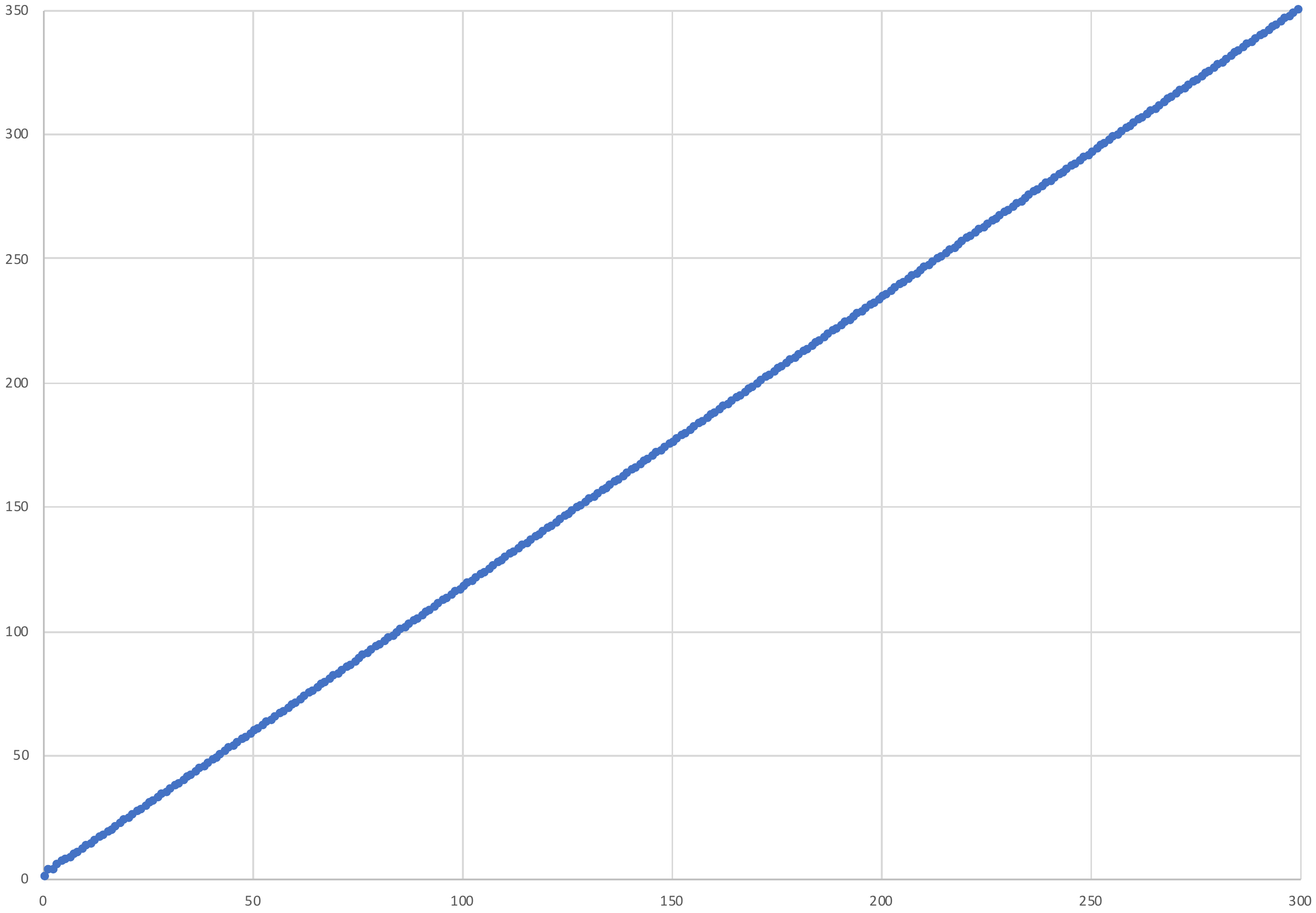}
\caption{Calculated values of $\rhon$ (left), $1\leq n\leq2000$, and $\rhod$ (right), $1\leq d\leq300$; see Section \ref{ss:max} for more details.}
\label{f:rho}
\end{center}
\end{figure}

\begin{table}[ht]
\begin{center}
{\footnotesize
\begin{tabular}{|r|r||r|r|r|r|r|r|}
\hline
\multicolumn{1}{|c|}{$n$} & \multicolumn{1}{|c||}{$d^*(n)$} & \multicolumn{1}{|c|}{$\tn$} & \multicolumn{1}{|c|}{$\mun$} & \multicolumn{1}{|c|}{$\rhon$} & \multicolumn{1}{|c|}{$\rhon/n$} & \multicolumn{1}{|c|}{$D_n$} & \multicolumn{1}{|c|}{$C_n$} \\
\hline\hline
$	1000	$ & $	238	$ & $	4360687860	$ & $	61261647	$ & $	71.181368	$ & $	0.071181368	$ & $	16323.42663	$ & $	229321.6190	$ \\
$	2000	$ & $	476	$ & $	139654346301	$ & $	980662960	$ & $	142.408097	$ & $	0.071204049	$ & $	16315.49335	$ & $	229137.1579	$ \\
$	3000	$ & $	714	$ & $	1060727392377	$ & $	4965253737	$ & $	213.630048	$ & $	0.071210016	$ & $	16313.36570	$ & $	229088.0784	$ \\
$	4000	$ & $	952	$ & $	4470309160343	$ & $	15693529026	$ & $	284.850473	$ & $	0.071212618	$ & $	16312.45589	$ & $	229066.9310	$ \\
$	5000	$ & $	1190	$ & $	13642977397892	$ & $	38315438250	$ & $	356.069982	$ & $	0.071213996	$ & $	16311.96271	$ & $	229055.5726	$ \\
$	6000	$ & $	1428	$ & $	33949118928429	$ & $	79452345303	$ & $	427.289072	$ & $	0.071214845	$ & $	16311.66450	$ & $	229048.6541	$ \\
$	7000	$ & $	1666	$ & $	73378891579018	$ & $	147197080983	$ & $	498.507790	$ & $	0.071215399	$ & $	16311.46477	$ & $	229044.0703	$ \\
$	8000	$ & $	1904	$ & $	143066182224551	$ & $	251113882084	$ & $	569.726297	$ & $	0.071215787	$ & $	16311.32443	$ & $	229040.8501	$ \\
$	9000	$ & $	2142	$ & $	257812568218126	$ & $	402238432945	$ & $	640.944641	$ & $	0.071216071	$ & $	16311.22106	$ & $	229038.4849	$ \\
$	10000	$ & $	2380	$ & $	436611276762080	$ & $	613077830990	$ & $	712.162885	$ & $	0.071216288	$ & $	16311.14272	$ & $	229036.6862	$ \\
\hline
$	11000	$ & $	2618	$ & $	703171140799375	$ & $	897610575535	$ & $	783.381079	$ & $	0.071216462	$ & $	16311.08233	$ & $	229035.2813	$ \\
$	12000	$ & $	2856	$ & $	1086440561201385	$ & $	1271286727639	$ & $	854.599153	$ & $	0.071216596	$ & $	16311.03318	$ & $	229034.1588	$ \\
$	13000	$ & $	3094	$ & $	1621131464654093	$ & $	1751027590016	$ & $	925.817202	$ & $	0.071216708	$ & $	16310.99371	$ & $	229033.2450	$ \\
$	14000	$ & $	3332	$ & $	2348243259709892	$ & $	2355226146934	$ & $	997.035152	$ & $	0.071216797	$ & $	16310.96022	$ & $	229032.4896	$ \\
$	15000	$ & $	3570	$ & $	3315586803814407	$ & $	3103746427345	$ & $	1068.253120	$ & $	0.071216875	$ & $	16310.93299	$ & $	229031.8562	$ \\
$	16000	$ & $	3808	$ & $	4578308347767338	$ & $	4017924323570	$ & $	1139.471025	$ & $	0.071216939	$ & $	16310.90950	$ & $	229031.3191	$ \\
$	17000	$ & $	4046	$ & $	6199413506818430	$ & $	5120566904917	$ & $	1210.688899	$ & $	0.071216994	$ & $	16310.88931	$ & $	229030.8589	$ \\
$	18000	$ & $	4284	$ & $	8250291212969114	$ & $	6435952631653	$ & $	1281.906764	$ & $	0.071217042	$ & $	16310.87207	$ & $	229030.4610	$ \\
$	19000	$ & $	4522	$ & $	10811237680756863	$ & $	7989831502801	$ & $	1353.124616	$ & $	0.071217085$ & $	16310.85711	$ & $	229030.1141	$ \\
$	20000	$ & $	4760	$ & $	13971980347794463	$ & $	9809425098449	$ & $	1424.342427	$ & $	0.071217121	$ & $	16310.84374	$ & $	229029.8097	$ \\
\hline\hline
$	30000	$ & $	7140	$ &  & $	49660446083302	$ & $		$ &  & $	16310.76770	$ & $		$ \\
\hline\hline
$	40000	$ & $	9520	$ &  & $	156951842767924	$ & $		$ &  & $	16310.78743	$ & $		$ \\
\hline
\end{tabular}%
}
\caption{Calculated values of $\tn$ and associated numbers defined in Section \ref{ss:max}.}
\label{t:rhon}
\end{center}
\end{table}

\begin{table}[ht]
\begin{center}
{\footnotesize
\begin{tabular}{|r|r||r|r|r|r|r|r|r|}
\hline
\multicolumn{1}{|c|}{$d$} & \multicolumn{1}{|c||}{$n^*(d)$} & \multicolumn{1}{|c|}{$\td$} & \multicolumn{1}{|c|}{$\mud$} & \multicolumn{1}{|c|}{$\rhod$} & \multicolumn{1}{|c|}{$\rhod/d$} & \multicolumn{1}{|c|}{${}^dD$} & \multicolumn{1}{|c|}{${}^dC$} \\
\hline\hline
$	100	$ & $	450	$ & $	256866619	$ & $	2205518	$ & $	116.465438	$ & $	1.164654376	$ & $	45.34082243	$ & $	38.93071057	$ \\
$	200	$ & $	900	$ & $	8227353208	$ & $	35302466	$ & $	233.053215	$ & $	1.165266076	$ & $	45.32261287	$ & $	38.89464715	$ \\
$	300	$ & $	1350	$ & $	62496428392	$ & $	178753878	$ & $	349.622784	$ & $	1.165409278	$ & $	45.31370223	$ & $	38.88222195	$ \\
$	400	$ & $	1800	$ & $	263399396125	$ & $	564995595	$ & $	466.197256	$ & $	1.165493140	$ & $	45.31008777	$ & $	38.87632299	$ \\
$	500	$ & $	2250	$ & $	803900006590	$ & $	1379446698	$ & $	582.769894	$ & $	1.165539789	$ & $	45.30802103	$ & $	38.87299383	$ \\
$	600	$ & $	2700	$ & $	2000468396580	$ & $	2860504080	$ & $	699.341214	$ & $	1.165568690	$ & $	45.30669993	$ & $	38.87089650	$ \\
$	700	$ & $	3150	$ & $	4323958989350	$ & $	5299548768	$ & $	815.910784	$ & $	1.165586834	$ & $	45.30574404	$ & $	38.86947134	$ \\
$	800	$ & $	3600	$ & $	8430487428682	$ & $	9040920238	$ & $	932.481120	$ & $	1.165601400	$ & $	45.30512262	$ & $	38.86845248	$ \\
$	900	$ & $	4050	$ & $	15192308794063	$ & $	14481975268	$ & $	1049.049492	$ & $	1.165610546	$ & $	45.30459332	$ & $	38.86769338	$ \\
$	1000	$ & $	4500	$ & $	25728695195597	$ & $	22072986511	$ & $	1165.619124	$ & $	1.165619124	$ & $	45.30424551	$ & $	38.86710898	$ \\
\hline
$	1100	$ & $	4950	$ & $	41436812404716	$ & $	32317280628	$ & $	1282.187474	$ & $	1.165624976	$ & $	45.30393559	$ & $	38.86664795	$ \\
$	1200	$ & $	5400	$ & $	64022597756042	$ & $	45771110797	$ & $	1398.755605	$ & $	1.165629670	$ & $	45.30368531	$ & $	38.86627671	$ \\
$	1300	$ & $	5850	$ & $	95531638650235	$ & $	63043701776	$ & $	1515.324068	$ & $	1.165633899	$ & $	45.30349455	$ & $	38.86597207	$ \\
$	1400	$ & $	6300	$ & $	138380047521949	$ & $	84797295248	$ & $	1631.892233	$ & $	1.165637310	$ & $	45.30333177	$ & $	38.86571870	$ \\
$	1500	$ & $	6750	$ & $	195385341804382	$ & $	111747102464	$ & $	1748.460027	$ & $	1.165640018	$ & $	45.30318808	$ & $	38.86550511	$ \\
$	1600	$ & $	7200	$ & $	269797320709960	$ & $	144661271196	$ & $	1865.027996	$ & $	1.165642497	$ & $	45.30307211	$ & $	38.86532295	$ \\
$	1700	$ & $	7650	$ & $	365328940733369	$ & $	184360988174	$ & $	1981.595696	$ & $	1.165644527	$ & $	45.30296828	$ & $	38.86516620	$ \\
$	1800	$ & $	8100	$ & $	486187194335920	$ & $	231720360895	$ & $	2098.163461	$ & $	1.165646367	$ & $	45.30288128	$ & $	38.86503022	$ \\
$	1900	$ & $	8551	$ & $	637103991780086	$ & $	287666521356	$ & $	2214.731102	$ & $	1.165647948	$ & $	45.30280388	$ & $	38.86491110	$ \\
$	2000	$ & $	9000	$ & $	823367026746276	$ & $	353179561833	$ & $	2331.298625	$ & $	1.165649312	$ & $	45.30273472	$ & $	38.86480629	$ \\
\hline\hline
$	3000	$ & $	13501	$ & $		$ & $	1787987440725	$ & $		$ & $		$ & $	45.30233164	$ & $		$ \\
\hline\hline
$	4000	$ & $	18001	$ & $		$ & $	5650944679691	$ & $		$ & $		$ & $	45.30215999	$ & $		$ \\
\hline\hline
$	5000	$ & $	22501	$ & $		$ & $	13796279726716	$ & $		$ & $		$ & $	45.30206783	$ & $		$ \\
\hline\hline
$	6000	$ & $	27001	$ & $		$ & $	28608001167605	$ & $		$ & $		$ & $	45.30201157	$ & $		$ \\
\hline\hline
$	7000	$ & $	31502	$ & $		$ & $	52999897734423	$ & $		$ & $		$ & $	45.30197420	$ & $		$ \\
\hline
\end{tabular}%
}
\caption{Calculated values of $\td$ and associated numbers defined in Section \ref{ss:max}.}
\label{t:rhod}
\end{center}
\end{table}

\section{Conclusion}

Enumeration of integer triangles of given perimeter is a classical problem \cite{JWW1979,Honsberger1985}, and polygons have recently been treated as well \cite{EN2019}.  The current article considered the corresponding problem for tetrahedra, which is considerably more difficult.  The central theme is the calculation of the numbers~$\tn$,~$\td$ and~$\tnd$ of integer tetrahedra with perimeter $n$ and/or diameter $d$ (as appropriate), up to congruence.
In Section \ref{s:BL} we set up a framework in which these numbers could (in principle) be calculated via an application of Burnside's Lemma, and made some partial progress by finding explicit formulas for some of the relevant fix-set parameters.  
The complexity of the other parameters---and of the numbers~$\tn$,~$\td$ and~$\tnd$ themselves---is highlighted by the computations discussed in Section \ref{s:computed}, and is visible in many of the figures therein.  
Nevertheless, a number of apparent patterns emerged from an exploration of the data.  Several of these were discussed in Section \ref{s:obs}, which contains a number of conjectures and avenues for future exploration; we hope that future studies will shed further light on the situation.


\footnotesize
\def\bibspacing{-1.1pt}
\bibliography{biblio}
\bibliographystyle{abbrv}

\end{document}